\documentclass[11pt]{article}     
\usepackage{amsthm,amsfonts,amsmath,amscd,amssymb,latexsym}      
\usepackage{epsfig,graphicx,color}       
\usepackage{color}      
\usepackage[all]{xy}
\SelectTips{cm}{}
\usepackage[titles]{tocloft}    
\usepackage{mathrsfs}
      
\title{\vspace{-2.5cm}  Lefschetz fibrations, intersection numbers,\\     
      and representations of the framed braid group}     

\author{    
Gw\'ena\"el~Massuyeau\\
IRMA, Strasbourg  \and
Alexandru~Oancea\thanks   
{Partially supported by the Swiss National 
Science Foundation, ANR project 
``Floer Power" ANR-08-BLAN-0291-03 and ERC StG-259118-Stein.} \\
IMJ, Paris  \and   
Dietmar~A.~Salamon \\    
ETH Z\"urich    
}

\date{30 August 2013}

\newtheorem{PARA}{}[section] 
\newtheorem{theorem}[PARA]{Theorem} 
\newtheorem{corollary}[PARA]{Corollary}      
\newtheorem{lemma}[PARA]{Lemma}      
\newtheorem{proposition}[PARA]{Proposition}      
\newtheorem{definition}[PARA]{Definition}      
\newtheorem{remark}[PARA]{Remark}      
\newtheorem{example}[PARA]{Example} 
\newcommand{\para}{\begin{PARA}\rm}      
\newcommand{\arap}{\end{PARA}\rm}      
\newcommand{\dfn}{\begin{definition}\rm}      
\newcommand{\nfd}{\end{definition}\rm}      
\newcommand{\rmk}{\begin{remark}\rm}      
\newcommand{\kmr}{\end{remark}\rm}      
\newcommand{\xmpl}{\begin{example}\rm}      
\newcommand{\lpmx}{\end{example}\rm}      
\newcommand{\cA}{\mathcal{A}}      
\newcommand{\cB}{\mathcal{B}}      
\newcommand{\cC}{\mathcal{C}}

\newcommand{\cG}{\mathcal{G}}      
\newcommand{\cH}{\mathcal{H}}      
     
\newcommand{\cJ}{\mathcal{J}}      
      
\newcommand{\cL}{\mathcal{L}}      
\newcommand{\cM}{\mathcal{M}}      
\newcommand{\cN}{\mathcal{N}}      
\newcommand{\cP}{\mathcal{P}}

\newcommand{\cS}{\mathcal{S}}      
\newcommand{\cT}{\mathcal{T}}      
      
\newcommand{\cV}{\mathcal{V}}      
\newcommand{\cW}{\mathcal{W}}

\newcommand{\cZ}{\mathcal{Z}}      

\newcommand{\sN}{\mathscr{N}}

\newcommand{\sT}{\mathscr{T}}

\newcommand{\C}{{\mathbb{C}}}      
\newcommand{\D}{{\mathbb{D}}}      
      
\renewcommand{\H}{{\mathbb{H}}}      
\renewcommand{\L}{{\mathbb{L}}}

\newcommand{\R}{{\mathbb{R}}}      
\newcommand{\bS}{{\mathbb{S}}}      
\newcommand{\T}{{\mathbb{T}}}      
      
\newcommand{\Z}{{\mathbb{Z}}}      
\newcommand{\gS}{{\mathfrak{S}}}      
\newcommand{\tcB}{{\widetilde{\mathcal B}}}      
\newcommand{\tcC}{{\widetilde{\mathcal C}}}      
\newcommand{\tcG}{{\widetilde{\mathcal G}}}

\newcommand{\tn}{{\widetilde n}}      
\newcommand{\ts}{{\widetilde s}}

\newcommand{\tPhi}{{\widetilde\Phi}}

\newcommand{\tchi}{{\widetilde\chi}}

\newcommand{\tal}{{\widetilde\alpha}}      
\newcommand{\tbe}{{\widetilde\beta}}      
\newcommand{\tga}{{\widetilde\gamma}}      
\newcommand{\teps}{{\widetilde\eps}}      
\newcommand{\tcP}{{\widetilde\cP}}
\newcommand{\lk}{{\ell\mathrm{k}}}
\newcommand{\one}      
{{{\mathchoice \mathrm{ 1\mskip-4mu l} \mathrm{ 1\mskip-4mu l}      
\mathrm{ 1\mskip-4.5mu l} \mathrm{ 1\mskip-5mu l}}}}      
\newcommand{\tr}{\mathrm{ tr }}        

\newcommand{\G}{{\mathrm{G}}}       
\newcommand{\GL}{{\mathrm{GL}}}      
\newcommand{\id}{\mathrm{ id}}         
\newcommand{\Id}{\mathrm{ Id}}      
       
\newcommand{\PSL}{{\mathrm{PSL}}}               
\newcommand{\SO}{{\mathrm{SO}}}               
               
\newcommand{\Aut}{\mathrm{ Aut}}          
\newcommand{\Diff}{\mathrm{ Diff}}        
\newcommand{\Vect}{\mathrm{ Vect}}        
\newcommand{\eps}{{\varepsilon}}

\newcommand{\inner}[2]{\left\langle #1, #2\right\rangle}      

\def\NABLA#1{{\mathop{\nabla\kern-.5ex\lower1ex\hbox{$#1$}}}}      
\def\Nabla#1{\nabla\kern-.5ex{}_{#1}}      
\def\Tabla#1{\Tilde\nabla\kern-.5ex{}_{#1}}      
\def\abs#1{\mathopen|#1\mathclose|}      
\def\Abs#1{\left|#1\right|}

\renewcommand{\Tilde}{\widetilde}

\newcommand{\p}{{\partial}}

\begin{document}      
 
\maketitle      

      
\vspace{-1cm}    
 
\begin{abstract}     
We examine the action of the fundamental group $\Gamma$ 
of a Riemann surface with $m$ punctures on the middle dimensional     
homology of a regular fiber in a Lefschetz fibration,   
and describe to what extent this action can be recovered    
from the intersection numbers of vanishing cycles.   
Basis changes for the vanishing cycles result in    
a nonlinear action of the framed braid group $\tcB$ on $m$ strings  
on a suitable space of $m\times m$ matrices.  
This action is determined by a family of cohomologous $1$-cocycles 
$\cS_c:\tcB\to\GL_m(\Z[\Gamma])$ parametrized by  
distinguished configurations $c$ of embedded  
paths from the regular value to the critical values. 
In the case of the disc, we compare this family of cocycles with the Magnus 
cocycles given by Fox calculus and consider some abelian reductions 
giving rise to linear representations of braid groups. 
We also prove that, still in the case of the disc,  
the intersection numbers along straight lines, 
which conjecturally make sense in infinite dimensional situations, 
carry all the relevant information.
\end{abstract}      

\vspace{-0.5cm}    

\setlength{\cftbeforesecskip}{0.4ex}    
{\small \tableofcontents}    
      
      
\section{Introduction} \label{sec:intro}    

Picard--Lefschetz theory can be viewed as a complexification    
of Morse theory with the stable and unstable manifolds replaced    
by vanishing cycles and the count of connecting orbits in the     
Morse--Witten complex replaced by the intersection numbers    
of the vanishing cycles along suitable paths in the base.     
Relevant topological information that can be recovered from these 
data includes, in Morse theory, the homology of the underlying
manifold and, in Picard--Lefschetz theory, the monodromy action of 
the fundamental group of the base on the middle dimensional homology
of a regular fiber.  That the monodromy action on the vanishing cycles
can be recovered from the intersection numbers follows from the     
{\bf Picard--Lefschetz formula}    
\begin{equation}\label{eq:monodromy}    
(\psi_L)_*\alpha = \alpha - (-1)^{n(n+1)/2}\inner{L}{\alpha}L.    
\end{equation}    
Here $X\to\Sigma$ is a Lefschetz fibration over a Riemann surface
$\Sigma$ with fibers of complex dimension $n$, meaning that $X$ is a complex manifold and the map $X \to \Sigma$ is holomorphic and has only nondegenerate critical points.
We assume moreover that each singular fiber contains exactly one critical point,
and denote by $M$   a regular fiber over $z_0\in\Sigma$.
In equation \eqref{eq:monodromy},
$L\subset M$ is an oriented vanishing cycle associated to a curve     
from $z_0$ to a singular point $z\in\Sigma$, $\psi_L\in\Diff(M)$  is the Dehn--Arnold--Seidel twist about $L$ 
obtained from the (counterclockwise) monodromy around the singular fiber    
along the same curve, $(\psi_L)_*$ is the induced action on $H_n(M)$    
and $\inner{\cdot}{\cdot}$ denotes the intersection form.     
Equation~\eqref{eq:monodromy} continues to hold when    
$X\to\Sigma$ is a symplectic Lefschetz fibration as introduced     
by Donaldson~\cite{D2,D3}. In either case the vanishing     
cycles are embedded Lagrangian spheres and so their    
self-intersection numbers are $2(-1)^{n/2}$ when $n$ is even    
and $0$ when $n$ is odd. See~\cite[Chapter~I]{AGV}
and~\cite[\S2.1]{AGLV} for a detailed account of Picard--Lefschetz
theory and an exhaustive reference list.  

The object of the present paper is to study an algebraic setting,    
based on equation~\eqref{eq:monodromy}, which allows one to describe    
the monodromy action of the fundamental group in terms of     
intersection matrices. This requires the choice of a distinguished     
basis of vanishing cycles and the ambiguity in this choice leads    
to an action of the braid group on distinguished basis~\cite{AGV},
which in turn determines an action on a suitable space of
matrices. In the case of the disc, this action was previously
considered by Bondal~\cite{Bondal} in the context of mirror
symmetry. Our motivation is different and comes from an attempt 
of two of the authors (A.O. and D.S.) 
to understand complexified Floer homology in the spirit of
Donaldson--Thomas theory~\cite{DT}.  In this theory 
the complex symplectic action or Chern--Simons functional     
is an infinite dimensional analogue of a Lefschetz fibration.    
While there are no vanishing cycles, one can (conjecturally)    
still make sense of their intersection numbers along straight lines
and build an intersection matrix whose orbit under the braid group     
might then be viewed as an invariant. Another source of inspiration    
for the present paper is the work of Seidel~\cite{S1,S2} about     
vanishing cycles and mutations.     

\begin{figure}[htp]      
\centering       
\includegraphics[scale=0.6]{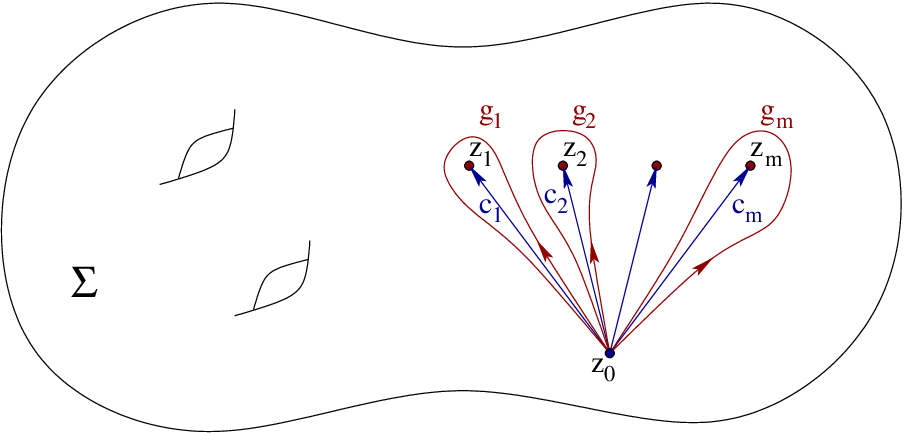}      
\caption{{A distinguished configuration.}}\label{fig:1}      
\end{figure}

We assume throughout this paper that the base $\Sigma$ 
of our Lefschetz fibration is a compact Riemann surface, 
possibly with boundary, not diffeomorphic to the $2$-sphere.  
(In particular, $\Sigma$ is oriented.)
Let $Z\subset\Sigma\setminus\p\Sigma$ be the set of critical 
values and $z_0$ be a regular value.  If $\Sigma$ is diffeomorphic 
to the unit disc $\D=\{z\in\C||z|\le 1\}$ we assume that $z_0\in\p\Sigma$. 
To assemble the intersection numbers into algebraic data it is
convenient to choose a collection $c=(c_1,\dots,c_m)$  of ordered
embedded paths from a regular value $z_0\in\Sigma$ to the critical
values (see Figure~\ref{fig:1}). Following~\cite{AGV,S1} we call such
a collection a {\bf distinguished configuration} and denote 
by $\cC$ the set of homotopy classes of distinguished configurations. 
Any distinguished configuration determines an    
ordering $\{z_1,\dots,z_m\}$ of the set $Z$ of critical     
values by $z_i:=c_i(1)$. It also determines $m$ vanishing cycles     
$L_1,\dots,L_m\subset M$ as well as $m$ special elements    
of the fundamental group    
$$    
g_1,\dots,g_m\in\Gamma := \pi_1(\Sigma\setminus Z,z_0)   
$$    
(obtained by encircling $z_i$ counterclockwise along $c_i$).   
The orientation of $L_i$ is not determined by the path $c_i$ and    
can be chosen independently. However, when $n$ is even,    
the monodromy along $g_i$ changes the orientation of~$L_i$.    
Thus, to choose orientations consistently, we fix nonzero    
tangent vectors $v_z\in T_z\Sigma$ for $z\in Z$,    
choose orientations of the vanishing cycles    
in the directions of these vectors, and consider only    
distinguished configurations $c$ that are tangent to the vectors    
$-v_z$ at their endpoints. We call these {\bf marked distinguished    
configurations} and denote by $\tcC$ the set of homotopy classes     
of marked distinguished configurations. The oriented vanishing cycles
determine homology classes, still denoted by    
$$   
L_1,\dots,L_m\in H_n(M).   
$$   
These data give rise to a {\bf monodromy character} $\cN^X_c:\Gamma\to
\Z^{m\times m}$ via 
\begin{equation}\label{eq:NXc}      
\cN^X_c(g):=\bigl(n_{ij}(g)\bigr)_{i,j=1}^m,\qquad    
n_{ij}(g) := \inner{L_i}{\rho(g)L_j}.
\end{equation}
Here $\rho:\Gamma\to \Aut(H_n(M))$ denotes the monodromy action of the
fundamental group. Any such function $\cN:\Gamma\to \Z^{m\times m}$
satisfies the conditions 
\begin{eqnarray}
n_{ij}(g^{-1}) & = & (-1)^n n_{ji}(g), \nonumber \\
n_{ii}(1) & = & \left\{\begin{array}{ll}      
0,&\mbox{if }n\mbox{ is odd},\\      
2(-1)^{n/2},&\mbox{if }n\mbox{ is even},
\end{array}\right. \label{eq:nij} \\
n_{ij}(gg_kh)&=&n_{ij}(gh)-(-1)^{n(n+1)/2} n_{ik}(g)n_{kj}(h)  \nonumber
\end{eqnarray}    
for $g,h\in\Gamma$ and $i,j,k=1,\dots,m$. The last equation
in~\eqref{eq:nij} follows from~\eqref{eq:monodromy}. 
Our convention for the composition 
is that $gh\in\Gamma$ first traverses
a representative of $h$, and then a representative of $g$. 

The part of $H_n(M)$ that is generated by the vanishing cycles 
under the action of $\Gamma$ can be recovered as the quotient 
$$
\H_\cN:=\Lambda/\ker \ \cN, \qquad \Lambda:=\Z[\Gamma]^m.
$$
Here $\Z[\Gamma]$ is the group ring of $\Gamma$, whose elements are thought of as maps $\lambda:\Gamma \to \Z$ with finite support
and whose multiplication is the convolution product $(\lambda_1 \lambda_2)(h):=\sum_g \lambda_1(hg^{-1})\lambda_2(g)$;
the map $\cN:\Gamma\to \Z^{m\times m}$ is regarded as an endomorphism $\cN:\Lambda\to \Lambda$ by the
convolution product $\lambda\mapsto \big(h \mapsto \sum_g\cN(hg^{-1})\lambda(g)\big)$. 
The $\Z$-module $\H_\cN$ is equipped with an intersection pairing, 
a $\Gamma$-action $\rho_\cN:\Gamma\to\Aut(\H_\cN)$, and  
special elements $\L_1,\dots,\L_m\in\H_\cN$, defined by  
\begin{equation} \label{eq:HN1}  
\inner{\mu}{\lambda}  
:= \sum_{g,h\in\Gamma}\mu(h)^T \cN(hg^{-1})\lambda(g), \quad 
\rho_\cN(g)[\lambda] := [g_*\lambda],\quad 
\L_i := [\delta_i]. 
\end{equation} 
Here $\Gamma$ acts on $\Lambda$  
by 
$$
(g_*\lambda)(h):=\lambda(hg),
$$
$e_i\in\Z^m$ denotes the standard basis, and  
$\delta_i:\Gamma\to\Z^m$ is defined by 
$\delta_i(1) := e_i$ and $\delta_i(g):=0$ for $g\ne 1$. 
With these structures $\H_\cN$ is isomorphic to the  
submodule of $H_n(M)$ generated by the vanishing cycles modulo  
the kernel of the intersection form. The isomorphism is induced  
by the map that assigns to every $\lambda\in\Lambda$
the homology class $\sum_{g,i}\lambda_i(g)\rho(g^{-1})L_i\in H_n(M)$.

The monodromy character $\cN^X_c:\Gamma\to\Z^{m\times m}$ depends  
on the choice of a marked distinguished configuration $c$ 
and this dependence gives rise to an action of the framed 
braid group $\tcB$ of $\Sigma\setminus\{z_0\}$ on $m$ strings, 
based at $Z$, on the space of monodromy characters. 
More precisely, our distinguished configuration $c$ 
determines $m$ special elements $g_i=:g_{i,c}\in\Gamma$ 
obtained by encircling ${z_i:=c_i(1)}$ counterclockwise 
along $c_i$. Denote by  
$$ 
\sN_c := \left\{\cN=(n_{ij}):\Gamma\to\Z^{m\times m}\,|\,
\eqref{eq:nij}\right\} 
$$ 
the space of monodromy characters on $(\Gamma,g_{1,c},\dots,g_{m,c})$.  
The framed braid group~$\tcB$, interpreted as the mapping class 
group of diffeomorphisms in $\Diff_0(\Sigma,z_0)$ 
that preserve the set $Z$ and the collection of vectors $\{v_z\}_{z\in
Z}$, acts freely and transitively on the space $\tcC$ of homotopy
classes of marked distinguished configurations (see Sections~\ref{sec:braid}
and~\ref{sec:bases}). Here $\Diff_0(\Sigma,z_0)$ denotes the identity
component of the group of diffeomorphisms of $\Sigma$ that fix
$z_0$. The framed braid group also acts on the fundamental group
$\Gamma$ and, for every $\tau\in\tcB$ and every $c\in\tcC$, the
isomorphism $\tau_*:\Gamma\to\Gamma$ maps $g_{i,c}$ to $g_{i,\tau_*c}$.
This action actually determines the (unframed) braid $\tau$
(see Section~\ref{sec:bases}).

Our first theorem asserts that there is a canonical family of
isomorphisms   
$$
\sT_{\tau,c}:\sN_c\to\sN_{\tau_*c}
$$
which extends the geometric correspondence  
${\cN^X_c\mapsto\cN^X_{\tau_*c}}$ between mono\-dromy  
characters associated to different choices of the distinguished  
configuration. It is an open question if every element $\cN\in\sN_c$
can be realized by a (symplectic) Lefschetz fibration  
$X\to\Sigma$ with critical fibers over $Z$. 
This question is a refinement of the Hurwitz problem of 
finding a branched cover with given combinatorial data.
The isomorphisms $\sT_{\tau,c}$ are determined by a family 
of cocycles 
$$
\cS_c:\tcB\to \GL_m(\Z[\Gamma])
$$
with values in the group of invertible
matrices over the group ring. To describe them we denote by  
$\pi_{\sigma,c}\in\gS_m$ the permutation associated to the 
action of $\sigma$ on the ordering determined by~$c$.  
Then the $(i,j)$ entry of the matrix is 
$$
s_{j,c}(\sigma):=c_i^{-1}\cdot \sigma_*c_j\in\Gamma
$$
(first $\sigma_*c_j$, second $c_i^{-1}$)
for $i=\pi_{\sigma,c}(j)$ and it is zero for
$i\ne\pi_{\sigma,c}(j)$. We emphasize that  
$
\sigma_*g_{j,c} = s_{j,c}(\sigma)^{-1}g_{i,c}s_{j,c}(\sigma)
$
for $i:=\pi_{\sigma,c}(j)$.
See Figure~\ref{fig:SiEpsCoc} below for some examples.

In the next theorem we think of an element
$\cM\in\GL_m(\Z[\Gamma])$ as a function  
$\cM:\Gamma\to\Z^{m\times m}$ with finite support, denote 
$\cM^t(g):=\cM(g^{-1})^T$, and multiply matrices
using the convolution product (see Section~\ref{sec:monodromy}).

\smallskip
\noindent {\bf Theorem~A.} 
{\bf (i)} 
{\it The maps $\cS_c:\tcB\to\GL_m(\Z[\Gamma])$
are injective and satisfy the cocycle and coboundary 
conditions}
\begin{equation}\label{eq:cocb} 
\cS_c(\sigma\tau) = \cS_c(\sigma)\sigma_*\cS_c(\tau),\qquad
\cS_{\tau_*c}(\sigma) = \cS_c(\tau)^{-1}\cS_c(\sigma)\sigma_*\cS_c(\tau). 
\end{equation} 

\noindent{\bf (ii)} 
{\it The maps $\cS_c$ in~(i) determine bijections}
$$
\sT_{\tau,c}:\sN_c\to \sN_{\tau_*c}, \qquad 
\cN\mapsto\cS_c(\tau)^t\cN\cS_c(\tau) .
$$

\noindent{\bf (iii)} 
{\it Given a Lefschetz fibration $X\to \Sigma$ and elements
$c\in\tcC$, $\tau \in\tcB$, we have}
\begin{equation}\label{eq:Tpsic} 
\cN^X_{\tau_*c} = \cS_c(\tau)^t\cN^X_c\cS_c(\tau).
\end{equation} 

\smallskip

In terms of Serre's definition of non-abelian    
cohomology~\cite[Appendix to Chapter VII]{Serre}   
the first equation in~\eqref{eq:cocb} asserts that~$\cS_c$    
is a cocycle for every $c\in\tcC$. The second equation in~\eqref{eq:cocb}    
asserts that the cocycles $\cS_c$ are all cohomologous    
and hence define a canonical $1$-cohomology class       
$$
[\cS_c]\in H^1(\tcB,\GL_m(\Z[\Gamma])).    
$$
We call it the {\bf Picard--Lefschetz monodromy class}. 

In the case $\Sigma=\D$ there is another well known cocycle 
arising from a topological context, namely the Magnus cocycle 
$$
\cM_c:\cB\to\GL_m(\Z[\Gamma]),\qquad c\in\cC.
$$
Here $\cB$ denotes the usual braid group (with no framing), 
which we view as a subgroup of $\tcB$ 
using the framing determined by 
a vector field~$v$ on~$\Sigma$ whose only singularity 
is an attractive point at $z_0$ and such that 
${v(z)=v_z}$ for all $z\in Z$ (see Section~\ref{sec:magnus}).
The Magnus cocycle is also related to an intersection 
pairing~\cite{Turaev,Perron} and its dependence on the choice 
of the distinguished configuration $c\in\cC$ is similar 
to the dependence of the Picard--Lefschetz cocycle. 
The cohomology classes $[\cM_c]$ and $[\cS_c|_\cB]$ are distinct 
and nontrivial in $H^1(\cB,\GL_m(\Z[\Gamma]))$. 
After reduction of $\Gamma$ to the infinite cyclic group, 
both cocycles define linear representations of the braid group 
with coefficients in $\Z[t,t^{-1}]$. In the case of the Magnus cocycle, 
this is the famous Burau representation~\cite{BIRMAN,KT}. 
In the case of the Picard--Lefschetz cocycle, this representation 
was first discovered by Tong--Yang--Ma~\cite{TYM} and is a 
key ingredient in the classification of $m$-dimensional 
representations of the braid group $\cB$ on $m$ 
strings~\cite{Sysoeva}. 
For pure braids, i.e$.$ braids which do not permute the elements of $Z$, the Tong--Yang--Ma representation 
is determined by linking numbers (see Section~\ref{sec:magnus}).

\smallbreak

We continue our discussion of the planar case $\Sigma=\D$.
The group $\Gamma$ is then isomorphic to the free group 
$\Gamma_m$ generated by $g_1,\dots,g_m$, and it is convenient 
to switch from the geometric picture in Theorem~A to generators 
and relations. We denote by $\tcB_m$ the abstract group 
generated by $\sigma_2,\dots,\sigma_m$, $\eps_1,\dots,\eps_m$ 
with relations    
\begin{equation}\label{eq:frabraid}   
\sigma_i\sigma_{i+1}\sigma_i=\sigma_{i+1}\sigma_i\sigma_{i+1},\qquad   
\eps_i\sigma_i=\sigma_i\eps_{i-1},\qquad   
\eps_{i-1}\sigma_i=\sigma_i\eps_i.   
\end{equation}   
All other pairs of generators commute. The choice of an element
$c\in\tcC$ determines an isomorphism $\tcB_m\to \tcB$ 
obtained by identifying the generators $\sigma_i$, $\eps_i$ 
with the moves depicted in Figure~\ref{fig:SiEpsCoc} 
(\cite{KS}, see also Section~\ref{sec:bases}).
\begin{figure}[htp]
\centering       
\input{figure-se-coc.pstex_t}      
\caption{{The generators of $\tcB_m$.}}\label{fig:SiEpsCoc}        
\end{figure} 
This gives rise to a contravariant free
and transitive action of $\tcB_m$ on $\tcC$ denoted by 
$\tcB_m\times\tcC\to\tcC:(\sigma,c)\mapsto\sigma^*c$, and to an action
of $\tcB_m$ on $\Gamma_m$ via 
\begin{equation}\label{eq:BmGm}   
(\sigma_i)_*:   
\left\{\begin{array}{ccl}     
  g_{i-1} &\mapsto &g_{i-1} g_i g_{i-1}^{-1}, \\    
  g_i & \mapsto &g_{i-1},    
\end{array}\right.    
\end{equation}
$(\sigma_i)_*g_j=g_j$ for $j\ne i-1,i$, and $(\eps_i)_*=\Id$.

The third equation in~\eqref{eq:nij} shows that, in the case of the
disc, we can switch from matrix valued functions to actual matrices. 
More precisely, a monodromy character
$\cN:\Gamma_m\to\Z^{m\times m}$ is uniquely determined by the single
matrix $N:=\cN(1)$. The latter satisfies
\begin{equation}\label{eq:nij-N} 
n_{ij}=(-1)^nn_{ji},  \qquad  
n_{ii}= \left\{\begin{array}{ll}     
0,&\mbox{if }n\mbox{ is odd},\\     
2(-1)^{n/2},&\mbox{if }n\mbox{ is even}.
\end{array}\right.  
\end{equation} 
We denote by $\sN_m$ the space of matrices satisfying~\eqref{eq:nij-N}. 
The map $\cN$ is explicitly given by $\cN(g)=N\rho_N(g)$,
the homomorphism $\rho_N:\Gamma_m\to\GL_m(\Z)$ 
being defined on generators by
\begin{equation}\label{eq:rhoN}   
\rho_N(g_i)=\one-(-1)^{n(n+1)/2}E_iN.   
\end{equation}   
Here $E_i\in\Z^{m\times m}$ is the matrix with the    
$i$-th entry on the diagonal equal to one and zeroes elsewhere. 
The representation $\rho_N$ induces an action of $\Gamma_m$ 
on the quotient module $\H_N:=\Z^m/\ker N$ which preserves 
the intersection form $\inner{\lambda}{\mu}:=\lambda^TN\mu$.     
Moreover, the triple $(\H_N,\rho_N,\inner{\cdot}{\cdot})$ becomes
isomorphic to $(\H_\cN,\rho_\cN,\inner{\cdot}{\cdot})$:
see Section~\ref{sec:monodromy}.   The next result rephrases 
Theorem~A for the particular case of the disc, 
and strengthens it with an additional uniqueness statement. 
The first part of~(ii) has been already proved by
Bondal~\cite[Proposition~2.1]{Bondal} for upper triangular
matrices instead of (skew-)symmetric ones. 
For $k=2,\dots,m$ we denote by $\Sigma_k$ the permutation    
matrix of the transposition $(k-1,k)$    
and, for $i=1,\dots,m$, we denote by $D_i$ the diagonal matrix    
with $i$-th entry on the diagonal equal to $(-1)^{n+1}$    
and the other diagonal entries equal to~$1$.
 
\medskip\noindent{\bf Theorem~B.}   
{\bf (i)}  {\it There is a unique function
$S:\tcB_m\times \sN_m \to \GL_m(\Z)$   
satisfying the following conditions.}   
   
\smallskip\noindent{\bf (Cocycle)}    
{\it For all $N\in\sN_m$ and $\sigma,\tau\in\tcB_m$ we have}   
\begin{equation} \label{eq:Cocycle}    
S(\sigma\tau,N)    
= S(\sigma,N)S(\tau, S(\sigma,N)^T N S(\sigma,N)).   
\end{equation}    
   
\smallskip\noindent{\bf (Normalization)}    
{\it For all $N\in\sN_m$ we have $S(1,N)=\one$ and}   
\begin{equation} \label{eq:normalize}    
\begin{split}   
S(\sigma_k,N) &=\Sigma_k - (-1)^{n(n+1)/2}n_{k-1,k}E_{k-1},\qquad   
k=2,\dots,m,  \\   
S(\eps_i,N) &= D_i,\qquad   
i=1,\dots,m.   
\end{split}   
\end{equation}  

\smallskip\noindent{\bf (ii)}    
{\it The function $S$ in~(i) determines a contravariant 
group action of $\tcB_m$ on $\sN_m$ via
$$
\sigma^*N:=S(\sigma,N)^T N S(\sigma,N)
$$
for $N\in\sN_m$ and $\sigma\in\tcB_m$.  This action 
is compatible with the action of $\tcB_m$ on the 
space of marked distinguished configurations 
in the sense that, for every Lefschetz fibration $X\to\D$ 
with singular fibers over $Z$, every ${\sigma\in\tcB_m}$, 
and every $c\in\tcC$, we have
$$
N^X_{\sigma^*c} = \sigma^*N^X_c,
$$
where $N^X_c := \cN^X_c(1)$.} 

\bigbreak
 
\smallskip\noindent{\bf (iii)} 
{\it For every $\sigma\in \tcB_m$ and every $N\in\sN_m$    
the matrix $S(\sigma,N)$ induces an isomorphism from    
$\H_{\sigma^*N}$ to $\H_N$ which preserves the bilinear pairings    
and satisfies}   
\begin{equation} \label{eq:Rho}    
\rho_{\sigma^*N}(g)     
= S(\sigma,N)^{-1}\rho_N(\sigma_*g) S(\sigma,N).   
\end{equation}    
   
\medskip   
   
By Theorem~B every symplectic Lefschetz fibration $f:X\to\D$    
with critical fibers over $Z$ determines a $\tcB_m$-equivariant map   
$$   
\tcC\to\sN_m:c\mapsto N^X_c   
$$    
which can be viewed as an algebraic invariant of $X$.    
Our next theorem asserts that this invariant is uniquely determined    
by the matrix   
$$   
Q^X:Z\times Z\to\Z   
$$   
of intersection numbers of vanishing cycles along straight lines.    
Here we assume that no straight line connecting two points   
in $Z$ contains another element of $Z$; such a set $Z$   
is called {\bf admissible}.   Denote by $\sN_Z$
the space of matrices $Q:Z\times Z\to\Z$ satisfying   
\begin{equation}\label{eq:Q}   
Q(z,z') = (-1)^nQ(z',z),\qquad   
Q(z,z) = \left\{\begin{array}{ll}      
0,&\mbox{if }n\mbox{ is odd},\\      
2(-1)^{n/2},&\mbox{if }n\mbox{ is even}.   
\end{array}\right.   
\end{equation}   
We define a map   
$   
\tcC\times\sN_m\to\sN_Z:(c,N)\mapsto Q_{c,N}   
$   
by   
\begin{equation}\label{eq:QNc}   
Q_{c,N}(z_i,z_j) := (N\rho_N(g))_{ij}   
\end{equation}   
where $z_1,\dots,z_m$ is the ordering of $Z$ given by $c$,   
$\Gamma_m$ is identified with~$\Gamma$, and   
$$   
g := c_i^{-1}\cdot s_{ij}\ \cdot c_j\in\Gamma.   
$$   
Here the right hand side denotes the based loop obtained by   
first traversing $c_j$, then moving clockwise near $z_j$   
until reaching the straight line $s_{ij}$ from $z_j$ to $z_i$,   
then following $s_{ij}$, then moving counterclockwise   
near $z_i$ until reaching $c_i$, and finally traversing $c_i$    
in the opposite direction    
\begin{figure}[htp]      
\centering       
\includegraphics[scale=0.6]{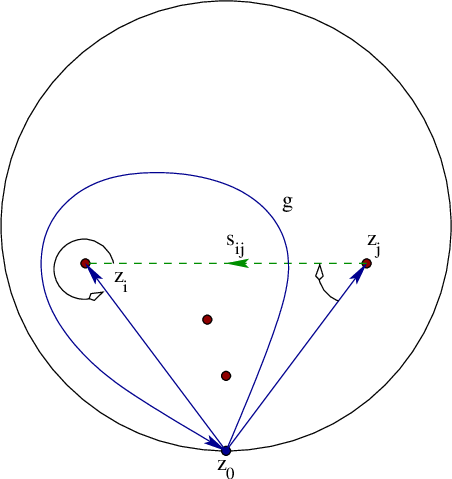}      
\caption{{Intersection numbers along straight lines.}}
\label{fig:straight} 
\end{figure} 
(see Figure~\ref{fig:straight}). Geometrically, this means 
that the intersection matrix $Q^X:=Q_{c,N^X_c}$ assigns to a pair    
$(z,z')\in Z\times Z$ the intersection number of the vanishing cycles    
along the straight line from $z$ to $z'$, where the orientations at the    
endpoints are determined by moving the oriented vanishing    
cycles in the directions $v_z$ and $v_{z'}$ clockwise towards    
the straight line. (Given a marked distinguished configuration   
$c$ and the loop $g$ as above, the straight line $s_{ji}$ from    
$z_i$ to $z_j$ corresponds to the curve $c_j\cdot g^{-1}\cdot c_i^{-1} $.)   
   
\medskip\noindent{\bf Theorem~C.}   
{\it The map $(c,N)\mapsto Q_{c,N}$ defined by~\eqref{eq:QNc}    
is invariant under the diagonal action of $\tcB_m$ on    
$\tcC\times\sN_m$.  Moreover, for every $Q\in\sN_Z$,   
there is a unique equivariant map   
$$
\tcC\to\sN_m:c\mapsto N_c
$$ 
such that    
$$
Q_{c,N_c} = Q   
$$
for every $c\in\tcC$.}   
   
\medskip   
Let $f:X\to\D$ be a symplectic Lefschetz fibration with critical fibers    
over an admissible set $Z\subset{\mathrm{int}}(\D)$.   
Let $z,z'\in Z$ and $x\in f^{-1}(z)$, $x'\in f^{-1}(z')$ be the   
associated critical points of $f$.  Then the number    
$Q^X(z,z')$ is the algebraic count of negative gradient flow lines   
\begin{equation}\label{eq:grad}   
\dot u + \nabla f^\theta(u) = 0,\quad   
f^\theta := \cos(\theta)\mathrm{Re}\,f + \sin(\theta)\mathrm{Im}\,f,\quad   
\theta := \arg(z'-z),   
\end{equation}   
from $x'=\lim_{s\to-\infty}u(s)$ to $x=\lim_{s\to\infty}u(s)$.    
According to Donaldson--Thomas~\cite{DT} this count of gradient   
flow lines is (conjecturally) still meaningful in suitable infinite
dimensional settings.  It thus gives rise to an intersection matrix
$Q$ and hence, by Theorem~C, also to an equivariant map $c\mapsto N_c$.    
A case in point, analogous to symplectic Floer theory, is where $f$    
is the complex symplectic action on the loop space of a complex    
symplectic manifold.    
   
The paper is organized as follows.     
In Section~\ref{sec:monodromy} we explain an algebraic setting   
for monodromy representations, Section~\ref{sec:braid} discusses the   
framed braid group $\tcB$, and in Section~\ref{sec:bases} we prove that   
$\tcB$ acts freely and transitively on the space $\tcC$ of distinguished   
configurations. Theorem~A is contained in Theorem~\ref{thm:Sc} from
Section~\ref{sec:coc}. We compare in Section~\ref{sec:magnus} the 
Picard--Lefschetz cocycle with the Magnus cocycle, and discuss some 
related linear representations of the braid group. 
Theorem~B is proved in Section~\ref{sec:proofB}, 
in Section~\ref{sec:groupoids} we introduce monodromy groupoids, 
in Section~\ref{sec:proofC} we prove Theorem~C, 
and Section~\ref{sec:ex} illustrates the monodromy representation 
by an example.   We include a brief discussion of some basic 
properties of Lefschetz fibrations in Appendix~\ref{app:L}, 
summarizing relevant facts from~\cite[Chapter~I]{AGV}. 

 
\medskip     
\noindent{\bf Acknowledgements.}     
We are grateful to Theo B\"uhler for pointing out to us Serre's    
definition of a non-abelian cocycle, and to Ivan Marin 
for indicating to us the references~\cite{Tits,TYM}.  
    
    
\section{Monodromy representations}\label{sec:monodromy}      
     
We examine algebraic structures that are relevant in the study      
of monodromy representations associated to Lefschetz      
fibrations.       
    
\begin{definition}\label{def:Gamma}     
Fix a positive integer $n$.    
Let $\Gamma$ be a group and $g_1,\dots,g_m$ be pairwise     
distinct elements of $\Gamma\setminus\{1\}$.      
A {\bf monodromy character} on $(\Gamma,g_1,\dots,g_m)$     
is a matrix valued function      
$     
\cN=(n_{ij}):\Gamma\to\Z^{m\times m}     
$     
satisfying      
\begin{equation}\label{eq:N1}      
n_{ij}(g^{-1})=(-1)^nn_{ji}(g),      
\end{equation}      
\begin{equation}\label{eq:N2}      
n_{ii}(1) = \left\{\begin{array}{ll}      
0,&\mbox{if }n\mbox{ is odd},\\      
2(-1)^{n/2},&\mbox{if }n\mbox{ is even},      
\end{array}\right.      
\end{equation}      
\begin{equation}\label{eq:N3}      
n_{ij}(gg_kh) = n_{ij}(gh)       
- (-1)^{n(n+1)/2}n_{ik}(g)n_{kj}(h)      
\end{equation}      
for all $g,h\in\Gamma$ and $i,j,k\in\{1,\dots,m\}$.      
A {\bf monodromy representation} of $(\Gamma,g_1,\dots,g_m)$     
is a tuple $(\cH,\rho,\cL_1,\dots,\cL_m)$ consisting of a $\Z$-module     
$\cH$ with a nondegenerate bilinear pairing,     
a representation $\rho:\Gamma\to\Aut(\cH)$ that preserves the      
bilinear pairing, and elements $\cL_1,\dots,\cL_m\in\cH$,     
satisfying  
\begin{equation}\label{eq:GH1}      
\inner{\cL}{\cL'} = (-1)^n\inner{\cL'}{\cL},     
\end{equation}      
\begin{equation}\label{eq:GH2}      
\inner{\cL_i}{\cL_i}       
= \left\{\begin{array}{ll}      
0,&\mbox{if }n\mbox{ is odd},\\      
2(-1)^{n/2},&\mbox{if }n\mbox{ is even},     
\end{array}\right.      
\end{equation}      
\begin{equation}\label{eq:GH3}      
\rho(g_i)\cL = \cL -       
(-1)^{n(n+1)/2}\inner{\cL_i}{\cL}\cL_i     
\end{equation}      
for all $\cL,\cL'\in\cH$ and $i\in\{1,\dots,m\}$. The 
automorphisms
$\rho(g_i)$ are called {\bf Dehn twists} and the $\cL_i$
are called
{\bf vanishing cycles}.     
\end{definition}     
    
\begin{remark}\label{rmk:N-inverse}\rm     
Every monodromy character $\cN:\Gamma\to\Z^{m\times m}$      
satisfies     
\begin{equation}\label{eq:N4}      
n_{ij}(gg_k^{-1}h) = n_{ij}(gh)       
- (-1)^n(-1)^{n(n+1)/2}n_{ik}(g)n_{kj}(h),      
\end{equation}      
\begin{equation}\label{eq:N5}      
n_{ij}(g_ig)=n_{ij}(gg_j)=(-1)^{n+1}n_{ij}(g).      
\end{equation}      
Every monodromy representation $(\cH,\rho,\cL_1,\dots,\cL_m)$     
satisfies     
\begin{equation}\label{eq:GH4}      
\rho(g_i^{-1})\cL = \cL -       
(-1)^n(-1)^{n(n+1)/2}\inner{\cL_i}{\cL}\cL_i,     
\end{equation}      
\begin{equation}\label{eq:GH5}      
\rho(g_i)\cL_i = (-1)^{n+1}\cL_i.     
\end{equation}      
Equation~(\ref{eq:N4}) follows from~(\ref{eq:N1}) and~(\ref{eq:N3}) 
by replacing $g,h$ with $h^{-1},g^{-1}$ and interchanging $i$ and $j$.
To prove the second equation in~(\ref{eq:N5}) use 
equation~(\ref{eq:N3}) with $k=j$ and $h=1$; then use~(\ref{eq:N2}). 
The proofs of~(\ref{eq:GH4}) and~(\ref{eq:GH5}) are similar.     
\end{remark}     
    
Every monodromy representation gives     
rise to a monodromy character and vice versa.    
If $(\cH,\rho,\cL_1,\dots,\cL_m)$ is a monodromy representation    
then the associated monodromy character      
$\cN_\rho:\Gamma\to\Z^{m\times m}$ assigns to every $g\in\Gamma$     
the intersection matrix     
\begin{equation}\label{eq:N}     
n_{ij}(g) := \inner{\cL_i}{\rho(g)\cL_j}.     
\end{equation}    
Conversely, every monodromy character $\cN$ gives rise to a monodromy     
representation $(\H_\cN,\rho_\cN,\L_1,\dots,\L_m)$ as follows.    
   
Denote by $\Z[\Gamma]$ the group ring of $\Gamma$.    
One can think of an element $\lambda\in\Z[\Gamma]$     
either as a function $\lambda:\Gamma\to\Z$ with finite     
support or as a formal linear combination    
$\lambda=\sum_{g\in\Gamma}\lambda(g)g$.   
With the first viewpoint, the multiplication in $\Z[\Gamma]$ 
is the convolution product    
$$    
(\lambda\mu)(h):=\sum_{g\in\Gamma}\lambda(g)\mu(g^{-1}h).    
$$    
The group $\Gamma$ acts on the group ring by the formula    
$(g_*\lambda)(h):= \lambda(hg)$ for $\lambda\in\Z[\Gamma]$    
and $g,h\in\Gamma$.  In the formal sum notation we have    
$$    
g_*\lambda = \sum_h\lambda(hg)h=\sum_h    
\lambda(h)hg^{-1} = \lambda g^{-1}.    
$$     
For any function $\cN:\Gamma\to\Z^{m\times m}$ we introduce     
the $\Z$-module    
\begin{equation} \label{eq:HN}     
\H_\cN:= \Lambda/\left\{\lambda\in\Lambda\,|\,\cN\lambda=0\right\},\qquad    
\Lambda:= \Z[\Gamma]^m,    
\end{equation}    
where $(\cN\lambda)(h):=\sum_{g}\cN(hg^{-1})\lambda(g)$     
is the convolution product. This abelian group
is equipped with a bilinear pairing     
\begin{equation} \label{eq:quadratic}     
\inner{\mu}{\lambda}     
:= \sum_{g,h\in\Gamma}\mu(h)^T \cN(hg^{-1}) \lambda(g)     
= \sum_{h\in\Gamma}\mu(h)^T(\cN\lambda)(h)    
\end{equation}     
and a group action $\rho_\cN:\Gamma\to\Aut(\H_\cN)$ defined by     
\begin{equation}\label{eq:rhocN}     
\rho_\cN(g)[\lambda] := [g_*\lambda],    
\end{equation}     
which preserves the bilinear pairing.     
The special elements $\L_1,\dots,\L_m\in\H_\cN$ are given by    
\begin{equation} \label{eq:Li}     
\L_i := [\delta_i],\qquad      
\delta_i(g) := \left\{\begin{array}{ll}      
e_i,&\mbox{if }g=1,\\      
0,&\mbox{if }g\ne 1.     
\end{array}\right.      
\end{equation}      
These structures are well defined     
for any function $\cN:\Gamma\to\Z^{m\times m}$.      
The next lemma asserts that the tuple $(\H_\cN,\rho_\cN,\L_1,\dots,\L_m)$     
is a monodromy representation whenever $\cN$ is a monodromy character.     
    
\begin{lemma}\label{le:iso}     
{\bf (i)}     
Let $\cN:\Gamma\to\Z^{m\times m}$ be a monodromy character.     
Then the tuple $(\H_\cN,\rho_\cN,\L_1,\dots,\L_m)$      
defined by~(\ref{eq:HN}-\ref{eq:Li}) is a monodromy     
representation whose associated character is $\cN$.     
    
\smallskip\noindent{\bf (ii)}     
Let $(\cH,\rho,\cL_1,\dots,\cL_m)$ be a monodromy      
representation and $\cN$ be its     
character~(\ref{eq:N}). Then the map      
\begin{equation}\label{eq:Lambda}    
\Lambda \to\cH: \lambda \mapsto      
\sum_{g\in \Gamma}\sum_{i=1}^m \lambda_i(g) \rho(g^{-1})\cL_i     
\end{equation}    
induces an isomorphism of monodromy representations     
from $\H_\cN$ to $\cV/\cW$, where $\cV\subset\cH$ denotes     
the submodule generated by the vanishing cycles $\rho(g)\cL_i$     
and $\cW\subset\cV$ is the kernel of the intersection form.     
\end{lemma}     
    
\begin{proof} 
The $\Z$-module $\H_\cN$ is isomorphic to the       
image of the homomorphism $\lambda\mapsto\cN\lambda$
and hence is torsion free.
That the bilinear pairing in~\eqref{eq:quadratic}      
is nondegenerate follows directly from the definition.      
That it satisfies~(\ref{eq:GH1}) follows from~(\ref{eq:N1})      
and that it satisfies~(\ref{eq:GH2}) follows from~(\ref{eq:N2})      
and the identity  $\inner{\cL_i}{\cL_i}=n_{ii}(1)$.      
To prove~(\ref{eq:GH3}) fix an index $i$ and an element       
$\lambda\in\Lambda$. Abbreviate 
$$
\eps:=(-1)^{n(n+1)/2}
$$
and define $\lambda'\in\Lambda$ by       
$$      
\lambda' := (g_i)_*\lambda - \lambda       
+ \eps\inner{\delta_i}{\lambda}\delta_i.     
$$      
Then     
$$     
\lambda'(h) = \lambda(hg_i)-\lambda(h)      
+\eps\left(\sum_{g\in\Gamma}e_i^T\cN(g^{-1})\lambda(g)\right)\delta_i(h)     
$$     
and hence     
\begin{eqnarray*}      
(\cN\lambda')(k)      
&=& \sum_h\cN(kh^{-1})      
\biggl(\lambda(hg_i)-\lambda(h)     
+\eps \sum_ge_i^T\cN(g^{-1})\lambda(g)\delta_i(h)     
\biggr) \\      
&=&       
\sum_h \biggl(\cN(kg_ih^{-1})      
- \cN(kh^{-1}) +\eps \cN(k)E_i\cN(h^{-1})      
\biggr)\lambda(h) \\      
&=& 0.      
\end{eqnarray*}      
The last equation follows from~(\ref{eq:N3}).      
Hence the equivalence class of $\lambda'$ vanishes       
and so the tuple $(\H_\cN,\rho_\cN,\L_1,\dots,\L_m)$ 
satisfies~(\ref{eq:GH3}). Thus we have proved that 
$(\H_\cN,\rho_\cN,\L_1,\dots,\L_m)$ is a monodromy 
representation.  Its character is      
$$      
\inner{\delta_i}{g_*\delta_j}      
= \sum_{h,k}      
\delta_i(k)^T\cN(kh^{-1})\delta_j(hg)      
=n_{ij}(g).     
$$      
This proves~(i).      
 
We prove~(ii). The homomorphism~\eqref{eq:Lambda} is obviously surjective.    
That the preimage of $\cW$ under~\eqref{eq:Lambda} is the subspace      
$\left\{\lambda\in\Lambda\,|\,\cN\lambda=0\right\}$ follows from     
the identity    
$$    
\inner{\rho(h^{-1})\cL_i}{\Phi_\cN(\lambda)}    
= \sum_{g,j}n_{ij}(hg^{-1})\lambda_j(g) = (\cN\lambda)_i(h),
$$
where
$$
\Phi_\cN(\lambda):=\sum_{g,j}\lambda_j(g)\rho(g^{-1})\cL_j.
$$
Hence~(\ref{eq:Lambda}) induces an isomorphism      
$\Phi_\cN:\H_\cN\to\cV/\cW$ and it follows directly from    
the definitions that it satisfies the requirements of the lemma.    
\end{proof} 

\begin{remark}\label{rmk:monodromy}\rm 
Our geometric motivation for Lemma~\ref{le:iso} is the following. 
Let $\cN=\cN^X_c:\Gamma\to\Z^{m\times m}$ be 
the function associated to a Lefschetz fibration $X\to\Sigma$ 
and a distinguished configuration $c$ via~\eqref{eq:NXc}. 
Let $V\subset H_n(M)$ be the submodule 
generated by the vanishing cycles $\rho(g)L_i$ and
$W\subset V$ be the kernel of the intersection form on~$V$. 
Then the homomorphism
$
\Lambda\to V:\lambda \mapsto 
\sum_{g,i}\lambda_i(g) \rho(g^{-1})L_i 
$
descends to a $\Gamma$-equivariant   
isomorphism of $\Z$-modules
$
\Phi_\cN:\H_\cN\to V/W  
$
that identifies the pairing in~\eqref{eq:HN1} with 
the intersection form and maps the element $\L_i$ 
defined by~\eqref{eq:HN1} to the equivalence 
class $[L_i]\in V/W$.
\end{remark} 
    
\begin{example}\label{ex:free}\rm    
Assume that $\Gamma$ is generated freely by $g_1,\dots,g_m$.     
In this case the function $\cN:\Gamma\to\Z^{m\times m}$  
in Definition~\ref{def:Gamma} is completely determined 
by $N:=\cN(1)$. This matrix satisfies    
$$   
N^T = (-1)^nN,\qquad      
n_{ii} = \left\{\begin{array}{ll}      
0,&\mbox{if }n\mbox{ is odd},\\      
2(-1)^{n/2},&\mbox{if }n\mbox{ is even}.     
\end{array}\right.      
$$   
It determines a monodromy representation    
$$    
\H_N:=\Z^m/\ker N,\qquad     
\inner{\mu}{\lambda}:=\mu^TN\lambda,    
$$    
with special elements associated to the standard basis     
of $\Z^m$ and the $\Gamma$-action $\rho_N:\Gamma\to\Aut(\H_N)$    
uniquely determined by     
\begin{equation}\label{eq:rhongi}
\rho_N(g_i)=\one-(-1)^{n(n+1)/2}E_iN.     
\end{equation}
Here $E_i\in\Z^{m\times m}$ denotes the matrix  with the 
$i$-th entry on the diagonal equal to~$1$ and zeroes elsewhere. 
The function $\cN$ can be recovered from the matrix $N$ via 
$$
\cN(g) = N\rho_N(g),
$$
where $\rho_N:\Gamma\to \GL_m(\Z)$ is defined by 
equation~\eqref{eq:rhongi}. Moreover, the monodromy 
representations $\H_N$ and $\H_\cN$ are isomorphic. 
The isomorphism $\H_N\to\H_\cN$ assigns to each $[v]\in\H_N$      
the equivalence class of the function $\lambda_v:\Gamma\to\Z^m$      
with value $v$ at $1$ and zero elsewhere.     
Its inverse is induced by the map     
$
\Z[\Gamma]^m\to\Z^m:\lambda\mapsto     
\sum_g\rho_N(g^{-1})\lambda(g).     
$
\end{example}  

The $\Z$-module $M_m(\Z[\Gamma])$ of $m\times m$-matrices with entries 
in the group ring is naturally an algebra over the group ring. 
One can think of an element $\cM\in M_m(\Z[\Gamma])$ 
also as a function $\cM:\Gamma\to\Z^{m\times m}$     
with finite support. Then the product is given by    
$$    
\cM\cN(g) = \sum_{h\in\Gamma}\cM(gh^{-1})\cN(h)    
$$    
and this formula continues to be meaningful when only one    
of the factors has finite support.     
The group ring $\Z[\Gamma]$ is equipped with an involution    
${\lambda\mapsto\bar\lambda}$ given by    
$    
\bar\lambda(g):=\lambda(g^{-1}).    
$    
The {\bf conjugate transpose} of a matrix     
$\cM\in M_m(\Z[\Gamma])$ is defined by     
$$    
\cM^t(g) := \overline{\cM}(g)^T = \cM(g^{-1})^T    
$$    
and it satisfies $(\cM\cN)^t=\cN^t\cM^t$.     
    
Let $\cB$ be a group that acts covariantly on $\Gamma$ and denote 
the action by     
$    
\cB\times\Gamma\to\Gamma:(\sigma,g)\mapsto\sigma_*g.     
$  
In the intended application $\Gamma$ is the fundamental group    
of a Riemann surface with $m$ punctures and $\cB$ is the     
braid group on $m$ strings in the same Riemann surface     
with one puncture. The action of $\cB$ on $\Gamma$ extends linearly to    
an action on $\Z[\Gamma]$ by algebra automorphisms given by    
$$    
\sigma_*\lambda     
= \sum_{g\in\Gamma}(\sigma_*\lambda)(g)g
:= \sum_{g\in\Gamma}\lambda(g)\sigma_*g,    
\qquad    
(\sigma_*\lambda)(g) := \lambda((\sigma^{-1})_*g).    
$$    
This action extends to the $\Z$-module $\Z[[\Gamma]]$ of 
formal sums of elements of $\Gamma$ with integer coefficients. 
These correspond to arbitrary integer valued functions on $\Gamma$.     
So $\cB$ acts on $M_m(\Z[\Gamma])$ componentwise,     
or equivalently by    
\begin{equation}\label{eq:phi*M}    
\sigma_*\cM := \cM\circ (\sigma^{-1})_*:\Gamma\to\Z^{m\times m}.    
\end{equation}    
We then have      
\begin{equation} \label{eq:sigma-matrices}    
\sigma_*(\cM\cN)=(\sigma_*\cM)(\sigma_*\cN), \qquad     
\sigma_*(\cM^t)=(\sigma_*\cM)^t    
\end{equation}    
for $\cM,\cN\in M_m(\Z[\Gamma])$. The action of $\cB$     
on $M_m(\Z[\Gamma])$ induces an action on the group    
$\GL_m(\Z[\Gamma])$ of invertible elements of $M_m(\Z[\Gamma])$.      

The following notion plays a crucial role in this paper. 
In our intended application $\G$ is the braid group,
or the framed braid group, and $A$ is the group $\GL_m(\Z[\Gamma])$
of invertible matrices with entries in the group ring of $\Gamma$. 

\begin{definition}[Serre~\cite{Serre}]\label{def:cocycle}    
Let $\G$ and $A$ be groups. Suppose that    
$\G$ acts covariantly on $A$ and denote the action by    
$  
\G\times A\to A:(g,a)\mapsto g_*a.    
$   
A map $s:\G\to A$ is called a {\bf cocycle} if    
\begin{equation} \label{eq:cocycle}      
s(gh) = s(g)g_*s(h).
\end{equation}    
Two cocycles $s_0,s_1:\G\to A$ are called {\bf cohomologous}     
if there is an element $a\in A$ such that    
\begin{equation} \label{eq:coboundary}      
s_1(g) = a^{-1}s_0(g)g_*a.    
\end{equation}    
The set of equivalence classes of cocycles is     
denoted by $H^1(\G,A)$.    
\end{definition}    
    
\begin{remark}\label{rmk:cocycle}\rm    
The semidirect product $\G\ltimes A$    
is equipped with the group operation    
$$    
(g,a)\cdot(h,b) :=  (gh,a g_*b).    
$$    
A map $s:\G\to A$ is a cocycle if and only if the map
$   
\G\to\G\ltimes A:g\mapsto(g,s(g))    
$
is a homomorphism (and the homomorphisms associated     
to two cohomologous cocycles are conjugate by the element    
$(1,a)\in\G\ltimes A$). Observe that the cocycle condition implies $s(1)=1$.   
\end{remark}   
    
\begin{lemma}\label{le:cocycle}    
Every cocycle $\cS:\cB\to\GL_m(\Z[\Gamma])$ induces a 
contravariant action of $\cB$ on the space of all functions     
$\cN:\Gamma\to\Z^{m\times m}$ via    
\begin{equation}\label{eq:BM}    
\sigma^*\cN:=\big(\cS(\sigma)^t\cN\cS(\sigma)\big) \circ \sigma_*.     
\end{equation}    
Moreover, two cohomologous cocycles induce conjugate actions. 
\end{lemma}    
 
\begin{proof}    
For $\sigma,\tau\in\cB$ we have     
\begin{eqnarray*}    
(\sigma\tau)^*\cN      
& = & \big(\cS(\sigma\tau)^t \cN    
\cS(\sigma\tau)\big) \circ (\sigma \tau)_* \\    
& = & \big(\sigma_*(\cS(\tau))^t  \cdot     
\cS(\sigma)^t\cN\cS(\sigma)\cdot \sigma_*(\cS(\tau))\big)     
\circ\sigma_*\circ \tau_* \\    
& = & \big(\cS(\tau)^t \cdot \sigma^*\cN \cdot    
\cS(\tau)\big)\circ \tau_*  = \tau^*\sigma^*\cN.    
\end{eqnarray*}  
Here we have used the cocycle condition     
and~\eqref{eq:sigma-matrices}.   

Assume two cocycles $\cS$ and $\cS'$ are cohomologous, 
i.e.\ there 
is a matrix
$\cA\in\GL_m(\Z[\Gamma])$ such that 
$\cS'(\sigma)=\cA^{-1}\cS(\sigma)\sigma_*\cA$ for all 
$\sigma\in\cB$. Denote by $\sigma^*$ and $\sigma^{*'}$ 
the actions defined by $\cS$ and $\cS'$ respectively, and denote 
\begin{equation}\label{eq:CA}
C_\cA:\textrm{Map}(\Gamma,\Z^{m\times m})\to
\textrm{Map}(\Gamma,\Z^{m\times m}),\qquad \cN\mapsto\cA^t\cN\cA.
\end{equation}
A straightforward verification shows that we have 
$\sigma^{*'}=C_\cA\circ\sigma^*\circ C_{\cA}^{-1}$.
\end{proof}    
 
\begin{proposition}\label{prop:HN}    
Let $\cS:\cB\to\GL_m(\Z[\Gamma])$ be a     
cocycle. For every ${\sigma\in\cB}$ and    
every function $\cN:\Gamma\to\Z^{m\times m}$, the isomorphism     
$$
\Lambda\to\Lambda:\lambda\mapsto\cS(\sigma)(\sigma_*\lambda)
$$ 
descends to an isomorphism      
\begin{equation}
\label{eq:S_iso}
\bS_\sigma:\H_{\sigma^*\cN}\longrightarrow \H_\cN    
\end{equation}
which preserves the bilinear pairings and fits 
into a commutative diagram    
\begin{equation} \label{eq:action_rho}      
\xymatrix    
@C=50pt    
@R=25pt    
{    
\Gamma \ar[r]^-{\rho_{\sigma^*\cN}} \ar[d]_{\sigma_*}     
& \mathrm{Aut}(\H_{\sigma^*\cN}) \ar[d] \\    
\Gamma \ar[r]_-{\rho_\cN}     
& \mathrm{Aut}(\H_\cN)    
}
\end{equation}     
with the second vertical arrow given by    
$\phi\mapsto\bS_\sigma\phi\bS_\sigma^{-1}$.  
Moreover, two cohomologous cocycles induce 
equivalent isomorphisms~\eqref{eq:S_iso}.
\end{proposition} 
 
\begin{proof} 
For any $\sigma\in\cB$ and any function $\cM:\Gamma\to\Z^{m\times m}$    
the isomorphism $\sigma_*:\Lambda\to\Lambda$    
descends to an isomorphism    
$
R_\sigma:\H_{\cM\circ \sigma_*}\longrightarrow \H_\cM   
$
because 
$(\cM\circ \sigma_*)\lambda = (\sigma^{-1})_*(\cM\cdot\sigma_*\lambda)$.    
That $R_\sigma$ preserves the bilinear pairings follows from    
the identity    
$    
{(\sigma_*\mu)^t\cM (\sigma_*\lambda)    
= \sigma_*( \mu^t (\cM\circ \sigma_*)\lambda)}    
$    
and the formula $\inner{\lambda}{\mu}=(\lambda^t\cM\mu)(1)$    
for the pairing~(\ref{eq:quadratic}) on $\H_\cM$.     
    
For any $\sigma\in\cB$ we also have an isomorphism     
$
Q_\sigma:\H_{\cS(\sigma)^t\cN\cS(\sigma)} \longrightarrow \H_\cN    
$
induced by the map $\lambda\mapsto\cS(\sigma)\lambda$.    
This follows from the associativity of the convolution product     
$(\cS^t\cN\cS)\lambda = \cS^t(\cN\cS\lambda)$ and the fact that 
$\cS(\sigma)$ is invertible. That $Q_\sigma$ preserves 
the bilinear pairing follows from the identity     
$(\cS\mu)^t\cN(\cS\lambda) = \mu^t(\cS^t\cN\cS)\lambda$.    
 
The isomorphism $\bS_\sigma:\H_{\sigma^*\cN}\to\H_\cN$     
is the composition $\bS_\sigma:= Q_\sigma \circ R_\sigma$    
with $\cM:=\cS(\sigma)^t\cN\cS(\sigma)$.     
By what we have already proved, this isomorphism     
preserves the bilinear pairings. The commutativity of the    
diagram~\eqref{eq:action_rho} is equivalent to the equation     
\begin{equation} \label{eq:NgS}     
\rho_\cN(\sigma_*g)\circ \bS_\sigma 
= \bS_\sigma \circ \rho_{\sigma^*\cN}(g),    
\qquad g\in \Gamma.     
\end{equation}     
To prove~\eqref{eq:NgS}  note that    
$$    
(\rho_\cN(\sigma_*g)\circ \bS_\sigma)(\lambda)     
= \rho_\cN(\sigma_*g)\big(\cS(\sigma)\sigma_*\lambda\big)     
= \bigl(\cS(\sigma)\sigma_*\lambda\bigr)(\sigma_*g)^{-1},     
$$
and     
$$
(\bS_\sigma\circ \rho_{\sigma^*\cN}(g))(\lambda)     
= \bS_\sigma(\lambda g^{-1})     
= \cS(\sigma)\sigma_*(\lambda g^{-1})     
= \bigl(\cS(\sigma)\sigma_*\lambda\bigr)(\sigma_*g)^{-1}.     
$$

To prove the last 
assertion, we assume that $\cS$ and $\cS'$ are
two cohomologous cocycles. Thus there 
is a matrix 
$\cA\in\GL_m(\Z[\Gamma])$ 
such that $\cS'(\sigma)=\cA^{-1}\cS(\sigma)\sigma_*\cA$ for all $\sigma\in\cB$. 
Denote by $\sigma^*$ and $\sigma^{*'}$ the actions defined by 
$\cS$ and $\cS'$ respectively. We proved in Lemma~\ref{le:cocycle} 
the relation $\sigma^{*'}=C_\cA\circ\sigma^*\circ C_{\cA}^{-1}$, 
where the map $C_\cA$ is defined by~\eqref{eq:CA}. 
We then have a commutative diagram 
$$
\xymatrix{
\H_{\sigma^{*'}(\cN)}\ar[r]^-{L_\cA} 
\ar[d]_{\bS'_\sigma}Ê& \H_{\sigma^{*'}(\cN) \cA^{-1}} = 
\H_{\sigma^*((\cA^{-1})^t\cN \cA^{-1})} \ar[d]^{\bS_\sigma} \\
\qquad\qquad\qquad\H_\cN=\H_{(\cA^{-1})^t\cN} 
\ar[r]_-{L_\cA} & \H_{(\cA^{-1})^t\cN \cA^{-1}}.
}
$$
Here the two isomorphisms denoted $L_\cA$ are induced by left multiplication 
$\Lambda\to\Lambda: \lambda\mapsto \cA\lambda$. The commutativity 
of the diagram is checked at the level of $\Lambda$ 
using the relation between $\cS'$ and $\cS$. 
\end{proof}    


\section{The framed braid group}\label{sec:braid}

In this section we recall the well known correspondence     
between the braid group and the mapping class group
(see~\cite{BIRMAN,KT}) and extend it to the framed     
braid group introduced by Ko and Smolinsky~\cite{KS}.

Let $\Sigma$ be a compact oriented $2$-manifold, 
possibly with boundary, let $Z\subset\Sigma\setminus\p\Sigma$ 
be a finite set consisting of $m$ points, 
and choose a base point $z_0\in\Sigma\setminus Z$.
We assume throughout that $\Sigma$ is not diffeomorphic to the 
$2$-sphere and that $z_0\in\p\Sigma$ whenever $\Sigma$ is 
diffeomorphic to the $2$-disc.
Denote by $\Diff(\Sigma,z_0)$ the group of all      
diffeomorphisms of $\Sigma$ that fix $z_0$ and by      
$
\Diff_0(\Sigma,z_0)\subset\Diff(\Sigma,z_0)
$
the identity component. Define 
\begin{equation*}    
\begin{split}    
\cG&:=\left\{\phi\in\Diff_0(\Sigma,z_0)\,|\,\phi(Z)=Z\right\}, \\    
\cG_0&:=\left\{\phi\in\cG\,|\,\exists\phi_t\in\cG    
\,\,\mathrm{s.t.}\,\,\phi_0=\id,\,\phi_1=\phi    
\right\}.     
\end{split}    
\end{equation*}    
Here $[0,1]\to\Diff(\Sigma,z_0):t\mapsto\phi_t$ is a      
smooth isotopy of diffeomorphisms fixing the base point.     
Thus $\cG_0$ is the identity component of $\cG$.     
We refer to the quotient     
$$     
\cG/\cG_0 = \pi_0(\cG)  
$$     
as the {\bf mapping class group}. 
It is naturally isomorphic to the braid group.    
  
Fix an ordering $Z=\{z_1,\dots,z_m\}$ and     
let $\gS_m$ denote the group of permutations.    
The {\bf braid group} $\cB$ on $m$ strings in $\Sigma\setminus\{z_0\}$    
based at $Z$ is defined as the fundamental group      
of the configuration space of $m$ unordered distinct      
points in $\Sigma\setminus\{z_0\}$.      
Think of a braid as an $m$-tuple of smooth paths      
$\beta_i:[0,1]\to\Sigma$, $i=1,\dots,m$, which avoid $z_0$,      
are pairwise distinct for each $t$, and satisfy      
$\beta_i(0)=z_i$, ${\beta_i(1)=z_{\pi(i)}}$     
for some permutation $\pi\in\gS_m$.      
Thus $\cB$ is the group of homotopy classes of braids.     
The composition law is $[\beta]\cdot[\alpha]:=[\beta\alpha]$,      
where $\beta\alpha$ denotes the braid obtained by first running      
through $\alpha$ and then through $\beta$.      
     
The isomorphism      
$
\Phi:\cB\to\cG/\cG_0     
$
is defined as follows. Given a braid $\beta$ choose     
a smooth isotopy $\{\phi_t\}_{0\le t\le 1}$      
in $\Diff(\Sigma,z_0)$ with $\phi_0=\id$ satisfying      
$\phi_t(z_i)=\beta_i(t)$ for $i=1,\dots,m$ and define      
$    
\Phi([\beta]):=[\phi_1].     
$     
To see that this map is well defined choose     
a smooth family of vector fields $v_t\in\Vect(\Sigma)$     
satisfying $v_t(z_0)=0$ and     
$$     
v_t(\beta_i(t))=\dot\beta_i(t),\qquad i=1,\dots,m,   
$$     
and let $\phi_t$ be the isotopy generated by $v_t$ via     
$\p_t\phi_t = v_t\circ\phi_t$ and $\phi_0=\id$.     
The existence of $v_t$ follows from an easy argument      
using cutoff functions, and that the isotopy class     
of $\phi_1$ is independent of the choices of $\beta$      
and $v_t$ follows from a parametrized version of      
the same argument respectively from taking      
convex combinations of vector fields.      
We claim that $\Phi$ is an isomorphism. 
(See \cite[Theorem 4.3]{BIRMAN} for a slightly different statement.)

That $\Phi$ is a surjective group homomorphism is obvious.
Thus we have an exact sequence
\begin{equation}
\label{eq:bg_to_mcg}
\pi_1\left(\Diff_0(\Sigma,z_0),\id\right) \longrightarrow
\cB \stackrel\Phi\longrightarrow \cG/\cG_0 \longrightarrow 1,
\end{equation}
where the first map sends an isotopy $\{\phi_t\}_{0 \leq t \leq 1}$ 
with $\phi_0=\phi_1 =\id$ to the braid in $\Sigma \setminus \{z_0\}$ 
defined by $t \mapsto \phi_t(Z)$.  Hence the injectivity of $\Phi$ 
follows from the fact that $\Diff_0(\Sigma,z_0)$  is simply connected. 
(This is why we exclude the $2$-sphere and the $2$-disc with 
a base point in the interior: see Remark \ref{rmk:S2D} below.) 
In fact it is contractible.   To see this in the case 
$\chi(\Sigma)<0$ we use the fact that the identity component 
$\Diff_0(\Sigma)$ of $\Diff(\Sigma)$ is contractible~\cite{EE,ES} 
and consider the fibration     
$$
\Diff_0(\Sigma,z_0)\hookrightarrow 
\Diff_0(\Sigma)\to\widetilde{\Sigma}.
$$
The map $\Diff_0(\Sigma)\to\widetilde{\Sigma}$ assigns to every
diffeomorphism $\phi\in\Diff_0(\Sigma)$ the homotopy class of the 
path $[0,1]\to\Sigma:t\mapsto\phi_t(z_0)$ with fixed endpoints,
where $[0,1]\to\Diff_0(\Sigma):t\mapsto\phi_t$ is a smooth isotopy
with $\phi_0=\id$ and $\phi_1=\phi$. This map is well defined because 
$\Diff_0(\Sigma)$ is simply connected.
If the Euler characteristic is zero $\Sigma$ is either diffeomorphic 
to the $2$-torus or to the annulus.  In both cases
$\Diff_0(\Sigma,z_0)=\Diff(\Sigma,z_0)\cap\Diff_0(\Sigma)$
acts freely on $\cJ(\Sigma)$ and there is a diffeomorphism      
$\Diff_0(\Sigma,z_0)\times\cT(\Sigma)\to\cJ(\Sigma)$ 
where $\cT(\Sigma):=\cJ(\Sigma)/\Diff_0(\Sigma)$
is diffeomorphic to the upper half space $\H$ in the case 
of the torus and to an open interval in the case of the annulus.      
In the case of the disc $\D$ the group $\Diff_0(\D)$
of orientation preserving diffeomorphisms acts transitively 
on $\cJ(\D)$ with isotropy subgroup $\PSL(2,\R)$.      
This action gives rise to a diffeomorphism from $\cJ(\D)$     
to the subgroup of all diffeomorphisms that fix a point      
on the boundary and a point in the interior.  Hence the group     
$\Diff_0(\D,z_0)$ is contractible for every $z_0\in\p\D$.     

\begin{remark}\label{rmk:S2D}\rm
In the cases $\Sigma=S^2$ and $\Sigma=\D$ with $z_0\notin\p\D$
the group $\Diff_0(\Sigma,z_0)$ is not contractible but homotopy
equivalent to the circle.  It can be deduced from 
the exact sequence~\eqref{eq:bg_to_mcg} that $\Phi$ 
is not injective and that, instead, $\cB$ is in these two cases 
a central extension of the mapping class group $\cG/\cG_0$ 
by an infinite cyclic group $\cZ$.  If $\Sigma=S^2$ 
then $\cB$ can be interpreted as the braid group 
on $m$ strings in $\C$, and the subgroup $\cZ$ 
is the center of $\cB$ for $m\geq 3$~\cite{BIRMAN,KT}. 
If $\Sigma=\D$ with $z_0\notin\p\D$ then 
$\cB$ can be interpreted as the subgroup of the braid group
on $m+1$ strings in $\D$ that fixes the point $z_0$, and  
$\cZ$ coincides with the center of $\cB$ for all $m\ge1$.
\end{remark}

The framed braid group is an extension of the braid group $\cB$.       
Choose nonzero tangent vectors $v_z\in T_z\Sigma$     
for $z\in Z$ and define     
\begin{equation*}     
\begin{split}    
\tcG &:=\left\{\phi\in\Diff_0(\Sigma,z_0)\,|\,     
\phi(Z)=Z,\,d\phi(z)v_z=v_{\phi(z)}\,\,\forall z\in Z    
\right\}, \\    
\tcG_0 &:=\left\{\phi\in\tcG\,|\,      
\exists\phi_t\in\tcG\,\,\mathrm{s.t.}\,\,     
\phi_0=\id,\,\phi_1=\phi\right\}.    
\end{split}    
\end{equation*}     
Thus $\tcG_0$ is the identity component of $\tcG$.     
The {\bf marked mapping class group} is the quotient     
$$    
\tcG/\tcG_0=\pi_0(\tcG).    
$$
It is naturally isomorphic to the {\bf framed braid group}      
$\tcB$ on $m$ strings in $\Sigma\setminus\{z_0\}$.      
Fix again an ordering $Z=\{z_1,\dots,z_m\}$     
and denote $v_i:=v_{z_i}$ for $i=1,\dots,m$.     
The framed braid group is defined as the fundamental group      
of the configuration space of $m$ unordered      
points in the complement of the zero section in    
the tangent bundle of $\Sigma\setminus\{z_0\}$    
whose projections to the base are pairwise distinct.     
Think of a framed braid as an $m$-tuple of smooth paths      
$(\beta_i,\xi_i):[0,1]\to T\Sigma$ for $i=1,\dots,m$     
such that $(\beta_1,\dots,\beta_m)$ is a braid in    
$\Sigma\setminus\{z_0\}$ satisfying    
$\beta_i(0)=z_i$ and $\beta_i(1)=z_{\pi(i)}$    
for some permutation $\pi\in\gS_m$, and each $\xi_i$ is a    
nowhere vanishing vector field along $\beta_i$ such that    
$\xi_i(0)=v_i$ and $\xi_i(1)=v_{\pi(i)}$.    
Thus $\tcB$ is the group of homotopy classes of framed braids.     
 
The isomorphism      
$     
\tPhi:\tcB\to\tcG/\tcG_0     
$     
is defined as follows. Given a framed braid $(\beta,\xi)$ choose     
a smooth isotopy $\{\phi_t\}_{0\le t\le 1}$      
in $\Diff(\Sigma,z_0)$ with $\phi_0=\id$ satisfying      
\begin{equation}\label{eq:framed}    
\phi_t(z_i)=\beta_i(t),\qquad    
d\phi_t(z_i)v_i = \xi_i(t),\qquad     
i=1,\dots,m,     
\end{equation}    
and define      
$    
\tPhi([\beta,\xi]):=[\phi_1].     
$    
To see that this map is well defined choose     
a smooth family of vector fields $v_t\in\Vect(\Sigma)$     
satisfying $v_t(z_0)=0$ and     
\begin{equation}\label{eq:framedv}    
v_t(\beta_i(t))=\dot\beta_i(t),\qquad     
\Nabla{\xi_i(t)}v_t(\beta_i(t)) = \Nabla{t}\xi_i(t),\qquad    
i=1,\dots,m.    
\end{equation}    
(Here $\nabla$ is a torsion free connection on $T\Sigma$ but     
the second equation in~(\ref{eq:framedv}) is independent of this choice.)    
Now let $\phi_t$ be the isotopy generated by $v_t$ via     
$\p_t\phi_t = v_t\circ\phi_t$ and $\phi_0=\id$.     
Then $\phi_t$ satisfies~(\ref{eq:framed}).    
As above, the existence of $v_t$ follows from an easy argument      
using cutoff functions, and that the isotopy class     
of $\phi_1$ is independent of the choices of $\beta$      
and $v_t$ follows from a parametrized version of      
the same argument respectively from taking      
convex combinations of vector fields.      
That $\tPhi$ is a surjective group homomorphism is    
obvious and that it is injective follows again from the     
definitions and the fact that $\Diff_0(\Sigma,z_0)$     
is simply connected.     
   
\begin{remark}\label{rmk:braid}\rm   
There is an obvious action of the mapping class   
group on the braid group induced by the    
action of $\cG$ on $\cB$ via    
$$   
\phi_*[\beta_1,\dots,\beta_m]   
:=[\phi\circ\beta_{\pi(1)},\dots,\phi\circ\beta_{\pi(m)}],   
$$   
for $\phi\in\cG$ and a braid $\beta=(\beta_1,\dots,\beta_m)$,     
where $\pi\in\gS_m$ is the permutation defined by    
$\phi(z_{\pi(i)})=z_i$.  On the other hand    
we have seen that the mapping class group can be    
identified with the braid group. The resulting    
action of the braid group on itself is given by   
inner automorphisms.  In other words   
$$   
\Phi(\alpha)_*\beta = \alpha\beta\alpha^{-1}   
$$   
for $\alpha,\beta\in\cB$. The same holds for the framed braid group.   
\end{remark}   

The framed braid group fits into an exact sequence     
\begin{equation}\label{eq:exact}    
0\to\Z^m\to\tcB\to\cB\to 1.     
\end{equation}   
This extension splits by choosing a nowhere vanishing vector field 
$w$ on $\Sigma\setminus\{z_0\}$ such that $w_z=v_z$ at each $z\in Z$. 
The splitting depends on the homotopy class of $w$ relatively to $Z$,
so that it is not unique in general.
In the following we shall not distinguish in notation      
between the mapping class group $\pi_0(\cG)$ and     
the braid group $\cB$, nor between $\pi_0(\tcG)$ and $\tcB$.

    
\section{Distinguished configurations}\label{sec:bases}     
     
Let $\Sigma,Z,z_0$ be as in Section~\ref{sec:braid}.    
An $m$-tuple $c=(c_1,\dots,c_m)$ of smooth paths      
$c_i:[0,1]\to\Sigma$ is called a     
{\bf distinguished configuration} if      
\begin{description}     
\item[(i)]      
each $c_i$ is an embedding with $c_i(0)=z_0$ and, for $i\neq j$,     
the paths $c_i$ and $c_j$ meet only at $z_0$;     
\item[(ii)]      
$\{c_1(1),\ldots,c_m(1)\}=Z$;     
\item[(iii)]      
the vectors $\dot c_1(0),\ldots,\dot c_m(0)$ are pairwise     
linearly independent and are ordered clockwise in $T_{z_0} \Sigma$.       
\end{description}    
Two distinguished configurations $c^0$ and $c^1$ are called      
{\bf homotopic} if there is a smooth homotopy      
$\{c^\lambda\}_{0\le\lambda \le 1}$ of distinguished      
configurations from $c^0$ to $c^1$.      We write $c^0\sim c^1$
if $c^0$ is homotopic to $c^1$ and denote the homotopy class 
of a distinguished configuration $c$ by $[c]$.  The set of homotopy classes 
of distinguished configurations will be denoted by $\cC$. 
Note that each distinguished configuration $c$ determines 
an ordering $Z=\{z_1,\dots,z_m\}$ via $z_i:=c_i(1)$.     

\begin{theorem} \label{thm:BC}     
The braid group $\cB$ acts freely and transitively on $\cC$ via      
\begin{equation}\label{eq:BC}     
\big( [\phi], [c_1,\ldots,c_m] \big) \mapsto      
[\phi\circ c_1,\ldots,\phi \circ c_m]=:[\phi_* c].     
\end{equation}     
\end{theorem}     
     
\begin{proof}     
We prove that the action is transitive.     
Given two distinguished configurations $c$ and $c'$ we need to     
construct an element $\psi\in \Diff_0(\Sigma,z_0)$ such that     
$\psi_*c$ is homotopic to $c'$. Up to homotopy we can assume that      
there is a constant $\eps>0$ such that $c_i(t)=c_i'(t)$ for $0\le t\le\eps$.     
Now construct an isotopy      
$[\eps,1]\to\Diff_0(\Sigma,z_0):\lambda\mapsto\phi_\lambda$     
satisfying $\phi_1=\id$ and      
$$     
\phi_\lambda(c_i(t)) = c_i(\lambda t),\qquad     
\eps\le\lambda\le1,\quad i=1,\dots,m,     
$$     
by choosing an appropriate family of vector fields.     
Choose an analogous isotopy $\phi_\lambda'$ for $c'$      
and define    
$$
\psi:=(\phi_\eps')^{-1}\circ\phi_\eps  \in \cG.    
$$
Then $\psi(c_i(t))=c_i'(t)$ as required.     

We prove that the action is free when $z_0\notin\p\Sigma$.
The case $z_0\in\p\Sigma$ is similar.
Let $c$ be a distinguished     
configuration and $\phi\in\cG$ such that $\phi_*c$ is      
homotopic to $c$. We prove in five steps that $\phi\in\cG_0$.     
     
\medskip\noindent{\bf Step~1.}     
{\it We may assume that $d\phi(z_0)=\one$.}     
     
\medskip\noindent    
It suffices to prove that, for every matrix     
$A\in\GL^+(n,\R)$, there is a diffeomorphism     
$\psi:\R^n\to\R^n$, supported in the unit ball      
and isotopic to the identity through      
diffeomorphisms with support in the unit ball,      
that satisfies $\psi(0)=0$ and $d\psi(0)=A$.      
If $A$ is symmetric and positive definite      
we may assume that $A$ is a diagonal matrix     
and choose $\psi$ in the form      
$$
\psi(x)=(\psi_1(x_1),\dots,\psi_n(x_n))
$$
where each $\psi_i$ is a suitable monotone      
diffeomorphism of $\R$. If $A$ is orthogonal     
we choose a smooth path $[0,1]\to\SO(n):r\mapsto A_r$,     
constant near the ends, with $A_0=A$ and $A_1=\one$     
and define 
$$
\psi(x):=A_{\abs{x}}x. 
$$
The general case follows by polar decomposition.      
     
\medskip\noindent{\bf Step~2.}     
{\it We may assume that $\phi$      
agrees with the identity near $z_0$.}     
     
\medskip\noindent     
Let $\phi$ be a diffeomorphism of $\R^n$      
with $\phi(0)=0$ and $d\phi(0)=\one$     
and choose a smooth nonincreasing cutoff function      
$\beta:[0,1]\to[0,1]$ equal to one near zero and     
vanishing near one. Then, for $\eps>0$ sufficiently small,      
the formula     
$$    
\phi_\lambda(x) := \lambda\beta(\abs{x}/\eps) x    
+(1-\lambda\beta(\abs{x}/\eps))\phi(x)     
$$    
defines an isotopy from $\phi_0=\phi$ to a      
diffeomorphism $\phi_1$ equal to the identity      
near the origin such that, for each $\lambda$,     
$\phi_\lambda$ agrees with $\phi$ outside the ball      
of radius $\eps$. Now choose a local coordinate chart     
near $z_0$ to carry this construction over to $\Sigma$.     
     
\medskip\noindent{\bf Step~3.}     
{\it We may assume that $\phi$      
agrees with the identity near $z_0$     
and $\phi_*c=c$.}    
     
\medskip\noindent     
Assume, by Step~2, that $\phi$ agrees with the identity      
near $z_0$. Then the homotopy $\lambda\mapsto c^\lambda$      
from $c^0=c$ to $c^1=\phi_*c$ can be chosen such that     
$c^\lambda(t)$ is independent of $t$ for $t$ sufficiently     
small. Hence there exists a family of vector fields      
$v^\lambda\in\Vect(\Sigma)$ satisfying      
$$     
v^\lambda(c^\lambda_i(t)) = \p_\lambda c^\lambda_i(t)     
$$     
for all $\lambda,t,i$ and $v^\lambda(z)=0$ for $z$ near $z_0$.      
Integrating this family of vector fields yields a      
diffeomorphism $\psi\in\cG_0$ such that $\psi_*\phi_*c=c$.      
     
\medskip\noindent{\bf Step~4.}     
{\it We may assume that $\phi$      
agrees with the identity near the      
union of the images of the $c_i$.}     
     
\medskip\noindent     
Assume, by Step~3, that $\phi(c_i(t))=c_i(t)$ for all $t$ and $i$     
and that $\phi$ agrees with the identity near $z_0$.      
If $d\phi(c_i(t))=\one$ for all $i$ and $t$ we can      
use an interpolation argument as in Step~2 to deform     
$\phi$ to a diffeomorphism that satisfies the      
requirement of Step~4. To achieve the condition      
$d\phi(c_i(t))=\one$ via a prior deformation     
we must solve the following problem. Given two     
smooth functions $a:[0,1]\to\R$ and $b:[0,1]\to(0,\infty)$     
find a diffeomorphism $\psi:\R^2\to\R^2$ that     
is equal to the identity outside a small neighborhood     
of the set $[0,1]\times\{0\}$ and satisfies     
$$     
\psi(x,0)=(x,0),\qquad     
d\psi(x,0) = \left(     
\begin{array}{cc}     
1 & a(x) \\ 0 & b(x)     
\end{array}     
\right),\qquad     
0\le x\le 1.     
$$     
It suffices to treat the cases $a(x)\equiv0$ and $b(x)\equiv1$.     
For $b(x)\equiv1$ one can use an interpolation argument      
as in the proof of Step~2. For $a(x)\equiv0$ one can use a      
parametrized version of the argument for the positive      
definite case in the proof of Step~1.      
     
\medskip\noindent{\bf Step~5.}     
{\it We prove that $\phi\in\cG_0$.}     
     
\medskip\noindent     
Assume, by Step~4, that there is a coordinate chart $u:\D\to\Sigma$      
such that $c_i(t)\in u(\D)$ for all $i$ and $t$ and $\phi\circ u=u$.      
Choose any isotopy     
$$
[0,1]\to\Diff_0(\Sigma,z_0):\lambda\mapsto\phi_\lambda
$$
from $\phi_0=\id$ to $\phi_1=\phi$.     
By a parametrized version of the argument in Step~1     
we may assume that $d\phi_\lambda(z_0)=\one$ for every $\lambda$.    
By a parametrized version of the argument in Step~2     
we may assume that there is an $\eps>0$ such that     
$\phi_\lambda\circ u$ agrees with $u$ on the disc      
of radius $\eps$ for every $\lambda$.     
Choose a diffeomorphism $\psi:\Sigma\to\Sigma$,      
supported in $u(\D)$, such that     
$$     
\psi(u(\left\{z\in\D\,|\,\abs{z}\le1-\eps\right\}))      
= u(\left\{z\in\D\,|\,\abs{z}\le\eps\right\}).    
$$     
Then $\lambda\mapsto\psi^{-1}\circ\phi_\lambda\circ\psi$     
is an isotopy in $\cG$ and so $\phi=\psi^{-1}\circ\phi\circ\psi\in\cG_0$.      
This proves the theorem.     
\end{proof}     

Recall from the introduction that every distinguished configuration $c$
determines elements $g_{1,c},\dots,g_{m,c}$ of the fundamental group
$\Gamma=\pi_1(\Sigma\setminus Z,z_0)$, where $g_{i,c}$ is the homotopy
class of the loop obtained by traversing $c_i$, encircling $z_i$ 
counterclockwise, and then traversing $c_i$ in the opposite direction.  
Clearly, the $g_{i,c}$ depend only on the homotopy class of $c$. 
Conversely, we have the following theorem.

\begin{theorem}\label{thm:C}
If $c,c'\in\cC$ satisfy $g_{i,c}=g_{i,c'}$ 
for $i=1,\dots,m$ then $c\sim c'$. 
\end{theorem}

The proof relies on the classical result of  Baer, Dehn and Nielsen that, in dimension two, isotopy coincides with homotopy.
Specifically, we need the following theorem of 
Epstein~\cite[Theorem~3.1]{Epstein} and Feustel~\cite{Feustel} 
about embedded arcs in $2$-manifolds.

\medskip\noindent {\bf The Epstein--Feustel Theorem.} 
{\it Let $S$ be a compact $2$-manifold with boundary and 
$\alpha,\beta:[0,1]\to S$ be smooth embeddings 
such that
$$
\alpha(0)=\beta(0),\qquad\alpha(1)=\beta(1),\qquad
\alpha^{-1}(\p S)=\beta^{-1}(\p S) = \left\{\alpha(0),\alpha(1)\right\}.
$$
If $\alpha$ and $\beta$ are smoothly homotopic with fixed endpoints 
then there is a smooth ambient isotopy 
$[0,1]\times S\to S:(\lambda,p)\mapsto\phi_\lambda(p)$ 
such that
$$
\phi_0=\id,\qquad
\phi_1\circ\alpha=\beta
$$
and $\phi_\lambda|_{\p S}=\id$ for all $\lambda\in[0,1]$.}

\medskip
In the work of Epstein and Feustel the $2$-manifold $S$ is triangulated
and the isotopy can be chosen piecewise linear whenever the arcs
are piecewise linear. In Feustel's theorem the homeomorphism 
$\phi_\lambda:S\to S$ fixes the endpoints of the arcs. 
In Epstein's theorem $S$ need not be compact and the 
$\phi_\lambda$ have uniform compact support
and are equal to the identity on the boundary of $S$. 
To obtain the smooth isotopy in the above formulation, 
we first approximate the embedded arcs by piecewise linear arcs, 
then use Epstein's version of the theorem in the piecewise linear setting, 
then approximate the piecewise linear isotopy by a smooth isotopy, 
and finally connect two nearby smooth arcs by a smooth isotopy.
 
\begin{proof}[Proof of Theorem~\ref{thm:C}]
Let $\zeta_z \subset \Gamma$ be the conjugacy class 
determined by a small loop encircling the puncture $z\in Z$.
Since $\Sigma \neq S^2$
we have  $\zeta_z \neq \zeta_{z'}$ whenever $z\neq z'$.
Since $g_{i,c} \in \zeta_{c_i(1)}$ and $g_{i,c'} \in \zeta_{c'_i(1)}$ 
for $i=1,\dots,m$, we deduce that $c$ and $c'$ 
determine the same ordering of $Z$:
$$
Z=\{z_1,\dots,z_m\}, \qquad z_i:=c_i(1)=c_i'(1).
$$
Performing an isotopy of the distinguished configuration $c'$,
if necessary, we may assume that $c_i'$ agrees with $c_i$
on the interval $[1-2\eps,1]$ for some $\eps>0$.  
Next we denote by $B_{2\eps}\subset\C$ the disc of radius $2\eps$ 
centered at zero and choose embeddings $\psi_i:B_{2\eps}\to\Sigma$ 
with disjoint images $U_i:=\psi_i(B_{2\eps})$ such that 
$\psi_i(t)=c_i(1-t)$ for $0\le t<2\eps$ and $c_i|_{[0,1-2\eps]}$ 
takes values in the complement of $U_1\cup\cdots\cup U_m$ for every $i$. 
Let $D_i:=\psi_i(B_\eps)$ and denote 
$$
D:=D_1\cup\cdots\cup D_m\subset\Sigma.
$$
Then $\Sigma\setminus D$ is a manifold with boundary and 
the inclusion ${\Sigma\setminus D\hookrightarrow\Sigma\setminus Z}$ 
induces an isomorphism of fundamental groups 
$\pi_1(\Sigma\setminus D,z_0)\cong\Gamma$.
We prove in four steps that the distinguished 
configurations $c$ and $c'$ are homotopic.

\medskip\noindent{\bf Step~1.}
{\it For $i=1,\dots,m$ let $h_i\in\Gamma$ be the homotopy class 
of the based loop that traverses $c_i|_{[0,1-\eps]}$ 
and then $c_i'|_{[0,1-\eps]}$ in the reverse direction.
Then, for every $i$, there is a $k_i\in\Z$ 
such that $h_i=g_{i,c}^{k_i}$.}

\medskip\noindent
By definition of $h_i$ we have
$$
h_ig_{i,c}h_i^{-1} = g_{i,c'} = g_{i,c} \;\in\;\Gamma
$$
for $i=1,\dots,m$. Thus $h_i$ commutes with $g_{i,c}$. 
Moreover, $\Gamma$ is a free group.  If $m>1$
or $\p\Sigma\ne\emptyset$ we can choose a basis of $\Gamma$ 
such that $g_{i,c}$ is one of the generators and the assertion follows.  

If $m=1$ and $\p\Sigma=\emptyset$ then $\Sigma$ has genus 
$g\ge 1$, by assumption. Hence the group $\Gamma$ is free 
of rank $2g$ and we can choose generators
$\alpha_1,\dots,\alpha_g,\beta_1,\dots,\beta_g$ such that
$
b:=g_{1,c}=\prod_{j=1}^g[\alpha_j,\beta_j]. 
$
Let $h:=h_1$. Since $\Gamma$ is free and $b$ commutes with $h$,
the subgroup of $\Gamma$ generated by $b$ and $h$ is 
free, by the Nielsen--Schreier theorem, and abelian and hence 
has rank one.  Thus there is a $d\in\Gamma$ such that $b=d^r$ 
and $h=d^k$ for some $r,k\in\Z$. Since $b\ne1$ 
we must have $r\ne 0$.  We claim that $r=\pm 1$. 
To see this, let $\Gamma=:\Gamma_1\supset \Gamma_2\supset\cdots$ 
be the lower central series of $\Gamma$, with 
$\Gamma_{\ell+1}:=[\Gamma_\ell,\Gamma_1]$ for $\ell\ge 1$.
The quotient $H:=\Gamma_1/\Gamma_2$ is 
free abelian of rank $2g$ with generators 
$[\alpha_1],\dots,[\alpha_g],[\beta_1],\dots,[\beta_g]$.
In this quotient the identity $d^r=b$ becomes $r\cdot [d]=[b]=0$. 
Since $r\ne 0$ and $H$ has no torsion we obtain $[d]=0$ in $H$, 
i.e.\ $d\in\Gamma_2$. Since $b\in\Gamma_2$ we can consider the 
identity $d^r=b$ in the quotient $\Gamma_2/\Gamma_3$. 
This quotient is canonically isomorphic to the second exterior 
power of $H$, by identifying the equivalence class of a commutator
$[u,v]$ with $[u]\wedge [v]$, for all $u,v\in\Gamma$. 
In $\Lambda^2 H$ the equivalence class of $b$ is equal to 
$\sum_{j=1}^g [\alpha_j]\wedge[\beta_j]$. 
This is a primitive element, and the equation $r\cdot [d]=[b]$ 
implies $r=\pm 1$. Thus Step~1 is proved.

\medskip \noindent {\bf Step~2.} 
{\it We may assume without loss of generality that, for each $i$, 
the paths $c_i|_{[0,1-\eps]}$ and $c_i'|_{[0,1-\eps]}$ are 
smoothly homotopic with fixed endpoints in $\Sigma\setminus D$.
Moreover, the path $c_i|_{[1-\eps,1]}$ agrees with 
the path $c_i'|_{[1-\eps,1]}$ for each $i$.}

\medskip\noindent 
Let $k_i$ be as in Step~1 and choose a smooth
cutoff function $\rho:[0,2\eps]\to\R$ such that 
$\rho(r)=1$ for $r\le\eps$ and $\rho(r)=0$ for $r\ge 3\eps/2$. 
Define the diffeomorphism $\phi:\Sigma\to\Sigma$ by 
$\phi(\psi_i(z)):=\psi_i(e^{-2\pi{\mathrm{i}}k_i\rho(\Abs{z})}z)$
for $z\in B_{2\eps}$ and $i=1,\dots,m$ and by
$\phi(p):=p$ for $p\notin U_1\cup\cdots\cup U_m$.
Replacing $c'$ by the equivalent distinguished configuration
$\phi_*c'$ and constructing $h_i$ as in Step~1, 
we obtain $h_i=1$ and this proves Step~2. 

\medskip\noindent{\bf Step~3.}
{\it If $z_0\notin\p\Sigma$ we may assume without loss of generality that 
there is a smooth embedding $\psi_0:B_{2\eps}\to\Sigma$ and real 
numbers $2\pi>\theta_1>\cdots>\theta_m\ge 0$
such that the following holds.}

\smallskip\noindent{\bf (a)}
{\it The closure of $U_0:=\psi_0(B_{2\eps})$
is disjoint from $\overline U_i$ for $i=1,\dots,m$.}

\smallskip\noindent{\bf (b)}
{\it For each $i$ we have 
$c_i(t)=c_i'(t)= \psi_0(e^{{\mathrm{i}}\theta_i}t)$ for $0\le t<\eps$,
$\Abs{\psi_0^{-1}(c_i(t))}=\Abs{\psi_0^{-1}(c_i'(t))}=t$ 
for $\eps\le t<2\eps$,
and $c_i(t),c_i'(t)\notin U_0$ for $2\eps\le t\le 1$.}

\smallskip\noindent{\bf (c)}
{\it The curves $c_i|_{[\eps,1-\eps]}$ and $c_i'|_{[\eps,1-\eps]}$
are smoothly homotopic with fixed endpoints in 
$\Sigma\setminus(D_0\cup D)$, where $D_0:=\psi_0(B_\eps)$.}

\smallskip\noindent
{\it In the case $z_0\in\p\Sigma$ we may assume the same with a smooth 
embedding $\psi_0:\left\{z\in B_{2\eps}\,|\,\mathrm{Im}\,z\ge 0\right\}\to\Sigma$
of a half disc and with $\pi>\theta_1>\cdots>\theta_m>0$.}

\medskip\noindent
Assume $z_0\notin\p\Sigma$ and choose any embedding 
$\psi_0:B_{2\eps}\to\Sigma$ with $\psi_0(0)=z_0$ such that~(a) holds. 
Reparametrizing $c_i$ near $t=0$ we may assume that 
$
{\Abs{\p_t(\psi_0^{-1}\circ c_i)(0)}=1}.
$
Rotating the embedding, if necessary, we may assume that
there are real numbers $2\pi>\theta_1>\theta_2>\cdots>\theta_m\ge 0$
such that 
$$
\p_t(\psi_0^{-1}\circ c_i)(0)=e^{{\mathrm{i}}\theta_i},\qquad
i=1,\dots,m. 
$$
Shrinking $\eps$, if necessary, we can then 
deform $c_i$ to a curve that satisfies 
$c_i(t)= \psi_0(e^{{\mathrm{i}}\theta_i}t)$ for $0\le t<2\eps$.
Shrinking $\eps$ again we may assume that $c_i(t)\notin U_0$ 
for $2\eps\le t\le1$. Applying the same argument to the $c_i'$, 
using partial Dehn twists supported in $U_0$, 
and shrinking $\eps$ again, we may assume that 
the arcs $c_i'$ satisfy the same two properties. 
Thus condition~(b) is fulfilled.

\smallbreak

By Step~2 the curves $c_i|_{[0,1-\eps]}$ and $c_i'|_{[0,1-\eps]}$
are homotopic with fixed endpoints in $\Sigma\setminus D$. 
Now consider the loop $\gamma_i$ in $\Sigma\setminus(D_0\cup D)$ 
obtained by first traversing $c_i|_{[\eps,1-\eps]}$ and then 
traversing $c_i'|_{[\eps,1-\eps]}$ in the reverse direction. 
This loop is contractible in $\Sigma\setminus D$ and hence, 
as a based loop in $\Sigma\setminus(D_0\cup D)$ with 
basepoint~$c_i(\eps)$, is homotopic to a multiple of $\p \overline D_0$.
Using a suitable multiple of a Dehn twist in the annulus $U_0\setminus D_0$ 
(as we did in Step~2), we can replace $c'$ by an equivalent distinguished 
configuration which still satisfies (b) and such that
$\gamma_1$ is now contractible in $\Sigma\setminus(D_0\cup D)$.
Hence the curves $c_1|_{[\eps,1-\eps]}$ and $c_1'|_{[\eps,1-\eps]}$
are homotopic with fixed endpoints in $\Sigma\setminus(D_0\cup D)$.
For $i\ge2$ the curves $c_i([\eps,1-\eps])$ and 
$c_i'([\eps,1-\eps])$ in $\Sigma\setminus(D_0\cup D)$
have endpoints different from $c_1(\eps)$ and $c_1(1-\eps)$. 
Hence they have well defined intersection numbers with 
$c_1|_{[\eps,1-\eps]}$.  By what we have just proved these 
agree with the intersection numbers with $c_1'|_{[\eps,1-\eps]}$.  
Since $c_i([\eps,1-\eps])$ 
is disjoint from $c_1([\eps,1-\eps])$ and $c_i'([\eps,1-\eps])$ 
is disjoint from $c_1'([\eps,1-\eps])$ we deduce that both
intersection numbers are zero. Hence the intersection number of 
$c_1|_{[\eps,1-\eps]}$ with $\gamma_i$ is zero for $i\ge 2$.
Since the loop $\gamma_i$ is a multiple of 
$\p \overline D_0$, we deduce that it is contractible 
in $\Sigma\setminus(D_0\cup D)$ for  $i\geq 2$.
This proves Step~3 in the case $z_0\notin\p\Sigma$. 
The proof in the case $z_0\in\p\Sigma$ is similar,
assertion (c) being simpler to prove.

\medskip\noindent{\bf Step~4.}
{\it The distinguished configurations $c$ and $c'$ are homotopic.}

\medskip\noindent 
We prove by induction on $\ell\in\{1,\dots,m\}$ that there is
an ambient isotopy 
$
[0,1]\times\Sigma\to\Sigma:(\lambda,p)\mapsto\phi_\lambda(p)
$
such that each $\phi_\lambda$ is the identity on $D_0\cup D$ and  
$\phi_1(c_i(t))=c_i'(t)$ for $0\le t\le 1$ and $i=1,\dots,\ell$.
For $\ell=1$ the existence of the isotopy follows immediately 
from Step~3, the Epstein--Feustel theorem, 
and the second assertion of Step~2.

Now suppose by induction that $\ell\in\{2,\dots,m\}$ 
and that $c_i'=c_i$ for $i=1,\dots,\ell-1$. 
By Step~3 the curves $c_\ell|_{[\eps,1-\eps]}$ and 
$c_\ell'|_{[\eps,1-\eps]}$ are homotopic with fixed 
endpoints in $\Sigma\setminus(D_0\cup D)$. 
Choose a smooth open disc $U\subset\Sigma$
(respectively half disc in the case $z_0\in\p\Sigma$) 
such that $\overline U$ is an embedded closed disc 
(respectively half disc) and 
\begin{equation*}
\begin{array}{ll}
D_0\cup D_1\cup\cdots\cup D_{\ell-1}\subset U,&
\bigcup_{i=1}^{\ell-1}c_i([0,1])\subset U, \\\\
(\overline D_\ell\cup\cdots\cup \overline D_m)
\cap\overline U=\emptyset,&
\bigcup_{i=\ell}^m c_i((\eps,1])\cap\overline U=\emptyset.
\end{array}
\end{equation*}
Then the inclusion of $\Sigma\setminus(U\cup D)$
into $\Sigma\setminus(D_0\cup D)$ induces an injection
of fundamental groups.  Hence the curves 
$c_\ell|_{[\eps,1-\eps]}$ and $c_\ell'|_{[\eps,1-\eps]}$ 
are homotopic with fixed endpoints in $\Sigma\setminus(U\cup D)$. 
Hence the existence of an ambient isotopy satisfying the assertion
for $\ell$ follows from the Epstein--Feustel theorem.
This proves Step~4 and the theorem.
\end{proof}

\begin{corollary}
\label{cor:Artin}
Let $N$ be the kernel of the homomorphism $\Gamma\to \pi_1(\Sigma,z_0)$ 
induced by the inclusion $\Sigma \setminus Z \hookrightarrow \Sigma$.
Then the homomorphism
\begin{equation}
\label{eq:Artin_representation}
\cB \longrightarrow \Aut(N):\sigma \longmapsto \sigma_*
\end{equation}
obtained by composing $\Phi: \cB \to \cG/\cG_0$
 with the canonical action of $\cG/\cG_0$ on the subgroup $N$, is injective.
\end{corollary}

\begin{proof}
That the canonical action of $ \cG/\cG_0$ on $\Gamma$ 
leaves the subgroup $N$ globally invariant follows from the definition of $\cG$.
To prove the injectivity, we consider a braid $\sigma \in \cB$ 
such that $\sigma_*$ is the identity of $N$. 
For an arbitrary configuration $c\in \cC$, we have
$$
g_{i,\sigma_*c}=\sigma_*g_{i,c}=g_{i,c},\qquad
i=1,\dots,m.
$$ 
Hence it follows from Theorem~\ref{thm:C} that $c$
and $\sigma_*c$ are equivalent distinguished configurations.
By Theorem~\ref{thm:BC} this implies that $\sigma$ is trivial.
\end{proof}

\begin{remark}\label{rmk:artin}\rm
Assume $\Sigma=\D$ with $z_0\in\p\D$. 
We have $N=\Gamma$ in this case and the map
(\ref{eq:Artin_representation}) is known as the {\bf Artin representation}. 
A classical theorem by Artin~\cite{Artin, BIRMAN,KT} 
asserts that it is injective and that its image consists 
of all automorphisms $\phi:\Gamma\to\Gamma$
that satisfy the following two conditions:
\begin{description}
\item[(i)]
$\phi$ permutes the $m$ conjugacy classes  in $\Gamma$ 
determined by small loops encircling the $m$ punctures;
\item[(ii)]  $\phi$ preserves the homotopy class of $\p\D$.
\end{description}
\noindent
Thus Corollary \ref{cor:Artin} is the injectivity part of Artin's theorem.
\end{remark}

Let us now choose a nonzero tangent vector $v_z \in T_{z} \Sigma$
at each puncture $z \in Z$. A {\bf marked distinguished configuration}     
is a distinguished configuration $c=(c_1,\ldots,c_m)$ satisfying     
$$
\dot c_i(1)=-v_{c_i(1)},\qquad i=1,\dots,m.
$$ 
Observe that the configuration $c$ induces an ordering on the set 
$\{v_z\}_{z\in Z}$ defined by $v_i := v_{c_i(1)}$ for $i=1,\dots,m$.
The notion of homotopy carries over to marked distinguished     
configurations and the set of homotopy classes will be 
denoted by $\tcC$. Now the proof of Theorem~\ref{thm:BC}      
carries over word by word to the present situation     
and shows the following.     
 
\begin{theorem}\label{thm:tBC} 
The framed braid group $\tcB$ acts freely and 
transitively on the set $\tcC$ via~(\ref{eq:BC}). 
\end{theorem} 

\begin{remark}\label{rmk:sigmaic}\rm
Given a marked distinguished configuration $c\in\tcC$, 
one can define elements 
$
\sigma_{2,c},\dots,\sigma_{m,c}, 
\eps_{1,c},\dots,\eps_{m,c}
$
in $\tcB$ as follows. 
\begin{itemize}
\item 
For $i=2,\dots,m$ we define the framed braid $\sigma_{i,c}$ 
as follows.  We choose an embedded arc 
$s_i:[0,1]\to\Sigma\setminus\{z_0\}$ from 
$z_{i-1}=c_{i-1}(1)=s_i(0)$ to $z_i=c_i(1)=s_i(1)$ 
by catenating $\left(c_{i-1}|_{[\eps,1]}\right)^{-1}$
with a clockwise arc from $c_{i-1}(\eps)$ to $c_i(\eps)$
and with $c_i|_{[\eps,1]}$. Given $s_i$ we choose a braid 
$\beta=(\beta_1,\dots,\beta_m)$ such that $\beta_{i-1}$ 
runs from $z_{i-1}$ to $z_i$ on the left of $s_i$, $\beta_i$
runs from $z_i$ to $z_{i-1}$ on the right of $s_i$,
and $\beta_j\equiv z_j=c_j(1)$ for $j\ne i-1,i$. 
The framing is determined by a vector field 
near the union of the curves $c_j$ which is 
tangent to the curves $c_j$ and has $z_0$ 
as an attracting fixed point. The mapping class
associated to $\sigma_{i,c}$ is represented by a 
diffeomorphism supported in an annulus around the 
geometric image of $\beta_{i-1}$ and $\beta_i$; 
it consists of two opposite half Dehn twists,
one in each half of this annulus, followed by
localized counterclockwise half turns centered
at $z_{i-1}$ and $z_i$. In terms of its 
action on $c$, the braid $\sigma_{i,c}$ preserves 
the curves $c_j$ for $j\ne i-1,i$ and replaces the pair    
$(c_{i-1},c_i)$ by $(c_ig_{i-1,c}^{-1},c_{i-1})$. 
\item 
For $i=1,\dots,m$ the framed braid $\eps_{i,c}$ is 
the trivial braid with the framing given by    
a counterclockwise turn about $z_i=c_i(1)$ and the trivial 
framing over $z_j$ for $j\ne i$. In terms of its action on $c$, 
the braid $\eps_{i,c}$ preserves the curves $c_j$ 
for $j\ne i$ and replaces $c_i$ by $c_ig_{i,c}$.  
\end{itemize}
In the case $\Sigma=\D$ there is an isomorphism 
$\phi_c:\tcB_m\to\tcB$, where $\tcB_m$ is the abstract
braid group with generators $\sigma_2,\dots,\sigma_m$ 
and $\eps_1,\dots,\eps_m$ subject to the 
relations~\eqref{eq:frabraid}, as introduced 
in the introduction.  The isomorphism sends 
$\sigma_i$ to $\sigma_{i,c}$ for $i=2,\dots,m$ 
and $\eps_j$ to $\eps_{j,c}$ for $j=1,\dots,m$. 
We refer to Figure~\ref{fig:SiEpsCoc} 
on page~\pageref{fig:SiEpsCoc} for a pictorial 
representation of these generators. 
\end{remark}

 
\section{The Picard--Lefschetz monodromy cocyle} \label{sec:coc} 

The main result of this section is Theorem~\ref{thm:Sc}, which
contains Theorem~A. We use the notations of the introduction. 
   
A marked distinguished configuration $c$ and a framed braid    
${\sigma=[\phi]\in\tcB}$ with $\phi\in\tcG$ determine   
a permutation $\pi_{\sigma,c}\in\gS_m$ such that    
$\phi(z_i) = z_{\pi_{\sigma,c}(i)}$ for $i=1,\dots,m$.     
These permutations satisfy   
\begin{equation}\label{eq:Pi}   
\pi_{\sigma\tau,c}=\pi_{\sigma,c}\circ\pi_{\tau,c},\qquad   
\pi_{\sigma,\tau_*c}=\pi_{\tau,c}^{-1}\circ\pi_{\sigma,c}\circ\pi_{\tau,c}.   
\end{equation}   
For $c\in\tcC$ and $j=1,\dots,m$ define the function    
$s_{j,c}:\tcB\to\Gamma$ by   
\begin{equation}\label{eq:sjc}   
s_{j,c}(\sigma) := c_i^{-1}\cdot\sigma_*c_j,\qquad    
i:=\pi_{\sigma,c}(j).   
\end{equation}   
Here the right hand side denotes the catenation of the paths   
$\sigma_*c_j$ and $c_i^{-1}$ pushed away from $z_i$    
in the common tangent direction $v_i$.   

\begin{remark}\label{rmk:sc}\rm   
Recall the elements $\sigma_{2,c},\dots,\sigma_{m,c}\in\tcB$ 
and $\eps_{1,c},\dots,\eps_{m,c}\in\tcB$ defined in Remark~\ref{rmk:sigmaic}. 
For $\sigma=\sigma_{k,c}$ the permutation   
$\pi_{\sigma,c}$ is the transposition of ${k-1}$ and $k$,   
for $\sigma=\eps_{i,c}$ it is the identity. Figure~\ref{fig:SiEpsCoc} 
shows that
$$   
s_{j,c}(\sigma_{k,c}) = \left\{\begin{array}{ll}   
1,& j\ne k-1, \\   
g_{k-1,c}^{-1},& j=k-1,   
\end{array}\right.   
\qquad   
s_{j,c}(\eps_{i,c}) = \left\{\begin{array}{ll}   
1,& j\ne i, \\   
g_{i,c},& j=i.
\end{array}\right.   
$$   
\end{remark}   
   
\begin{lemma}\label{le:s}   
The functions $s_{j,c}:\tcB\to\Gamma$ defined by~\eqref{eq:sjc}   
satisfy the conjugation condition   
\begin{equation}\label{eq:sjc-con}   
\sigma_*g_{k,c} = s_{k,c}(\sigma)^{-1}g_{i,c}s_{k,c}(\sigma),\qquad   
i:=\pi_{\sigma,c}(k),   
\end{equation}   
the cocycle condition   
\begin{equation}\label{eq:sjc-coc}   
s_{k,c}(\sigma\tau)    
= s_{j,c}(\sigma) \sigma_*s_{k,c}(\tau),\qquad   
j:=\pi_{\tau,c}(k),   
\end{equation}   
and the coboundary condition   
\begin{equation}\label{eq:sjc-cob}   
s_{k,c}(\sigma\tau)   
= s_{\ell,c}(\tau)s_{k,\tau_*c}(\sigma),\qquad   
\ell:=\pi_{\sigma,\tau_*c}(k)   
\end{equation}   
for $\sigma,\tau\in\tcB$ and $k=1,\dots,m$.   
\end{lemma}     
   
\begin{proof}   
To prove~\eqref{eq:sjc-con} we denote $i:=\pi_{\sigma,c}(k)$.   
Then   
$$   
\sigma_*g_{k,c} = g_{k,\sigma_*c}    
= (\sigma_*c_k)^{-1}\cdot c_i\cdot g_{i,c}\cdot c_i^{-1}\cdot \sigma_*c_k,   
$$   
where the middle term $c_i\cdot g_{i,c}\cdot c_i^{-1}$ represents a    
counterclockwise turn about~$z_i$. To prove~\eqref{eq:sjc-coc}   
we denote $j:=\pi_{\tau,c}(k)$ and $i:=\pi_{\sigma,c}(j)=\pi_{\sigma\tau,c}(k)$.   
Then    
\begin{eqnarray*}   
s_{j,c}(\sigma)\sigma_*s_{k,c}(\tau)   
&=&   
c_i^{-1}\cdot\sigma_*c_j\cdot \sigma_*(c_j^{-1}\cdot\tau_*c_k) \\   
&=&   
c_i^{-1}\cdot\sigma_*\tau_*c_k \\   
&=&   
s_{k,c}(\sigma\tau).   
\end{eqnarray*}   
To prove~\eqref{eq:sjc-cob} we denote $\ell:=\pi_{\sigma,\tau_*c}(k)$   
and ${i:=\pi_{\tau,c}(\ell)=\pi_{\sigma\tau,c}(k)}$.  Then    
\begin{eqnarray*}   
s_{\ell,c}(\tau)s_{k,\tau_*c}(\sigma)   
&=&   
c_i^{-1}\cdot\tau_*c_\ell\cdot (\tau_*c)_\ell^{-1}\cdot\sigma_*(\tau_*c)_k \\   
&=&   
c_i^{-1}\cdot\sigma_*\tau_*c_k \\   
&=&   
s_{k,c}(\sigma\tau).   
\end{eqnarray*}   
This proves the lemma.   
\end{proof}   

\begin{remark}\rm
It follows from the definition of the mapping class group 
$\tcG/\tcG_0 \cong \tcB$ that $\sigma_*g_{k,c}$ is conjugate 
to $g_{i,c}$ for some $i$. (When $\Sigma$ is the disc, 
this condition appears in Artin's theorem: see Remark~\ref{rmk:artin}.)
Thus Lemma~\ref{le:s} gives an explicit formula for a 
conjugating group element, namely $s_{k,c}(\sigma)$. 
\end{remark}

For every marked distinguished configuration $c=(c_1,\dots,c_m)$   
we define the map $\cS_c:\tcB\to\GL_m(\Z[\Gamma])$ by    
\begin{equation}\label{eq:Scij}   
(\cS_c(\sigma))_{ij} := \left\{   
\begin{array}{ll}   
s_{j,c}(\sigma),&\mbox{if }i=\pi_{\sigma,c}(j), \\   
0,&\mbox{if }i\ne\pi_{\sigma,c}(j).   
\end{array}\right.   
\end{equation}   

\begin{remark}\label{rmk:Sc}\rm   
By Remark~\ref{rmk:sc} we have    
\begin{equation*}   
{\scriptsize
\cS_c(\sigma_{k,c}):=    
\left(\begin{array}{cccc}     
\one_{k-2} &     0         & 0 &    0       \\    
   0       &     0         & 1 &    0       \\    
   0       &  g_{k-1,c}^{-1} & 0 &    0       \\    
   0       &     0         & 0 & \one_{m-k}     
\end{array}\right), \    
\cS_c(\eps_{i,c}) :=    
\left(\begin{array}{ccc}     
\one_{i-1} &  0  & 0 \\    
   0       & g_{i,c} & 0 \\    
   0       &  0  & \one_{m-i}    
\end{array}\right).   
}
\end{equation*}    
When $\Sigma$ is the disc, these conditions uniquely   
determine the cocycle $\cS_c$.   
\end{remark}   

The next theorem contains the statement of Theorem~A. 
The notion of cocycle has been introduced in
Definition~\ref{def:cocycle}. 

\begin{theorem}\label{thm:Sc}   
The maps $\cS_c$ with $c\in\tcC$ satisfy the following conditions.   
   
\smallskip\noindent{\bf (Homotopy)}    
If $c$ is homotopic to $c'$ then $\cS_c=\cS_{c'}$.  

\smallskip\noindent{\bf (Injectivity)}       
Each map $\cS_c:\tcB\to \GL_m(\Z[\Gamma])$ is injective. 
   
\smallskip\noindent{\bf (Cocycle)}    
For all $c\in\tcC$ and $\sigma,\tau\in\tcB$ we have    
\begin{equation}\label{eq:Coc}    
\cS_c(\sigma\tau) = \cS_c(\sigma)\sigma_*\cS_c(\tau).    
\end{equation}    
   
\smallskip\noindent{\bf (Coboundary)}    
For all $c\in\tcC$ and $\sigma,\tau\in\tcB$ we have    
\begin{equation}\label{eq:Cob}    
\cS_{\tau_*c}(\sigma) = \cS_c(\tau)^{-1}\cS_c(\sigma)\sigma_*\cS_c(\tau).    
\end{equation}   
   
\smallskip\noindent{\bf (Monodromy)}    
For each $c\in\tcC$ the formula    
\begin{equation}\label{eq:sigma*N}   
\sigma^*\cN := \bigl(\cS_c(\sigma)^t\cN\cS_c(\sigma)\bigr)\circ\sigma_*    
\end{equation}
defines a contravariant group action of $\tcB$ on $\sN_c$.   
For $c\in\tcC$ and $\tau\in\tcB$ the formula   
\begin{equation}\label{eq:Ttauc}    
\sT_{\tau,c}(\cN) := \cS_c(\tau)^t\cN\cS_c(\tau)   
\end{equation}    
defines an equivariant isomorphism from $\sN_c$ to $\sN_{\tau_*c}$.    
The isomorphisms satisfy the composition rule     
$\sT_{\sigma,\tau_*c}\circ\sT_{\tau,c}=\sT_{\sigma\tau,c}$.    

\smallskip\noindent{\bf (Representation)}    
For all $c\in\tcC$ and $\tau\in\tcB$ the matrix    
$\cS_c(\tau)\in\GL_m(\Z[\Gamma])$ induces    
a collection of isomorphisms    
$$    
\cS_c(\tau):\H_{\sT_{\tau,c}(\cN)}\longrightarrow\H_\cN,\qquad     
\cN\in\sN_c,    
$$    
that preserve the structures~(\ref{eq:quadratic}-\ref{eq:Li}).     
   
\smallskip \noindent{\bf (Lefschetz)}    
If $X\to\Sigma$ is a Lefschetz fibration with singular fibers     
over $Z$ then    
$$    
\sT_{\tau,c}(\cN^X_c) = \cN^X_{\tau_*c}    
$$    
for all $c\in\tcC$ and $\tau\in\tcB$.  Moreover,    
if $\Phi^X_c:\H_{\cN^X_c}\to V/W$ denotes the isomorphism    
of Remark~\ref{rmk:monodromy}, then    
$   
\Phi^X_{\tau_*c}  = \Phi^X_c\circ\cS_c(\tau):   
\H_{\cN^X_{\tau_*c}}\to V/W.   
$   
    
\smallskip    
\noindent{\bf (Odd)}    
If $n$ is odd then the contravariant action of $\tcB$     
on $\sN_c$ descends to $\cB$.   
\end{theorem}   
 
\begin{proof}   
The {\it (Homotopy)} condition is obvious.
To prove the {\it (Cocycle)} condition~\eqref{eq:Coc}    
we denote $j:=\pi_{\tau,c}(k)$ and    
$i:=\pi_{\sigma,c}(j)=\pi_{\sigma\tau,c}(k)$.   
Then   
\begin{eqnarray*}   
(\cS_c(\sigma)\sigma_*\cS_c(\tau))_{ik}   
&=&    
\sum_{\nu=1}^m (\cS_c(\sigma))_{i\nu}\sigma_*((\cS_c(\tau))_{\nu k}) \\   
&=&    
(\cS_c(\sigma))_{ij}\sigma_*(\cS_c(\tau))_{jk} \\   
&=&    
s_{j,c}(\sigma)\sigma_*s_{k,c}(\tau) \\   
&=&    
s_{k,c}(\sigma\tau) \\   
&=&    
(\cS_c(\sigma\tau))_{ik}.   
\end{eqnarray*}   
Here we have used~\eqref{eq:sjc-coc} and~\eqref{eq:Scij}.   
For $i\ne\pi_{\sigma\tau,c}(k)$ the $(i,k)$ entry of both matrices    
$\cS_c(\sigma\tau)$ and $\cS_c(\sigma)\sigma_*\cS_c(\tau)$ is zero.  
Thus we have proved~\eqref{eq:Coc}.   
   
To prove the {\it (Coboundary)} condition~\eqref{eq:Cob}   
we denote $\ell:=\pi_{\sigma,\tau_*c}(k)$   
and ${i:=\pi_{\tau,c}(\ell)=\pi_{\sigma\tau,c}(k)}$.  Then    
\begin{eqnarray*}   
(\cS_c(\tau)\cS_{\tau_*c}(\sigma))_{ik}   
&=&    
\sum_{\nu=1}^m (\cS_c(\tau))_{i\nu}(\cS_{\tau_*c}(\sigma))_{\nu k} \\   
&=&    
(\cS_c(\tau))_{i\ell}(\cS_{\tau_*c}(\sigma))_{\ell k} \\   
&=&    
s_{\ell,c}(\tau)s_{k,\tau_*c}(\sigma) \\   
&=&    
s_{k,c}(\sigma\tau) \\   
&=&    
(\cS_c(\sigma\tau))_{j\ell}.   
\end{eqnarray*}   
Here we have used~\eqref{eq:sjc-cob} and~\eqref{eq:Scij}.   
For $i\ne\pi_{\sigma\tau,c}(k)$ the $(i,k)$ entry of both matrices    
$\cS_c(\sigma\tau)$ and $\cS_c(\tau)\cS_{\tau_*c}(\sigma)$  is zero.     
Thus we have deduced~\eqref{eq:Cob} from~\eqref{eq:Coc}.    

To prove the {\it (Injectivity)} condition we assume that 
$c\in\tcC$ is a marked distinguished configuration and 
$\sigma\in\tcB$ is a framed braid such that $\cS_c(\sigma)=\one$. 
Then the permutation $\pi_{\sigma,c}$ is the identity and 
we deduce from~\eqref{eq:sjc-con} that
$$
g_{i,\sigma_*c}
= \sigma_*g_{i,c} 
= s_{i,c}(\sigma)^{-1} g_{i,c} s_{i,c}(\sigma)
= g_{i,c},\qquad i=1,\dots,m.
$$
Hence it follows from Theorem~\ref{thm:C} that $c$
and $\sigma_*c$ descend to equivalent distinguished
configurations in $\cC$.  By Theorem~\ref{thm:BC} 
this implies that $\sigma$ is a lift of the trivial braid in $\cB$.  Thus 
$$
\sigma = \eps_{1,c}^{k_1}\cdots\eps_{m,c}^{k_m}
$$
for some integer vector $(k_1,\dots,k_m)\in\Z^m$.
Using again the fact that $\cS_c(\sigma)=\one$ we obtain
that $k_1=\cdots=k_m=0$ and hence $\sigma=1$. 
Thus we have proved that, for every $c\in\tcC$
and every $\sigma\in\tcB$, we have 
$$
\cS_c(\sigma)=\one\qquad\implies\qquad\sigma=1.
$$
Now let $c\in\tcC$ and $\sigma,\tau\in\tcB$ be given such that
$$
\cS_c(\sigma\tau)=\cS_c(\tau).
$$
Then it follows from the coboundary and cocycle conditions
that
$$
\cS_c(\tau)=\cS_c(\sigma\tau)=\cS_c(\sigma)\sigma_*\cS_c(\tau)
= \cS_c(\tau)\cS_{\tau_*c}(\sigma).
$$
Hence $\cS_{\tau_*c}(\sigma)=\one$ and so, by what we have 
already proved, it follows that $\sigma=1$.  This shows that
the map $\cS_c:\tcB\to\GL_m(\Z[\Gamma])$ is injective, as claimed.

To prove the {\it (Monodromy)} condition let $\cN=(n_{ij})\in\sN_c$   
and $\tau\in\tcB$. We must prove that    
$\sT_{\tau,c}(\cN) := \cS_c(\tau)^t\cN\cS_c(\tau)\in\sN_{\tau_*c}$.   
To see this denote the entries of  $\sT_{\tau,c}(\cN)$ by $\tn_{ij}$   
and observe that    
\begin{equation}\label{eq:tn}   
\tn_{ij}(g) = n_{i'j'}\bigl(s_{i,c}(\tau)gs_{j,c}(\tau)^{-1}\bigr),\qquad   
i':=\pi_{\tau,c}(i),\qquad j':=\pi_{\tau,c}(j).   
\end{equation}   
That the functions $\tn_{ij}$ satisfy~\eqref{eq:N1} and~\eqref{eq:N2}   
is obvious from this formula. To prove~\eqref{eq:N3} we abbreviate    
$\eps:=(-1)^{n(n+1)/2}$, $k':=\pi_{\tau,c}(k)$, and compute   
\begin{eqnarray*}   
\tn_{ij}(gg_{k,\tau_*c}h)   
&=&     
n_{i'j'}\bigl(s_{i,c}(\tau)g(\tau_*g_{k,c})hs_{j,c}(\tau)^{-1}\bigr)  \\   
&=&     
n_{i'j'}\bigl(s_{i,c}(\tau)g   
s_{k,c}(\tau)^{-1}g_{k',c}s_{k,c}(\tau)   
hs_{j,c}(\tau)^{-1}\bigr)  \\   
&=&     
n_{i'j'}\bigl(s_{i,c}(\tau)ghs_{j,c}(\tau)^{-1}\bigr)  \\   
&& - \eps   
n_{i'k'}\bigl(s_{i,c}(\tau)gs_{k,c}(\tau)^{-1}\bigr)   
n_{k'j'}\bigl( s_{k,c}(\tau)hs_{j,c}(\tau)^{-1}\bigr)  \\   
&=&   
\tn_{ij}(gh)-\eps\tn_{ik}(g)\tn_{kj}(h).   
\end{eqnarray*}   
Here the first and last equations follow from~\eqref{eq:tn},    
the second equation follows from~\eqref{eq:sjc-con},    
and the third follows from~\eqref{eq:N3} for $n_{i'j'}$.    
Thus we have proved that $\sT_{\tau,c}(\cN)\in\sN_{\tau_*c}$, as claimed,    
and hence also    
$$   
\tau^*\cN=\sT_{\tau,c}(\cN)\circ\tau_*\in\sN_c.   
$$    
That equation~\eqref{eq:sigma*N} defines a contravariant   
group action of $\tcB$ on $\sN_c$ for every $c\in\tcC$ 
is a consequence of Lemma \ref{le:cocycle}.
That the map $\sT_{\tau,c}:\sN_c\to\sN_{\tau_*c}$ 
is equivariant under the action of $\tcB$ follows from~\eqref{eq:Cob}.     
The composition rule follows from the definitions as well    
as~\eqref{eq:Coc} and~\eqref{eq:Cob}.    
This proves the {\it (Monodromy)} condition.   
   
The {\it (Representation)} condition follows from    
the proof of Proposition~\ref{prop:HN}.   
To prove the {\it (Lefschetz)} condition we fix a symplectic   
Lefschetz fibration $f:X\to\Sigma$ with critical fibers over $Z$   
as well as a marked distinguished configuration $c\in\tcC$    
and a framed braid $\tau\in\tcB$.  To emphasize 
the dependence of the vanishing cycles on the choice 
of distinguished configuration $c$, we denote their 
homology classes by $L_{1,c},\dots,L_{m,c}$. 
Let $\rho:\Gamma\to\Aut(H_n(M))$ be the associated    
monodromy representation and, for $g\in\Gamma$, let   
$$   
n_{ij}(g) := \inner{L_{i,c}}{\rho(g)L_{j,c}}   
$$   
denote the entries of the intersection matrix $\cN^X_c(g)$.   
Then    
$$   
L_{i,\tau_*c} = \rho(s_{i,c}(\tau)^{-1})L_{i',c},\qquad   
i':=\pi_{\tau,c}(i).   
$$   
Hence the entries of the matrix    
$\cN^X_{\tau_*c}(g)=:(\tn_{ij}(g))$ are    
\begin{eqnarray*}   
\tn_{ij}(g)    
&=&   
\inner{L_{i,\tau_*c}}{\rho(g)L_{j,\tau_*c}} \\   
&=&   
\inner{\rho\bigl(s_{i,c}(\tau)^{-1}\bigr)L_{i',c}}   
{\rho\bigl(gs_{j,c}(\tau)^{-1}\bigr)L_{j',c}} \\   
&=&   
\inner{L_{i',c}}{\rho\bigl(s_{i,c}(\tau)g   
s_{j,c}(\tau)^{-1}\bigr)L_{j',c}} \\   
&=&   
n_{i'j'}\bigl(s_{i,c}(\tau)gs_{j,c}(\tau)^{-1}\bigr).   
\end{eqnarray*}    
Hence it follows from~\eqref{eq:tn} that    
$\cN^X_{\tau_*c}=\sT_{\tau,c}(\cN^X_c)$ as claimed.   
Thus we have proved the {\it (Lefschetz)} condition.   
   
Now assume that $n$ is odd. Then it follows from~\eqref{eq:N5} that   
\begin{equation}\label{eq:odd}   
n_{ij}(g)=n_{ij}(g_{i,c}g)=n_{ij}(gg_{j,c})   
\end{equation}   
for $\cN=(n_{ij})\in\sN_c$.  Moreover, if $\tau\in\tcB$ belongs 
to the kernel of the homomorphism $\tcB\to\cB$ then 
$\pi_{\tau,c}=\id\in\gS_m$ and $s_{i,c}(\tau)=g_{i,c}^{k_i}$ 
for some $k_i\in\Z$. Hence equations~\eqref{eq:tn} 
and~\eqref{eq:odd} show that any such element $\tau$ 
acts trivially on $\sN_c$.  This proves the theorem.    
\end{proof}


 \section{Comparison with the Magnus cocycle}\label{sec:magnus}
 
Let $\Gamma$ be a free group of finite rank and 
denote by $\Aut(\Gamma)$ its group of automorphisms.
Any choice of a basis $g_1,\dots,g_m\in\Gamma$ determines a 
{\bf Magnus cocycle} 
$$
\cM:\Aut(\Gamma)\to\GL_m(\Z[\Gamma])
$$
defined by 
\begin{equation}\label{eq:magnus}
\cM(\psi) := \left(\overline{\frac{\p \psi(g_j)}{\p g_i}}
\right)_{i,j=1,\dots,m}.
\end{equation}
This formula is to be understood as follows.
A {\bf derivation} is a $1$-cocycle $d:\Gamma\to\Z[\Gamma]$, 
i.e$.$ a map that satisfies the equation
$$
d(gh) = d(g)+gd(h)
$$
for all $g,h\in\Gamma$.  In particular we have $d(1)=0$ 
and $d(g^{-1})=-g^{-1}d(g)$ for every $g\in\Gamma$. 
Examples of derivations are the {\bf Fox derivatives}~\cite{FOX}
$$
\frac{\p}{\p g_i}:\Gamma\to\Z[\Gamma],\qquad i=1,\dots,m,
$$
characterized by the condition
$$
\frac{\p g_j}{\p g_i} = \delta_i^j,\qquad i,j=1,\dots,m.
$$
Recall that the conjugation 
$\Z[\Gamma]\to\Z[\Gamma]:\lambda\mapsto\bar\lambda$ 
is the ring anti-ho\-mo\-mor\-phism defined by  
$\bar g:=g^{-1}$ for all $g\in\Gamma\subset\Z[\Gamma]$.
This explains the right hand side of~\eqref{eq:magnus}.

The map
$
\cM:\Aut(\Gamma)\to\GL_m(\Z[\Gamma])
$
satisfies the cocycle condition
\begin{equation}\label{eq:cocmagnus}
\cM(\psi\circ\phi) = \cM(\psi)\cdot \psi_*\cM(\phi)
\end{equation}
for $\phi,\psi\in\Aut(\Gamma)$.  
This is proved in Birman~\cite{BIRMAN} as a 
consequence of the chain rule for Fox calculus.   
Fox calculus has its origin in the theory of covering spaces. 
The matrix $\cM(\phi)$ represents the action of $\phi$ on the twisted 
homology of a bouquet of $m$ circles relative to a base point
with coefficients in $\Z[\Gamma]$. The resulting map is a cocycle 
(instead of a homomorphism) because the lift of a continuous map 
of the bouquet of circles to its universal cover is not 
$\Gamma$-equivariant.

This construction applies to the braid group of the disc as follows. 
We return to the geometric setting of Section~\ref{sec:braid}
with $\Sigma=\D$.  Thus $Z\subset\D$ is a set of $m$ points in the interior,
$\Gamma$ is the fundamental group of $\D\setminus Z$ based at $z_0\in\p\D$, 
and $\cB$ is the braid group on $m$ strings in $\D$ based at $Z$. 
The choice of a distinguished configuration $c\in\cC$ determines
a basis $g_{1,c},\dots,g_{m,c}$ of $\Gamma$.
Since $\cB$ acts on $\Gamma$ we obtain a  Magnus cocycle 
$$
\cM_c:\cB\to\GL_m(\Z[\Gamma])
$$
for every distinguished configuration $c\in\cC$. 

\begin{proposition}\label{prop:magnus}
The maps $\cM_c:\cB\to\GL_m(\Z[\Gamma])$ are injective 
and satisfy the cocycle and coboundary conditions
\begin{equation}\label{eq:magnusBC}
\cM_c(\sigma\tau) 
= \cM_c(\sigma)\sigma_*\cM_c(\tau),\qquad  
\cM_{\tau_*c}(\sigma) 
= \cM_c(\tau)^{-1}\cM_c(\sigma)\sigma_*\cM_c(\tau)
\end{equation}
for all $c\in\cC$ and $\sigma,\tau\in\cB$.
\end{proposition}

\begin{proof}
The first equation in~\eqref{eq:magnusBC} follows immediately
from~\eqref{eq:cocmagnus} and the second equation can also be derived 
from the chain rule in Fox calculus. For the sake of completeness we give 
the details.  The chain rule in Fox calculus has the following form.
If $g_1,\dots,g_m$ and $h_1,\dots,h_m$ are two basis of $\Gamma$
and $a\in\Gamma$ is an arbitrary element then
$$
\frac{\p a}{\p g_i} = \sum_{j=1}^m \frac{\p a}{\p h_j}\frac{\p h_j}{\p g_i}, \qquad i=1,\dots,m.
$$
To prove the first formula in~\eqref{eq:magnusBC} we take
$$
g_i:=g_{i,c},\qquad h_j := \sigma_*g_{j,c},\qquad a:=\sigma_*\tau_*g_{k,c}.
$$
Then the chain rule asserts that
$$
\frac{\p (\sigma\tau)_*g_{k,c}}{\p g_{i,c}}
= \sum_{j=1}^m \frac{\p \sigma_*\tau_*g_{k,c}}{\p \sigma_*g_{j,c}}
\frac{\p \sigma_*g_{j,c}}{\p g_{i,c}}
= \sum_{j=1}^m \left(\sigma_*\frac{\p \tau_*g_{k,c}}{\p g_{j,c}}\right)
\frac{\p \sigma_*g_{j,c}}{\p g_{i,c}}.
$$
The first equation in~\eqref{eq:magnusBC} follows by conjugation.  
To prove the second equation in~\eqref{eq:magnusBC} we choose 
$$
g_i:=g_{i,c},\qquad h_j := \tau_*g_{j,c},\qquad a:=\sigma_*\tau_*g_{k,c}.
$$
Then the chain rule asserts that
$$
\frac{\p (\sigma\tau)_*g_{k,c}}{\p g_{i,c}}
= \sum_{j=1}^m \frac{\p \sigma_*\tau_*g_{k,c}}{\p \tau_*g_{j,c}}
\frac{\p \tau_*g_{j,c}}{\p g_{i,c}}
= \sum_{j=1}^m \frac{\p \sigma_*g_{k,\tau_*c}}{\p g_{j,\tau_*c}}
\frac{\p \tau_*g_{j,c}}{\p g_{i,c}}.
$$
Hence conjugation gives $\cM_c(\sigma\tau)=\cM_c(\tau)\cM_{\tau_*c}(\sigma)$ 
and so the second equation in~\eqref{eq:magnusBC} follows from the first.  

\smallbreak

Injectivity of $\cM_c$ is a consequence of the fundamental formula
in Fox calculus~\cite{FOX}.  It has the form
\begin{equation}\label{eq:Fox}
a-1 = \sum_{i=1}^m \frac{\p a}{\p g_{i,c}}\left(g_{i,c}-1\right)
\end{equation}
for all $a\in\Gamma$. Applying this formula to $a=\sigma_*g_{j,c}$ 
and $a=\tau_*g_{j,c}$ we see that $\cM_c(\sigma)=\cM_c(\tau)$ 
if and only if $\sigma_*g=\tau_*g$ for all $g\in\Gamma$, 
which is equivalent to $\sigma=\tau$ by Artin's theorem 
(see Remark~\ref{rmk:artin}). 
This proves the proposition. 
\end{proof}

The Magnus cocycle $\cM_c$ is connected to the Reidemeister 
intersection pairing (in its relative version)
$$
\langle\cdot,\cdot\rangle:
H_1(\D\setminus Z,z_0;\Z[\Gamma])\times 
H_1(\D\setminus Z,z_0;\Z[\Gamma])\to \Z[\Gamma]
$$
or, equivalently, to the homotopy intersection pairing 
$$
\omega:\Gamma\times \Gamma\longrightarrow \Z[\Gamma] 
$$ 
introduced by Turaev~\cite{Turaev} and Perron~\cite{Perron}. 
Let $\Omega_c$ be the $m\times m$-matrix with 
coefficients in $\Z[\Gamma]$ which represents 
$\omega$ in the basis $(g_{1,c},\dots,g_{m,c})$. 
Perron shows that 
\begin{equation}\label{eq:Perron}
\cM_c(\tau)^t \cdot \Omega_c \cdot \cM_c(\tau)=\Omega_{\tau_*c}
\end{equation}
for any $c\in\cC$ and $\tau\in\cB$. 
This identity can be deduced from the topological interpretation 
of the Magnus cocycle, according to which $\cM_c(\tau)$ is
(after conjugation) the matrix representing the homomorphism
$$
\tau_*:H_1(\D\setminus Z,z_0;\Z[\Gamma])
\to H_1(\D\setminus Z,z_0;\Z[\Gamma])
$$
with respect to a  basis given by lifts of $g_{1,c},\dots,g_{m,c}$. 

\begin{remark}\label{rmk:Mc}\rm  
By the definition of Fox derivatives we have 
$$
\frac {\p (g_{i-1,c}g_{i,c}g_{i-1,c}^{-1})} {\p g_{i-1,c}} 
= 1-g_{i-1,c}g_{i,c}g_{i-1,c}^{-1},\qquad
\frac {\p (g_{i-1,c}g_{i,c}g_{i-1,c}^{-1})} {\p g_{i,c}} 
= g_{i-1,c}.
$$
Hence the Magnus cocycle satisfies    
\begin{equation*}   
{\scriptsize
\cM_c(\sigma_{i,c}):=    
\left(\begin{array}{cccc}     
\one_{i-2} & 0  & 0 & 0 \\    
   0 & 1-g_{i-1,c}g_{i,c}^{-1}g_{i-1,c}^{-1} & 1 & 0 \\    
   0 & g_{i-1,c}^{-1} & 0 &  0  \\    
   0 & 0 & 0 & \one_{m-i}     
\end{array}\right)
}
\end{equation*}
for all $i=2,\dots,m$.
\end{remark}   

In order to compare the Magnus cocycle with the Picard--Lefschetz cocycle,
we need to restrict the latter to an embedded image of $\cB$ in $\tcB$.
For this we fix, as in the previous sections, a  non-zero tangent vector 
$v_z$ at each $z\in Z$.  Moreover we fix a contractible neighborhood $U$ 
of $z_0$ in $\Sigma = \D$ whose closure does not meet $Z$,
and we fix a vector field $v$ on $U$ whose only singularity 
is an attractive point at $z_0$. Relative homotopy classes 
of nowhere vanishing vector fields on $\Sigma$ that restrict 
to $v$ on $Z\cup U$ are parametrized by 
$H^1(\Sigma,U\cup Z;\Z)\cong \Z^m$.
We {\bf choose} such a relative homotopy class $\xi$. 
The class $\xi$ defines a section
$$
\cB\stackrel\xi\hookrightarrow \tcB
$$
of the short exact sequence~\eqref{eq:exact},
which assigns to every braid the framing defined 
by any representative $w$ of $\xi$.  Furthermore, 
we can associate to each $\gamma \in \cC$
the homotopy class $[c] \in \tcC$ where $c$ is a 
distinguished configuration representing $\gamma$ 
which is tangent to $w$ for some representative $w$ of $\xi$. 
This defines a section 
$$
\cC\stackrel\xi\hookrightarrow \tcC
$$
of the canonical projection $\tcC \to \cC$.
The free and transitive action of $\tcB$ on $\tcC$ restricts
to a free and transitive action of $\xi(\cB)$ on $\xi(\cC)$.
In the sequel, we shall identify $\cB$ and $\cC$ with their 
respective images by $\xi$.

For every $c\in\cC$, the cocycles $\cM_c$ show similarities 
with the cocycles $\cS_c|_\cB$. Both of them are injective 
maps $\cB\to\GL_m(\Z[\Gamma])$ and, according to the 
formulas~\eqref{eq:cocb} and~\eqref{eq:magnusBC}, 
they behave in the same way under change of $c\in\cC$.
Moreover,  the formulas~\eqref{eq:Tpsic} and~\eqref{eq:Perron} 
show that they can both be interpreted as matrices of basis change 
for certain geometrically defined bilinear forms.
However, in contrast to $\cM_c$, the entries 
of $\cS_c$ are elements of $\Gamma\subset\Z[\Gamma]$.
The framed braids $\sigma_{2,c},\dots,\sigma_{m,c}$ 
belong to the embedded image of $\cB$ in $\tcB$ defined 
by the homotopy class of vector fields $\xi$.
We see from the formulas in Remarks~\ref{rmk:Sc} and~\ref{rmk:Mc} 
that $\cM_c(\sigma_{i,c})$ differs from $\cS_c(\sigma_{i,c})$ 
exactly in one entry.

The Magnus cocycle gives a unified framework for the definition 
of various linear representations of mapping class groups, 
including the Burau and Gassner representations~\cite{BIRMAN}. 
These representations are obtained by reducing $\Gamma$ 
to an abelian quotient. Let us apply the same 
reductions to the Picard--Lefschetz cocycle.

There is a natural homomorphism $\Gamma\to\Z$ which assigns 
to every loop $g$ in $\D\setminus Z$ the total winding number 
around the punctures. If we identify $\Z$ with the free group 
on one generator $t^{-1}$, this homomorphism sends $g_{i,c}$ to $t^{-1}$. 
There is an induced ring homomorphism $\Z[\Gamma]\to\Z[t,t^{-1}]$, 
and hence a group homomorphism $\GL_m(\Z[\Gamma])\to\GL_m(\Z[t,t^{-1}])$. 
Composition with this homomorphism turns every cocycle into 
a representation, and turns cohomologous cocycles into 
conjugate representations. The compositions of $\cS_c|_\cB$ 
and $\cM_c$ with the homomorphism 
$\GL_m(\Z[\Gamma])\to\GL_m(\Z[t,t^{-1}])$ 
will be denoted by 
$$
\overline\cS_c:\cB\to\GL_m(\Z[t,t^{-1}]),\qquad 
\overline\cM_c:\cB\to\GL_m(\Z[t,t^{-1}]).
$$
On the generators $\sigma_{2,c},\dots,\sigma_{m,c}$ of $\cB$ we have 
\begin{equation*}   
{\scriptsize
\overline\cS_c(\sigma_{i,c}):=    
\left(\begin{array}{cccc}     
\one_{i-2} &     0         & 0 &    0       \\    
   0       &     0    & 1 &    0       \\    
   0       &  t & 0 &    0       \\    
   0       &     0         & 0 & \one_{m-i}     
\end{array}\right),\quad
\overline\cM_c(\sigma_{i,c}):=    
\left(\begin{array}{cccc}     
\one_{i-2} &     0         & 0 &    0       \\    
   0       &     1-t   & 1 &    0       \\  
   0       &  t & 0 &    0       \\    
   0       &     0         & 0 & \one_{m-i}     
\end{array}\right).
}
\end{equation*}

These representations of the braid group are well known: 
$\overline\cM_c$ is the {\bf Burau representation}~\cite{Burau,BIRMAN,KT}
and $\overline\cS_c$ is the {\bf Tong--Yang--Ma representation} 
introduced in~\cite{TYM}. The Burau representation plays 
an important role in knot theory because of its deep connection 
with the Alexander polynomial~\cite{BIRMAN,KT}. In contrast 
to the Burau representation, the Tong--Yang--Ma representation 
is irreducible. Furthermore, Sysoeva shows in~\cite{Sysoeva} 
that any irreducible $m$-dimensional complex representation 
of the braid group on $m\ge 9$ strings is equivalent to the 
tensor product of a $1$-dimensional representation with a 
specialization of the latter for some $t\in\C\setminus\{0,1\}$. 
(See also \cite{FLSV} for the cases $m\in \{5,6,7,8\}$.)

\begin{proposition} \label{prop:McSc}
The cocycles $\cM_c$ and $\cS_c|_\cB$ define distinct and nontrivial 
cohomology classes in $H^1(\cB,\GL_m(\Z[\Gamma]))$. 
\end{proposition}

\begin{proof}
We have $\tr(\overline\cS_c(\sigma_{i,c}))=m-2$ and 
$\tr(\overline\cM_c(\sigma_{i,c}))=m-1-t$ for $i=2,\dots,m$, 
while $\tr(\mathrm{Id})=m$.  Thus $\overline\cS_c$ and $\overline\cM_c$
are neither conjugate to each other nor conjugate to the 
trivial representation. 
\end{proof}

\begin{remark}\label{rmk:TYMextended} \rm
Since the cocycle $\cS_c$ is defined on $\tcB$, 
it gives rise by reduction to an extension 
$$
\overline\cS_c :\tcB\to \GL_m(\Z[t,t^{-1}])
$$
of the Tong--Yang--Ma representation to the 
framed braid group. Explicitly, we have 
\begin{equation*}   
{\scriptsize
\overline\cS_c(\eps_{i,c}):=    
\left(\begin{array}{cccc}     
\one_{i-1} &     0         & 0       \\    
   0       &     t^{-1}    &    0       \\    
   0       &     0         &  \one_{m-i}     
\end{array}\right), \qquad i=1,\dots,m.
}
\end{equation*}  
This formula implies in particular that the 
Picard--Lefschetz monodromy class in 
$H^1(\tcB,\GL_m(\Z[\Gamma]))$ is nontrivial.
\end{remark}

Fix an ordering $Z=\{z_1,\dots,z_m\}$ of the punctures.
A {\bf pure braid} is a braid $\sigma\in\cB$ 
whose associated permutation of $Z$ is trivial.
We denote by $\cP\cB \subset \cB$ the subgroup of pure braids.
The ordering of $Z$ induces a natural isomorphism 
between the abelianization of $\Gamma$ and $\Z^m$
and hence a natural homomorphism from $\Gamma$ to $\Z^m$. 
If we identify $\Z^m$ with the free abelian group on $m$ generators 
$t_1^{-1},\dots,t_m^{-1}$, this homomorphism sends $g_{i,c}$ to $t_i^{-1}$ 
for every distinguished configuration $c$ that
determines the given ordering of~$Z$.
Thus, there is an induced group homomorphism 
$\GL_m(\Z[\Gamma])\to\GL_m(\Z[t_1^{\pm1},\dots,t_m^{\pm1}])$. 
The compositions of $\cS_c|_{\cP\cB}$ and $\cM_c|_{\cP\cB}$ with this
homomorphism will be denoted by 
$$
\widehat\cS_c : \cP\cB\to\GL_m(\Z[t_1^{\pm1},\dots,t_m^{\pm1}]), \qquad
\widehat\cM_c : \cP\cB\to\GL_m(\Z[t_1^{\pm1},\dots,t_m^{\pm1}]).
$$
These are representations of the pure braid group and $\widehat\cM_c$ 
is called the {\bf Gassner representation}~\cite{Gassner,BIRMAN}.
Observe that the embedding of the subgroup $\cP\cB$ in $\tcB$ is canonical
(i.e.\ it does not depend on the choice of $\xi$)
and, furthermore, the representation $\widehat\cS_c$ does not 
depend on $c$ but only on the chosen ordering of $Z$. 
Indeed, using equation~\eqref{eq:Cob} and the fact that 
$\widehat\cS_c(\tau)$ is a diagonal matrix for all $\tau\in\cP\cB$, 
we obtain  
$$
\widehat\cS_{\tau_*c}(\sigma)
=\widehat\cS_c(\tau)^{-1}\widehat\cS_c(\sigma)\widehat\cS_c(\tau)
=\widehat\cS_c(\sigma)
$$
for all $\sigma\in\cP\cB$. 

The next proposition gives an explicit formula 
for $\widehat\cS_c$ and shows that this representation
is completely determined by the linking numbers.

\begin{proposition}\label{prop:linking}
For every $\sigma\in\cP\cB$ we have
\begin{equation}\label{eq:linking}
\widehat\cS_c(\sigma) = {\mathrm{Diag}}\left(
\prod_{j\ne 1}t_j^{-\lk(1,j)},\dots,\prod_{j\ne m}t_j^{-\lk(m,j)}
\right)
\end{equation}
where $\lk(i,j)$ denotes the linking number of the $i$-th and $j$-th 
components of the closure of the braid.
\end{proposition}

\begin{figure}[htp]
\centering    
\includegraphics[scale=0.6]{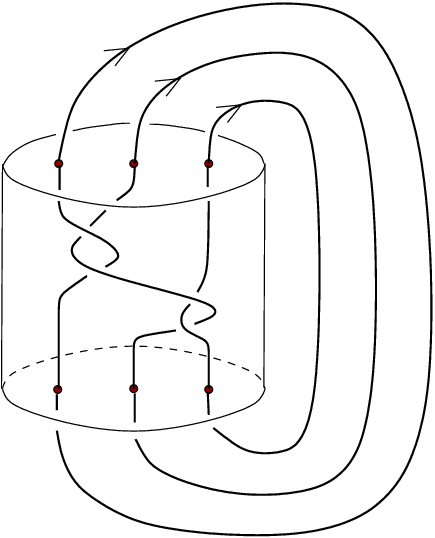} 
\caption{Closure of a pure braid.}\label{fig:closure}   
\end{figure}

\noindent 
Here the {\bf closure} of the braid $\sigma$ is defined 
in the usual way~\cite{BIRMAN,KT} by connecting the top 
and the bottom of $\D\times [0,1]$ without twisting 
(see Figure~\ref{fig:closure} for an illustration).

\begin{proof}[Proof of Proposition~\ref{prop:linking}]
By definition $\cS_c(\sigma)$ is a diagonal matrix with diagonal 
entries 
$
s_{j,c}(\sigma)=c_j^{-1}\cdot \sigma_*c_j\in\Gamma
$
(see equation~\eqref{eq:Scij}).  To understand the corresponding diagonal 
entry $\widehat s_{j,c}(\sigma)$ of $\widehat\cS_c(\sigma)$ we must express 
$s_{j,c}(\sigma)$ as a word in the generators $g_{i,c}$ and their inverses. 
The exponent of $t_i$ is then the total occurence of the factor 
$g_{i,c}^{-1}$ in this word and we claim that this number is 
$-\lk(i,j)$ for $i\ne j$ and is zero for $i=j$. Equivalently, 
if we denote by $\widehat g\in H_1(\D\setminus Z;\Z)$ the 
homology class of an element $g\in\Gamma$, we must prove that 
\begin{equation}
\label{eq:lk}
\widehat{s}_{j,c}(\sigma) = \sum_{i\ne j}\lk(i,j) \widehat{g}_{i,c}.
\end{equation}
To see this, let $N\subset\D\times[0,1]$ be the complement of the braid $\sigma$, 
viewed as a collection of strings running from $\D\times\{0\}$ 
to $\D\times\{1\}$. There is a deformation retract $r:N\to\D\setminus Z$ 
such that, for every $z\in\D\setminus Z$, we have $r(z,0)=z$ and 
$r(z,1)=\phi^{-1}(z)$ with $\phi\in\cG$ a diffeomorphism 
representing the element in the mapping class group 
corresponding to the braid $\sigma$. This map induces 
an isomorphism
$$
r_*:H_1(N;\Z)\to H_1(\D\setminus Z;\Z).
$$
The oriented meridians of the strings of $\sigma$ form a 
basis $\mu_1,\dots,\mu_m$ of $H_1(N;\Z)$ and their images 
under $r_*$ are represented by small loops encircling 
the elements of $Z$ counterclockwise; thus we have
$$
r_*(\mu_i)=\widehat{g}_{i,c},\qquad
i=1,\dots,m. 
$$
Using the distinguished configuration $c$,
we can view the closure of the braid inside $N$. 
We denote its components by $K_1,\dots,K_m$, 
and the corresponding homology classes by 
$\widehat K_1,\dots,\widehat K_m\in H_1(N;\Z)$. 
By the homological definition of the linking numbers, we have 
$$
\widehat K_j=\sum_{i\neq j}\lk(i,j) \mu_i. 
$$
The image by $r$ of the knot $\widehat{K_j}$ is the loop 
that first traverses $c_j|_{[0,1-\varepsilon]}$  
and then $\phi^{-1}\circ c_j|_{[0,1-\varepsilon]}$ 
in the reverse direction (for some small $\varepsilon>0$).
It follows that $r_*(\widehat K_j)=\widehat{s}_{j,c}(\sigma)$ 
for all $j\in\{1,\dots,m\}$ and equation~\eqref{eq:lk} follows. 

Alternatively, one can prove~\eqref{eq:lk} as follows.
This identity is equivalent to the formula
\begin{equation} \label{eq:wlk}
\mathrm{w}(s_{j,c}(\sigma),z_i) = 
\left\{\begin{array}{ll}
\lk(i,j) & \hbox{if } i\ne j\\
0 & \hbox{if } i=j
\end{array}\right.
\end{equation}
where $\mathrm{w}(\gamma,z)$ denotes the winding number of a loop 
$\gamma$ about $z$.  We choose a representative
of the braid $\sigma$ in which $z_i$ is constant. 
Then the $j$-th component of the braid has the same 
winding number about $z_i$ as $s_{j,c}(\sigma)$, 
and this winding number agrees with $\lk(i,j)$ by definition 
of the linking number as an intersection number.
\end{proof}

\begin{remark}\rm 
Proposition~\ref{prop:linking} implies, by specializing 
to $t_1=\dots=t_m=t$, that the Tong--Yang--Ma 
representation $\overline\cS_c$ is trivial on the 
commutator subgroup 
$
\cP\cB':=[\cP\cB,\cP\cB].
$ 
Thus $\overline\cS_c$ factors through the quotient
$\cB/\cP\cB'$, which is an instance of the
{\bf extended Coxeter group} in the sense of Tits~\cite{Tits}. 
In this case the Coxeter group is the symmetric group $\gS_m$, 
and $\cB/\cP\cB'$ is an extension of the latter by the 
abelian group 
$
\cP\cB/\cP\cB'\cong \Z^{m(m-1)/2}.
$
\end{remark}

      
\section{Proof of Theorem~B}\label{sec:proofB}      
    
We still specialize to the case where $\Sigma=\D$ is the closed   
unit disc.  In this case $\Gamma$ is isomorphic to the free group   
$\Gamma_m$ generated by $g_1,\dots,g_m$ and $\tcB$ is isomorphic   
to the abstract framed braid group $\tcB_m$ with generators    
$\sigma_2,\dots,\sigma_m,\eps_1,\dots,\eps_m$ and    
relations~\eqref{eq:frabraid}.  The isomorphisms depend on    
the choice of a marked distinguished configuration $c\in\tcC$   
and will be denoted by    
$$   
\iota_c:\Gamma_m\to\Gamma,\qquad \phi_c:\tcB_m\to\tcB.   
$$   
The isomorphism $\iota_c$ assigns to $g_i$ the special    
element $g_{i,c}$ obtained by encircling $z_i$ counterclockwise   
along $c_i$.  The isomorphism $\phi_c$ assigns to $\sigma_i$ and    
$\eps_i$ the generators $\sigma_{i,c}$ and $\eps_{i,c}$   
of $\tcB$ associated to $c$, as defined in Remark~\ref{rmk:sigmaic}. 
Recall the action of $\tcB_m$ on $\Gamma_m$ by~\eqref{eq:BmGm}.    

\begin{lemma}\label{le:BmGm}   
{\bf (i)} The isomorphisms $\iota_c$ and $\phi_c$ satisfy   
$$   
\iota_{\psi_*c}(g) = \psi_*\iota_c(g),\qquad   
\phi_{\psi_*c}(\sigma) = \psi\phi_c(\sigma)\psi^{-1},   
$$   
for $g\in\Gamma_m$, $\sigma\in\tcB_m$, $c\in\tcC$, and $\psi\in\tcB$.    
   
\smallskip\noindent{\bf (ii)}   
For every $c\in\tcC$ and every $\sigma\in\tcB_m$   
there is a commutative diagram   
$$    
\xymatrix    
@C=20pt    
@R=20pt    
{\Gamma_m \ar[r]^{\sigma_*} \ar[d]_{\iota_c} & \Gamma_m    
  \ar[d]^{\iota_c} \\    
\Gamma \ar[r]_{\phi_c(\sigma)_*} & \Gamma    
}.   
$$    
   
\smallskip\noindent{\bf (iii)}   
The formula   
$$
\sigma^*c := \phi_c(\sigma)_*c   
$$ 
defines a free and transitive contravariant action 
of $\tcB_m$ on $\tcC$ and    
$$ 
\iota_{\sigma^*c} = \iota_c\circ\sigma_*   
$$   
for every $\sigma\in\tcB_m$ and every $c\in\tcC$.    
\end{lemma}   
   
\begin{proof}   
Assertions~(i) and~(ii) follow immediately from the definitions
by checking them on the generators. To prove~(iii) we use~(i) 
with $\psi:=\phi_c(\sigma)$ to obtain   
$$    
\phi_c(\sigma\tau)    
= \phi_c(\sigma)\phi_c(\tau)    
= \phi_{\phi_c(\sigma)_*c}(\tau)\phi_c(\sigma)    
$$    
for $\sigma,\tau\in\tcB_m$ and $c\in\tcC$. Hence 
$$    
(\sigma\tau)^*c  =  \phi_c(\sigma\tau)_*c     
=  \phi_{\phi_c(\sigma)_*c}(\tau)_*\phi_c(\sigma)_*c      
=  \phi_{\sigma^*c}(\tau)_*(\sigma^*c)      
= \tau^*\sigma^*c.    
$$    
That the action is free and transitive follows from Theorem~\ref{thm:tBC}.    
To prove the last equation in~(iii) let $g\in\Gamma_m$.      
Then, by~(i) and~(ii), we have    
$$    
\iota_{\sigma^*c}(g) = \iota_{\phi_c(\sigma)_*c}(g)    
= \phi_c(\sigma)_*\iota_c(g)    
= \iota_c(\sigma_*g).    
$$    
This proves the lemma. 
\end{proof}   
   
\begin{proof}[Proof of Theorem~B]   
Uniqueness is clear. To prove existence,    
fix a marked distinguished configuration $c\in\tcC$,    
let    
$$   
\cS_c:\tcB\to\GL_m(\Z[\Gamma])
$$   
be the cocycle of Theorem~A, and define    
$S:\tcB_m\times\sN_m\to\GL_m(\Z)$ by   
\begin{equation}\label{eq:SsigmaN}   
S(\sigma,N) := \bigl(\rho_N\cdot(\cS_c(\phi_c(\sigma))\circ\iota_c)\bigr)(1)   
\end{equation}   
for $\sigma\in\tcB_m$ and $N\in\sN_m$.
Here $\cS_c(\phi_c(\sigma))\circ\iota_c\in\GL_m(\Z[\Gamma_m])$,   
the representation $\rho_N:\Gamma_m\to\GL_m(\Z)$ is given    
by~\eqref{eq:rhoN}, and the term  
$\rho_N\cdot (\cS_c(\phi_c(\sigma))\circ\iota_c)$    
is understood as the convolution product.    
   
We prove that $S$ satisfies~\eqref{eq:normalize}.    
For ${2\le k\le m}$ we have, by Remark~\ref{rmk:Sc},   
\begin{eqnarray*}   
S(\sigma_k,N)
&= & \bigl(\rho_N(\cS_c(\sigma_{k,c})\circ\iota_c)\bigr)(1) \\   
&=& \rho_N(1)\cS_c(\sigma_{k,c})(1)   
+ \rho_N(g_{k-1})\cS_c(\sigma_{k,c})(g_{k-1,c}^{-1}) \\     
&=&  (\Sigma_k-E_{k,k-1}) + (\one-\eps E_{k-1}N)E_{k,k-1} \\   
&=& \Sigma_k - \eps n_{k-1,k} E_{k-1}.   
\end{eqnarray*}   
Here 
$
\eps:=(-1)^{n(n+1)/2}
$ 
and $E_{k,k-1}$ is the matrix whose    
entries are all zero except for the entry     
$(k,k-1)$ which is equal to $1$. For $1\le i\le m$ we have    
\begin{eqnarray*}   
S(\eps_i,N)   
&= & \bigl(\rho_N(\cS_c(\eps_{i,c})\circ\iota_c)\bigr)(1) \\   
&=& \rho_N(1)\cS_c(\eps_{i,c})(1)    
+ \rho_N(g_i^{-1})\cS_c(\eps_{i,c})(g_{i,c}) \\     
&=&  (\one - E_i) + (\one-(-1)^n \eps E_i N)E_i \\   
&=&  \one - (-1)^n \eps n_{ii} E_i \\   
&=& D_i.   
\end{eqnarray*}   
Here the second equation uses the identity    
$$
\rho_N(g_i^{-1})=\one-(-1)^n\eps E_i N   
$$
and the fourth equation follows from~\eqref{eq:nij-N}.
We also have $(\rho_N\cdot(\cS_c(1)\circ \iota_c))(1)=\rho_N(1)=\one$, and
the {\it (Normalization)} condition~\eqref{eq:normalize} is proved.    

We prove that $S$ satisfies~\eqref{eq:Cocycle} and~\eqref{eq:Rho}   
for all $\sigma,\tau\in\tcB_m$ and $N\in\sN_m$. The proof    
is by induction on the word length of $\sigma$.    
The induction step relies on the following two observations.   
   
\medskip\noindent{\bf Claim~1.}   
{\it If~\eqref{eq:Rho} holds for $\sigma$ and $N$,   
then~\eqref{eq:Cocycle} holds for every $\tau$.}   
   
\medskip\noindent{\bf Claim~2.}   
{\it If~\eqref{eq:Cocycle} holds for $\sigma$,    
$\tau$, and $N$ and~\eqref{eq:Rho} holds for    
the pairs $(\sigma,N)$ and $(\tau,\sigma^*N)$,   
then~\eqref{eq:Rho} also holds for the    
pair $(\sigma\tau,N)$.}   
   
\medskip\noindent   
To carry out the induction argument we first observe    
that~\eqref{eq:Rho} holds for the generators $\sigma_k$   
and $\eps_i$ by direct verification. Hence, by Claim~1,    
\eqref{eq:Cocycle} also holds whenever $\sigma$ is a generator.    
Assume, by induction, that~\eqref{eq:Cocycle} and~\eqref{eq:Rho}   
hold whenever $\sigma$ is a word of length at most $k$.   
Let $\sigma$ be a word of length $k+1$.    
That~\eqref{eq:Rho} holds for every $N$    
follows from Claim~2 by decomposing $\sigma$   
as a product of two words of length at most $k$.    
Hence, by Claim~1, \eqref{eq:Cocycle} holds for every $\tau$.   
Thus it remains to prove Claims~1 and ~2.   
   
To prove Claim~1 it is convenient to abbreviate
$$   
\cM_\tau := \cS_c(\phi_c(\tau))\circ\iota_c\in\GL_m(\Z[\Gamma_m])   
$$   
for $\tau\in\tcB_m$.  Then it follows from Lemma~\ref{le:BmGm}~(ii)   
that the {\it (Cocycle}) condition~\eqref{eq:Coc}   
for $\cS_c$ takes the form   
$$   
\cM_{\sigma\tau} = \cM_\sigma (\sigma_*\cM_\tau) 
$$   
for $\sigma,\tau\in\tcB_m$.  Moreover, equation~\eqref{eq:Rho}    
can be written in the form   
$$   
\sigma_*(\rho_{\sigma^*N})     
= S(\sigma,N)^{-1}\rho_NS(\sigma,N).   
$$   
This implies    
\begin{eqnarray*}    
S(\tau,\sigma^*N)    
&=& \bigl(\rho_{\sigma^*N}\cM_\tau\bigr)(1) \\   
&=& \bigl((\sigma_*\rho_{\sigma^*N})(\sigma_*\cM_\tau)\bigr)(1) \\   
&=& S(\sigma,N)^{-1}    
\bigl(\rho_NS(\sigma,N) (\sigma_*\cM_\tau)\bigr)(1).   
\end{eqnarray*}   
On the other hand we have    
\begin{eqnarray*}    
S(\sigma\tau,N)    
&=&    
\bigl(\rho_N\cM_{\sigma\tau}\bigr)(1) \\   
&=&    
\bigl(\rho_N\cM_\sigma(\sigma_*\cM_\tau)\bigr)(1) \\   
&=&   
\bigl(\rho_NS(\sigma,N) (\sigma_*\cM_\tau)\bigr)(1).   
\end{eqnarray*}   
Here the third equation follows from the definition of the    
convolution product and the fact that $\rho_N$    
is a group homomorphism. This proves Claim~1.   
   
To prove Claim~2, we first observe that equation~\eqref{eq:Cocycle}   
for the triple $(\sigma,\tau,N)$ implies that    
$   
(\sigma\tau)^*N=\tau^*\sigma^*N.    
$   
With this understood Claim~2 follows immediately from the definitions.   
Thus we have proved assertions~(i) and~(iii) of Theorem~B.   
   
To prove~(ii) we recall from Example~\ref{ex:free} that   
any monodromy character ${\cN:\Gamma_m\to\Z^{m\times m}}$   
is uniquely determined by the matrix $N:=\cN(1)$ via   
$$   
\cN = N\rho_N   
$$    
where $\rho_N:\Gamma_m\to\GL_m(\Z)$ is given    
by~\eqref{eq:rhoN}.  This implies that, for every Lefschetz fibration   
$f:X\to\D$ with critical fibers over $Z$, we have    
$$   
\cN^X_c = N^X_c(\rho_{N^X_c}\circ\iota_c^{-1})   
$$   
Hence assertion~(ii) follows from Theorem~\ref{thm:Sc} and the identity   
\begin{equation}\label{eq:equiv}   
\sigma^*N   = \bigl(\sT_{\phi_c(\sigma),c}(\cN)\bigr)(1),   
\qquad \cN := N(\rho_N\circ\iota_c^{-1}),   
\end{equation}   
for $c\in\tcC$ and $\sigma\in\tcB_m$.  
To prove~\eqref{eq:equiv} we observe that
$$   
\rho_N(g)^T N \rho_N(g)=N   
$$
for $g\in\Gamma_m$ (as can be checked on the 
generators $g_1,\dots,g_m$), and that
$$   
S(\sigma,N) = \sum_{g\in\Gamma_m}    
\rho_N(g^{-1}) \cS_c(\phi_c(\sigma))(\iota_c(g))   
$$   
for $\sigma\in\tcB_m$. Hence   
\begin{eqnarray*}    
\sigma^*N    
&=& S(\sigma,N)^TNS(\sigma,N) \\   
&=& \sum_{g,h\in\Gamma_m}    
\cS_c(\phi_c(\sigma))(\iota_c(g))^T \rho_N(g^{-1})^T    
N \rho_N(h^{-1}) \cS_c(\phi_c(\sigma))(\iota_c(h)) \\   
&=& \sum_{g,h\in\Gamma_m} \cS_c(\phi_c(\sigma))(\iota_c(g))^T N    
\rho_N(gh^{-1}) \cS_c(\phi_c(\sigma))(\iota_c(h)) \\   
&=& \sum_{g,h\in\Gamma} \cS_c(\phi_c(\sigma))(g)^T    
\cN(gh^{-1}) \cS_c(\phi_c(\sigma))(h) \\   
&=& \bigl(\cS_c(\phi_c(\sigma))^t\cN \cS_c(\phi_c(\sigma))\bigr)(1) \\   
&=&  \bigl(\sT_{\phi_c(\sigma),c}(\cN)\bigr)(1).   
\end{eqnarray*}    
This proves~\eqref{eq:equiv} and Theorem~B.   
\end{proof}   
 
 
\section{The monodromy groupoid}\label{sec:groupoids} 
 
\begin{figure}[htp]   
\centering    
\includegraphics[scale=0.7]{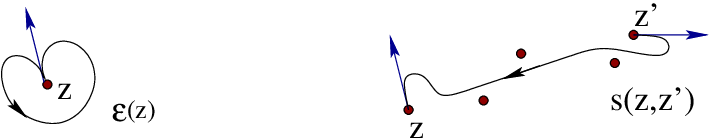}   
\caption{{Generators of the groupoid $\cP$.}}\label{fig:s}   
\end{figure}

Let $Z\subset\mathrm{int}(\D)$ be a finite set and 
$v:Z\to S^1$ be a collection of unit tangent vectors $v(z)\in S^1$ for $z\in Z$.  
The pair $(Z,v)$ is called {\bf admissible} if no three elements lie on  
a straight line and $v(z)\ne \Abs{z'-z}^{-1}(z'-z)$   
for all $z,z'\in Z$ with $z\ne z'$.  
We assume throughout that $(Z,v)$ is an admissible pair.  
Associated to this pair is a groupoid whose objects   
are the elements of $Z$ and the morphisms from   
$z_0$ to $z_1$ are homotopy classes of paths   
$\gamma:(0,1)\to\D\setminus Z$ satisfying    
\begin{equation} \label{eq:P}  
\gamma(t)=z_0+tv(z_0),\qquad    
\gamma(1-t)=z_1+tv(z_1),\qquad    
0\le t\le\delta,   
\end{equation}   
for $\delta>0$ sufficiently small.  (This condition is required to hold   
for a uniform constant $\delta>0$ along a homotopy.)  Composition   
is given by catenation, pushed away from the intermediate point
by the corresponding tangent vector. Denote the set of morphisms 
from $z_0$ to $z_1$ by $\cP(z_0,z_1)$ and write   
$$
\cP := \bigsqcup_{z_0,z_1\in Z}\cP(z_0,z_1).   
$$ 
It is convenient to present the groupoid $\cP$ in terms  
of generators and relations.  The generators are  
$$  
s(z,z')\in\cP(z',z),\qquad \eps(z)\in\cP(z,z)  
$$  
for $z,z'\in Z$ with $z\ne z'$. Geometrically, $\eps(z)$ represents   
a counterclockwise rotation about $z$ and $s(z,z')\in\cP(z',z)$   
represents the straight line from $z'$ to~$z$, modified at each end  
by a counterclockwise turn to match the boundary conditions   
(see Figure~\ref{fig:s}).
To describe the relations we need some definitions.  
An {\bf ordered triangle} is a triple $z_0,w,z_1$   
of pairwise distinct elements of $Z$; it is called   
{\bf local} if its convex hull contains no other   
elements of $Z$.  Associated to every ordered  
triangle is an index  
$$  
\mu(z_0,w,z_1)  
:= \left\{\begin{array}{rl}  
0,&\mbox{if }v(w)\mbox{ points out of the triangle}, \\  
\pm1,&\mbox{if }v(w)\mbox{ points into the triangle}.  
\end{array}\right.  
$$  
Here we choose $+1$ if $z_0,w,z_1$ are ordered   
counterclockwise around the boundary of   
their convex hull, and $-1$ if they are ordered clockwise.  
We have $\mu(z_1,w,z_0)=-\mu(z_0,w,z_1)$ and,  
when $w=e$ is an extremal point of~$Z$, 
i.e.\ when $e$ does not belong to the interior 
of the convex hull of $Z$, we have
\begin{equation}\label{eq:triangle}  
\mu(z_0,e,z_1)+\mu(z_1,e,z_2)=\mu(z_0,e,z_2).  
\end{equation}  
  
\begin{theorem}\label{thm:P}  
Let $(Z,v)$ be an admissible pair.  
Then the groupoid $\cP$ is generated 
by the morphisms $s(z,z')$ and $\eps(z)$
subject to the relations  
\begin{equation}\label{eq:rel1}  
s(z_0,z_1)s(z_1,z_0)=1  
\end{equation}  
and, for every local triangle $z_0,z_1,z_2$,  
\begin{equation}\label{eq:rel2}  
s(z_0,z_1)\eps(z_1)^{\mu_1}  
s(z_1,z_2)\eps(z_2)^{\mu_2}  
s(z_2,z_0)\eps(z_0)^{\mu_0}  
= 1,  
\end{equation}  
where $\mu_1:=\mu(z_0,z_1,z_2)$, $\mu_2:=\mu(z_1,z_2,z_0)$,  
$\mu_0:=\mu(z_2,z_0,z_1)$.   
\end{theorem}  
  
\begin{proof}  
That the generators satisfy these relations follows   
by inspection of local triangles   
(see Figure~\ref{fig:triangle}).
\begin{figure}[htp]   
\centering    
\includegraphics[scale=0.7]{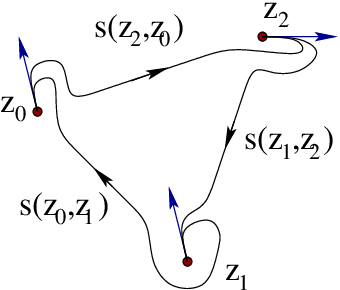}   
\caption{{A positive local triangle.}}\label{fig:triangle}  
\end{figure}  
To prove that there are no further relations and that   
every morphism is a composition of the generators, we choose   
a triangulation of the disc such that each element  
of $Z$ is a vertex and all other vertices are on the boundary.    
Any path $\gamma:(0,1)\to\D\setminus Z$ satisfying~\eqref{eq:P}   
can then be approximated by a smooth path $\gamma'$ intersecting   
the edges transversally and avoiding the vertices.   
Next one can homotop the path to a composition of the  
morphisms associated to the edges of the triangulation  
connecting two vertices in $Z$ and suitable rotations $\eps(z)$.  
To prove that there are no other relations one can study  
the combinatorial pattern of intersection points between   
the path~$\gamma'$ and the edges of the triangulation.   
First, the ambiguity in the choice of the word associated to   
the path~$\gamma'$ is governed by the local triangle   
relation.  Second, every word representing~$\gamma$ can be   
obtained by a suitable choice of~$\gamma'$.  Third, in a generic   
homotopy of~$\gamma'$ there are two kinds of phenomena occuring   
at discrete times.  Either two adjacent intersection points   
are created on the same edge or, conversely, two adjacent  
intersection points cancel. One can then check that these  
phenomena are again governed by the local triangle   
relation.  This completes the sketch of the proof.   
\end{proof}   
  
A {\bf monodromy character on $\cP$}   
is a map $\chi:\cP\to\Z$ satisfying   
\begin{equation}\label{eq:chi0}  
\chi(\gamma^{-1}) = (-1)^n\chi(\gamma),\qquad  
\chi(1) = \left\{\begin{array}{ll}  
0,&\mbox{if }n\mbox{ is odd}, \\  
2(-1)^{n/2},&\mbox{if }n\mbox{ is even},  
\end{array}\right.  
\end{equation}  
for all $\gamma\in\cP$ and   
\begin{equation}\label{eq:chi1}   
\chi(\gamma_{01}\eps_1\gamma_{12}) = \chi(\gamma_{01}\gamma_{12})   
- (-1)^{n(n+1)/2}\chi(\gamma_{01})\chi(\gamma_{12})   
\end{equation}   
for all $\gamma_{01}\in\cP(z_1,z_0)$ and $\gamma_{12}\in\cP(z_2,z_1)$,   
where $\eps_1:= \eps(z_1) \in\cP(z_1,z_1)$ 
denotes the counterclockwise turn about~$z_1$.    
As in Remark~\ref{rmk:N-inverse} one finds that every monodromy   
character satisfies  
\begin{equation}\label{eq:chi2}   
\chi(\gamma_{01}\eps_1^{-1}\gamma_{12})
= \chi(\gamma_{01}\gamma_{12})
- (-1)^n(-1)^{n(n+1)/2}\chi(\gamma_{01})\chi(\gamma_{12}),
\end{equation}   
\begin{equation}\label{eq:chi3}   
\chi(\eps_0\gamma_{01})
= \chi(\gamma_{01}\eps_1)
= (-1)^{n+1}\chi(\gamma_{01}).
\end{equation}   
We refer to the equations~\eqref{eq:chi1} and~\eqref{eq:chi2} 
as the {\bf reflection formulas}.  
The next theorem asserts that a monodromy character is uniquely   
determined by its values on the straight lines.   
Recall from the introduction the definition of $\sN_Z$   
as the set of all matrices ${Q:Z\times Z\to\Z}$   
satisfying~\eqref{eq:Q}.   
   
\begin{theorem}\label{thm:chi-Q}   
For every ${Q\in\sN_Z}$ there is a unique monodromy character   
${\chi^Q:\cP\to \Z}$ such that    
\begin{equation} \label{eq:chiQ}  
\chi^Q(s(z,z'))=Q(z,z')   
\end{equation}  
for all $z,z'\in Z$ with $z\ne z'$.   
\end{theorem}   
   
\begin{remark}\label{rmk:assoc}\rm 
By the reflection formulas~\eqref{eq:chi1} and~\eqref{eq:chi2}, 
the value of $\chi$ on any product of the form  
$$  
\gamma = \gamma_{01}\eps_1^{m_1}\gamma_{12}  
\cdots\eps_{k-1}^{m_{k-1}}\gamma_{k-1,k}  
$$  
is uniquely determined by the values of $\chi$ on the products  
$\gamma_{i,i+1}\cdots\gamma_{j-1,j}$ for $i<j$, regardless in which   
order we apply~\eqref{eq:chi1} and~\eqref{eq:chi2}.   
An example is the identity  
\begin{eqnarray*}  
\chi(\gamma_{01}\eps_1\gamma_{12}\eps_2\gamma_{23})  
&=&   
\chi(\gamma_{01}\gamma_{12}\gamma_{23})    
+  \chi(\gamma_{01})\chi(\gamma_{12})\chi(\gamma_{23}) \\  
&& -\ (-1)^{n(n+1)/2}\chi(\gamma_{01})\chi(\gamma_{12}\gamma_{23}) \\  
&& -\ (-1)^{n(n+1)/2}\chi(\gamma_{01}\gamma_{12})\chi(\gamma_{23}).  
\end{eqnarray*}  
\end{remark}  
   
\begin{remark}\label{rmk:consistent}\rm  
Equations~\eqref{eq:chi0} and~\eqref{eq:chi1} are consistent in the   
following sense.  Suppose the values of $\chi$ on the morphisms  
$\gamma_{01}$, $\gamma_{12}$, $\gamma_{01}\gamma_{12}$  
and their inverses all satisfy the first equation in~\eqref{eq:chi0}   
and that $\chi$ also satisfies the second equation in~\eqref{eq:chi0}  
on each identity morphism.  Suppose further that the values of   
$\chi$ on $\gamma:=\gamma_{01}\eps_1\gamma_{12}$ and its inverse are   
both given by~\eqref{eq:chi1} and~\eqref{eq:chi2}.    
Then $\chi(\gamma^{-1})=(-1)^n\chi(\gamma)$.   
\end{remark}  
  
\begin{proof}[Proof of Theorem~\ref{thm:chi-Q}] 
The proof is by induction on the number $N := \#Z$ of elements of $Z$.   
If~${N=1}$ there is only one monodromy character $\chi$   
given by $\chi(\eps^m)=(-1)^{m(n+1)}\chi(1)$ and so the assertion is obvious.    
In the case $N=2$ with $Z=\{z_0,z_1\}$ every morphism   
involving both vertices can be expressed as a   
composition of morphisms of the form  
$$  
\alpha(k) := s(z_1,z_0) \eps(z_0)^k,\qquad  
\beta(\ell) := s(z_0,z_1)\eps(z_1)^\ell
$$
where $k,l \in \Z$.  
Hence the map $\chi$ is uniquely determined by $\chi(s(z_0,z_1))$.    
Namely, the value of $\chi$ on the product   
$\alpha(k_1)\beta(\ell_1)\alpha(k_2)\beta(\ell_2)\cdots$ is obtained
by  applying equations~\eqref{eq:chi1} and~\eqref{eq:chi2} inductively.  
By Remark~\ref{rmk:assoc}, the answer does not depend on 
the order in which we apply this formula, and hence 
the resulting function $\chi$ satisfies~\eqref{eq:chi1}.   
That it also satisfies~\eqref{eq:chi0} follows from   
Remark~\ref{rmk:consistent}. Thus $\chi$ is a monodromy character.   
  
Now let $N\ge3$ and denote by $E\subset Z$  
the set of extremal points of the convex hull of $Z$.  
Then $E$ has at least three elements, because $N\ge3$ and $Z$   
is admissible.  For every $e\in E$ denote  
$Z_e:=Z\setminus\{e\}$ and let $\cP_e$ be the space   
of all morphisms $\gamma\in\cP$ that can be expressed   
as compositions of generators involving only vertices in $Z_e$.   
Then the induction hypothesis takes the following form.  
  
\medskip\noindent{\bf Induction hypothesis.}  
{\it For every $e\in E$ there is a unique   
monodromy character $\chi_e:\cP_e\to\Z$ that   
satisfies~\eqref{eq:chiQ} for all   
$z,z'\in Z_e$ with $z\ne z'$. Moreover, for all   
$e_1,e_2\in E$, the functions $\chi_{e_1}$   
and $\chi_{e_2}$ agree on $\cP_{e_1}\cap\cP_{e_2}$.}  
  
\medskip\noindent  
Assuming this we wish to prove that there exists a unique monodromy   
character $\chi:\cP\to\Z$ satisfying~\eqref{eq:chiQ}. This   
monodromy character will necessarily restrict to $\chi_e$ 
on $\cP_e$ for every $e\in E$, so that the uniqueness 
of $\chi$ is granted. To prove the existence of $\chi$ 
it is convenient to introduce another category  
with the same set $Z$ of objects, and  whose morphisms   
from $a$ to $z$ are sequences   
$$  
\tga = (z_0,m_0,z_1,m_1,\dots,z_k,m_k)  
$$  
with $m_i\in\Z$, $z_i\in Z$ such that $z_i\ne z_{i+1}$   
for all $i$, and $z_k=a$, $z_0=z$. In this category 
the space of morphisms of length at most $k$   
will be denoted by $\tcP^k(a,z)$ and we write  
$$  
\tcP(a,z) := \bigcup_k\tcP^k(a,z),\qquad  
\tcP := \bigsqcup_{a,z\in Z}\tcP(a,z).  
$$  
  
The notion of a monodromy character extends to the   
category $\tcP$ (with $\gamma^{-1}$ in~\eqref{eq:chi0}   
replaced by the sequence $\tga^\#:=(z_k,-m_k,\dots,z_0,-m_0)$).    
There is a functor $\pi:\tcP\to\cP$ given by   
\begin{equation}\label{eq:phi}  
\pi(\tga) := \eps(z_0)^{m_0}  
s(z_0,z_1)\eps(z_1)^{m_1}\cdots s(z_{k-1},z_k)\eps(z_k)^{m_k}.  
\end{equation}   
By Theorem~\ref{thm:P}, this functor is surjective and two   
morphisms in $\tcP(a,z)$ give the same morphism in $\cP(a,z)$   
if and only if they are related by a sequence of elementary   
moves corresponding to the relations~\eqref{eq:rel1} and~\eqref{eq:rel2}.   
If $\chi:\cP\to\Z$ is a monodromy character then so is the   
composition $\tchi:=\chi\circ\pi:\tcP\to\Z$.  We shall construct  
a monodromy character $\tchi:\tcP\to\Z$ and then show that   
it descends to $\cP$.  The proof has three steps.  
  
\medskip\noindent{\bf Step~1.}  
{\it If $e$ is an extremal point and $a,z\in Z_e$   
then $s(z,e)\eps(e)^\mu s(e,a)\in\cP_e$ where  
$\mu:=\mu(z,e,a)\in\{-1,0,1\}$.}  
  
\medskip\noindent  
When $z,e,a$ is a local triangle then it follows from   
the local triangle relation~\eqref{eq:rel2} that,   
for suitable integers $\mu_a,\mu_z\in\{-1,0,1\}$,  
we have  
$$ 
s(z,e)\eps(e)^\mu s(e,a)   
= \eps(z)^{\mu_z}s(z,a)\eps(a)^{\mu_a} \in\cP_e(a,z).  
$$  
In general we can find a sequence $z=z_0,\dots,z_k=a$ in $Z_e$  
such that $z_{i-1},e,z_i$ is a local triangle and hence  
$$  
\gamma_{i-1,i}:=s(z_{i-1},e)\eps(e)^{\mu_i}s(e,z_i)\in\cP_e(z_i,z_{i-1})  
$$  
for $\mu_i:=\mu(z_{i-1},e,z_i)$.   
Composing these morphisms we obtain  
$$  
s(z,e)\eps(e)^{\mu'} s(e,a)   
= \gamma_{01}\gamma_{12}\cdots\gamma_{k-1,k}\in\cP_e(a,z),\qquad  
\mu':=\mu_1+\cdots+\mu_k.  
$$ 
Hence the assertion of Step~1 follows from~\eqref{eq:triangle} 
which implies that $\mu=\mu'$. 
 
\medskip\noindent{\bf Step~2.}   
{\it There is a unique monodromy character   
$\tchi:\tcP\to\Z$ satisfying the following conditions.  
     
\smallskip\noindent{\bf (i)}   
If $z_0\ne z_1$ and $\ts(z_0,z_1):=(z_0,0,z_1,0)$ then   
$\tchi(\ts(z_0,z_1)) = Q(z_0,z_1)$.  
    
\smallskip\noindent{\bf (ii)}   
If $a,z\in Z$, $e\in E\setminus\{a,z\}$, and   
$$  
\tga = (z_0,m_0,\dots,z_k,m_k)\in\tcP(a,z),\qquad  
z_k=a,\qquad z_0=z,  
$$  
such that   
\begin{equation}\label{eq:mu}  
z_i=e\quad\implies\quad   
m_i=\mu_i:=  
\left\{\begin{array}{ll}  
\mu(z_{i-1},z_i,z_{i+1}),&\mbox{if }z_{i-1}\ne z_{i+1}, \\  
0,&\mbox{if }z_{i-1}=z_{i+1},  
\end{array}\right.  
\end{equation}   
then $\tchi(\tga):=\chi_e(\pi(\tga))$.  
(In this case $\pi(\tga)\in\cP_e$ by Step~1.) }  
    
\medskip\noindent  
We prove uniqueness of the value $\tchi(\tga)$  
by induction on the length $k$. For $k=0$ we must have  
$$  
\tchi(z_0,m_0) =   
(-1)^{m_0(n+1)}\left\{\begin{array}{ll}  
0,&\mbox{if }n\mbox{ is odd}, \\  
2(-1)^{n/2},&\mbox{if }n\mbox{ is even},  
\end{array}\right.  
$$  
and for $k=1$ we can argue as in the case $N=2$ above.  
Now let $k\ge 2$ and suppose, by induction, that uniqueness  
of $\tchi(\tal)$ has been verified for every morphism of length   
at most $k-1$. Fix two elements $a,z\in Z$, an extremal point   
$e\in E\setminus\{a,z\}$, and a morphism  
$$  
\tga = (z_0,m_0,\dots,z_k,m_k)\in\tcP(a,z),\qquad  
z_k=a,\qquad z_0=z.  
$$  
If $\tga$ and $e$ satisfy~\eqref{eq:mu} then $\tchi(\tga)$   
is determined by condition~(ii). If there is precisely one value   
$i$ with $z_i=e$ and $m_i\ne\mu_i$ then the value $\tchi(\tga)$   
is determined inductively by~(ii) and equation~\eqref{eq:chi1}  
for $\tchi$. Namely, with   
$$  
\tal:= (z_0,m_0,\dots,z_i,m_i),\qquad  
\tbe:=(z_i,0,z_{i+1},m_{i+1},\dots,z_k,m_k),  
$$ 
and $\teps_i:=(z_i,1)$ we have $\tga=\tal\tbe$ and
$$  
\tga':=\tal\teps_i\tbe  
=(z_0,m_0,\dots,z_i,m_i+1,\dots,z_k,m_k),  
$$  
\begin{equation}\label{eq:chi}  
\tchi(\tga') = \tchi(\tga)-(-1)^{n(n+1)/2}\tchi(\tal)\tchi(\tbe).  
\end{equation}  
This proves uniqueness when~\eqref{eq:mu} fails for   
precisely one value of $i$. In general one can repeat   
this argument inductively for all $i$.   
  
The same argument can be used to construct the value  
$\tchi(\tga)$ by induction on the length $k$.   
The induction hypothesis is that $\tchi^{k-1}:\tcP^{k-1}\to\Z$  
has been constructed as to satisfy~(i) and~(ii) as well   
as~\eqref{eq:chi0} and~\eqref{eq:chi1} (the latter for   
compositions of two morphisms with total length at most $k-1$).  
We must extend this map to one on $\tcP^k$ with the same   
properties. For this we fix a morphism $\tga\in\tcP(a,z)$ of   
length $k$ and use the above induction argument with the auxiliary  
choice of an extremal point $e$ to define $\tchi^k(\tga)$.   
Assume there are two such extremal points   
$e,e'\in E\setminus\{a,z\}$. If   
\begin{equation}\label{eq:muE}  
z_i\in\{e,e'\}\qquad\implies\qquad m_i=\mu_i  
\end{equation}  
then $\pi(\tga)\in\cP_e\cap\cP_{e'}$ and our value $\tchi(\tga)$   
is independent of $e$ by~(ii).  In general, the same induction argument   
as above, using~\eqref{eq:chi1} repeatedly for those values  
of $i$ for which~\eqref{eq:muE} fails, shows that the two definitions  
of $\tchi(\tga)$ (associated to $e$ and $e'$) agree for all values  
of the multiplicities $m_i$.  Thus our definition of $\tchi^k(\tga)$  
is independent of $e$.   
  
It remains to show that $\tchi^k$ satisfies the   
conditions~\eqref{eq:chi0} and~\eqref{eq:chi1}    
(for compositions of two morphisms with total   
length equal to $k$). Let   
$$  
\tal\in\tcP^i(a,b),\qquad \tbe\in\tcP^{k-i}(b,c)  
$$   
be two such morphisms.  We must prove that~\eqref{eq:chi} holds.  
If there is an extremal point $e\in E\setminus\{a,b,c\}$  
we can argue as above and verify~\eqref{eq:chi} for this  
pair by using Remark~\ref{rmk:assoc} and reducing the discussion  
to the case where both $\tal$ and $\tbe$ satisfy~\eqref{eq:mu}.   
In this case equation~\eqref{eq:chi} follows from the fact   
that $\chi_e$ is a monodromy character. If   
$  
E=\{a,b,c\}  
$  
we may choose $e=b$ and use the definition of   
$\tchi(\tal\tbe)$ and $\tchi(\tal\teps(b)\tbe)$ given above 
to establish the reflection formula. This shows that $\tchi^k$ 
satisfies~\eqref{eq:chi1}.  That it also satisfies~\eqref{eq:chi0}   
follows from Remark~\ref{rmk:consistent}. This proves Step~2.   
  
\medskip\noindent{\bf Step~3.}  
{\it Let $\tchi:\tcP\to\Z$ be as in Step~2.     
Then there is a unique monodromy character   
$\chi:\cP\to\Z$ such that $\tchi=\chi\circ\pi$.}  
  
\medskip\noindent  
We must prove that $\tchi$ descends to $\cP$.  
This means that for all $a,z\in Z$ and all   
$\tga,\tga'\in\tcP(a,z)$ we have  
\begin{equation}\label{eq:descend}  
\pi(\tga)=\pi(\tga')\qquad\implies\qquad  
\tchi(\tga)=\tchi(\tga').  
\end{equation}  
To see this fix an extremal point $e\in Z\setminus\{a,z\}$.  
If both $\tga$ and $\tga'$ satisfy condition~\eqref{eq:mu}  
then~\eqref{eq:descend} follows immediately from~(ii) in Step~2.   
In general, we must lower or raise the indices $m_i$ with   
$z_i=e$.  Inserting a term $s(z_i,w)s(w,z_i)$ into the word  
associated to $\tga$ does not affect this induction.   
When $w=e$ this is obvious and when $z_i=e$ the total number   
of induction steps on the left and right of the inserted term   
agrees with the number of steps at $z_i$ before inserting, because  
$$  
\mu(z_{i-1},e,z_{i+1}) = \mu(z_{i-1},e,w)+\mu(w,e,z_{i+1}).   
$$   
Replacing a term $s(z_i,z_{i+1})$ with the product  
$$  
s(z_i,z_{i+1})  
= \eps(z_i)^{\mu(z_{i+1},z_i,w)}  
s(z_i,w)\eps^{\mu(z_i,w,z_{i+1})}   
s(w,z_{i+1})\eps(z_{i+1})^{\mu(w,z_{i+1},z_i)}  
$$  
for a local triangle $z_i,w,z_{i+1}$ also does not affect this   
induction for a similar reason.  Here at most one of the   
terms $z_i,w,z_{i+1}$ can be equal to $e$.  If $w=e$ then   
the exponent $\mu(z_i,w,z_{i+1})$ already satisfies~\eqref{eq:mu}.  
If $z_i=e$ we use the fact that  
$\mu(z_{i-1},e,z_{i+1})+\mu(z_{i+1},e,w)=\mu(z_{i-1},e,w)$  
and similarly in the case $z_{i+1}=e$.    
  
Thus we have proved that $\tchi$ descends to $\cP$.    
That the resulting map ${\chi:\cP\to\Z}$ is a monodromy character   
is obvious and that it satisfies~\eqref{eq:chiQ} follows from the   
fact that $\pi(\ts(z_0,z_1))=s(z_0,z_1)$ for all $z_0,z_1\in Z$   
with $z_0\ne z_1$.  Uniqueness follows immediately from the   
uniqueness statement in Step~2. This completes the proof   
of Theorem~\ref{thm:chi-Q}.  
\end{proof}  
  

\section{Proof of Theorem~C}\label{sec:proofC}      

Fix an admissible pair $(Z,v)$ and let $\cP$ be the groupoid  
of Section~\ref{sec:groupoids}. Recall from Section~\ref{sec:bases} 
the notion of marked distinguished configurations
$c=(c_1,\dots,c_m)\in\tcC$ with $m=\#Z$ such that 
$\{c_1(1),\dots,c_m(1)\}=Z$ and $\dot c_i(1)=-v(c_i(1))$ 
for $i=1,\dots,m$.

\begin{lemma}\label{le:chi-N}   
{\bf (i)}   
For every $c\in\tcC$ and every $N=(n_{ij})_{i,j=1}^m\in\sN_m$    
there is a unique monodromy character $\chi_{c,N}:\cP\to\Z$   
such that   
$$   
\chi_{c,N}(c_i\cdot c_j^{-1}) = n_{ij}   
$$   
for all $i$ and $j$.   
   
\smallskip\noindent{\bf (ii)}   
For all $c\in\tcC$, $N\in\sN_m$, and $\sigma\in\tcB_m$ we have   
$$   
\chi_{c,N}=\chi_{\sigma^*c,\sigma^*N}.   
$$   
\end{lemma}   
   
\begin{proof}   
We prove~(i). 
Let $z_1,\dots,z_m$ be the ordering of $Z$ determined by $c$. 
Uniqueness follows from the definition of a monodromy   
character and from the fact that any element $\gamma\in\cP(z_i,z_j)$   
can be expressed uniquely as $c_i \cdot g \cdot c_j^{-1}$ for some   
$g\in \Gamma_m\cong \Gamma$, and $g$ can be expressed uniquely in    
reduced form as a product of the generators $g_1,\dots,g_m$ and their   
inverses. To prove existence, let $c\in\tcC$ and   
$N=(n_{ij})_{i,j=1}^m\in\sN_m$    
be given and define   
\begin{equation}\label{eq:chicN}   
\chi_{c,N}(\gamma)    
= \bigl(N\rho_N(\iota_c^{-1}(c_i^{-1}\cdot\gamma\cdot c_j))\bigr)_{ij}   
\end{equation}   
for $\gamma\in\cP(z_i,z_j)$. Denote    
$n_{ij}(g):= \bigl(N\rho_N(\iota_c^{-1}(g))\bigr)_{ij}$   
for $g\in\Gamma$.  These functions define a monodromy character   
on $\Gamma$ in the sense of Definition~\ref{def:Gamma}    
(see Example~\ref{ex:free}).    
Let $\gamma_{ij}\in\cP(z_j,z_i)$ and $\gamma_{jk}\in\cP(z_k,z_j)$   
and denote   
$$   
g:= c_i^{-1}\cdot\gamma_{ij}\cdot c_j,\qquad   
h:= c_j^{-1}\cdot\gamma_{jk}\cdot c_k.   
$$   
Abbreviate $\chi:=\chi_{c,N}$.  Then    
\begin{eqnarray*}   
\chi(\gamma_{ij}\eps_j\gamma_{jk})   
&=&     
n_{ik}(c_i^{-1}\cdot\gamma_{ij}\eps_j\gamma_{jk}\cdot c_k) \\   
&=&     
n_{ik}((c_i^{-1}\cdot\gamma_{ij}\cdot c_j)   
\cdot(c_j^{-1}\cdot\eps_j\cdot c_j)   
\cdot(c_j^{-1}\cdot\gamma_{jk}\cdot c_k)) \\   
&=&     
n_{ik}(gg_jh) \\   
&=&     
n_{ik}(gh) - (-1)^{n(n+1)/2}n_{ij}(g)n_{jk}(h)   \\   
&=&     
\chi(\gamma_{ij}\gamma_{jk}) - (-1)^{n(n+1)/2}   
\chi(\gamma_{ij})\chi(\gamma_{jk}).   
\end{eqnarray*}   
Hence $\chi$ is a monodromy character.   
Moreover, by~\eqref{eq:chicN} we have   
$$   
\chi(c_i\cdot c_j^{-1})    
=  \bigl(N\rho_N(\iota_c^{-1}(1))\bigr)_{ij}    
= n_{ij}   
$$   
for all $i$ and $j$. This proves~(i).   

We prove~(ii). Abbreviate $\tau:=\phi_c(\sigma)$ and   
$$   
n_{i'j'}(g):=\bigl(N\rho_N(\iota_c^{-1}(g))\bigr)_{i'j'},\qquad   
\tn_{ij}(g):=\bigl(   
(\sigma^*N)\rho_{\sigma^*N}(\iota_{\sigma^*c}^{-1}(g))   
\bigr)_{ij}.   
$$   
Let $i,j\in\{1,\dots,m\}$ and $i'=\pi_{\tau,c}(i)$, $j':=\pi_{\tau,c}(j)$.
Let $\gamma\in\cP(z_{j'},z_{i'})$.  Then    
\begin{eqnarray*}   
\chi_{\sigma^*c,\sigma^*N}(\gamma)   
&=&    
\tn_{ij}((\sigma^*c_i)^{-1}\cdot\gamma\cdot(\sigma^*c_j)) \\   
&=&    
\tn_{ij}((\tau_*c_i)^{-1}\cdot c_{i'}\cdot c_{i'}^{-1}   
\cdot\gamma\cdot c_{j'}\cdot c_{j'}^{-1}\cdot(\tau_*c_j)) \\   
&=&    
\tn_{ij}(s_{i,c}(\tau)^{-1}\cdot c_{i'}^{-1}   
\cdot\gamma\cdot c_{j'}\cdot s_{j,c}(\tau)) \\   
&=&    
n_{i'j'}(c_{i'}^{-1}\cdot\gamma\cdot c_{j'}) \\   
&=&   
\chi_{c,N}(\gamma).   
\end{eqnarray*}   
Here the third equation follows from~\eqref{eq:sjc}   
and the fourth equation follows from~\eqref{eq:tn}.     
This proves the lemma.   
\end{proof}    
   
\begin{proof}[Proof of Theorem~C]   
The invariance of the map $(c,N)\mapsto Q_{c,N}$ follows from   
Lemma~\ref{le:chi-N}~(ii) and the fact that    
$$   
Q_{c,N}(z,z') =\chi_{c,N}(s(z,z')).   
$$    
Given $Q\in\sN_Z$, the existence of an equivariant map   
$c\mapsto N_c$ with $Q_{c,N_c}=Q$ follows from    
Lemma~\ref{le:chi-N} and Theorem~\ref{thm:chi-Q}.   
We define   
$$   
N_c:=(n_{ij})_{i,j=1}^m,\qquad    
n_{ij} := \chi^Q(c_i\cdot c_j^{-1}).   
$$   
Then it follows from uniqueness in Lemma~\ref{le:chi-N}    
that $\chi_{c,N_c}=\chi^Q$ and hence   
$$   
Q_{c,N_c}(z_i,z_j) = \chi_{c,N_c}(s(z_i,z_j))   
= \chi^Q(s(z_i,z_j)) = Q(z_i,z_j).   
$$   
Here the first equation follows from the definition of   
$Q_{c,N_c}$ in the introduction, and the last equation
follows from the definition of $\chi^Q$ in Theorem~\ref{thm:chi-Q}.   
This proves existence. To prove uniqueness, let $c\mapsto  
N_c'=(n_{ij}')_{i,j=1}^m$ be another   
equivariant family satisfying $Q_{c,N_c'}=Q$, so that  
$$  
\chi_{c,N_c'}(s(z,z'))   
= Q_{c,N_c'}(z,z')  
= Q_{c,N_c}(z,z')  
= \chi_{c,N_c}(s(z,z')).  
$$  
By uniqueness in Theorem~\ref{thm:chi-Q} we have  
$\chi_{c,N_c'}=\chi_{c,N_c}=\chi^Q$, and therefore  
$$  
n'_{ij} = \chi_{c,N_c'}(c_i\cdot c_j^{-1})   
= \chi_{c,N_c}(c_i\cdot c_j^{-1}) = n_{ij}  
$$   
for all $i,j$. This completes the proof of Theorem~C.    
\end{proof}

      
\section{An Example}\label{sec:ex}      
    
\begin{figure}[htp]      
\centering       
\includegraphics[scale=0.6]{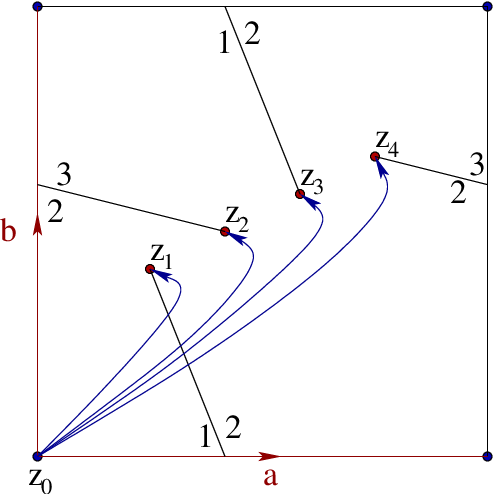}      
\caption{{A $3$-sheeted branched cover over the torus.}}\label{fig:5}      
\end{figure}

Let $f:\Sigma\to\T^2$ be a degree $3$ branched cover of a genus $3$      
surface over the torus with $4$ branch points of branching order $2$       
(see Figure~\ref{fig:5}).      
The map $f$ is a covering over the complement of the two       
slits from $z_1$ to $z_3$ and $z_2$ to $z_4$ depicted in the figure.       
The surface $\Sigma$ is obtained by gluing the three sheets       
labeled $1$, $2$ and $3$ along the slits as indicated.       
We identify the fiber $M$ over $z_0$ with the set $\{1,2,3\}$.       
Each vanishing cycle is an oriented $0$-sphere, i.e$.$       
consists of two points labeled by opposite signs       
(which can be chosen arbitrarily). The marked distinguished    
configuration indicated in the figure determines the vanishing     
cycles up to orientation. We choose the orientations as follows:      
$$      
L_1=\{1^-,2^+\}, \qquad L_2=L_4=\{1^+,3^-\}, \qquad L_3=\{1^+,2^-\}.      
$$      
The fundamental group $\Gamma$ of the $4$-punctured torus is     
generated by the elments $a,b,g_1,g_2,g_3,g_4$ where the $g_i$     
are determined by the distinguished configuration and $a$     
and $b$ are the horizontal and vertical loops     
indicated in Figure~\ref{fig:5}.  These generators     
are subject to the relation $b^{-1}a^{-1}ba=g_1g_2g_3g_4$.    
Understood as permutations, the generators $g_1$, $g_3$       
and $a$ act by the transposition $(1,2)$,       
while $g_2$ and $g_4$ act by the transposition $(1,3)$,       
and $b$ acts by the transposition $(2,3)$. Note that       
the commutator    
$$      
b^{-1}a^{-1}ba=g_1g_2g_3g_4=(231)      
$$      
is not the identity.       
    
The intersection matrix of the $L_i$ is       
$$      
\cN_c(1)=\left(\begin{array}{rrrr}       
2 & -1 & -2 & -1 \\      
-1 & 2 & 1 & 2 \\      
-2 & 1 & 2 & 1 \\      
-1 & 2 & 1 & 2      
\end{array}      
\right).      
$$     
Moreover, the images of the $L_i$ under $a$ are      
$$      
a(L_1)=\{1^+,2^-\}, \qquad      
a(L_2)=a(L_4)=\{2^+,3^-\}, \qquad      
a(L_3)=\{1^-,2^+\}    
$$      
and hence      
$$      
\cN_c(a)=\left(\begin{array}{rrrr}     
-2 & 1 & 2 & 1 \\      
1 & 1 & -1 & 1 \\      
2 & -1 & -2 & -1 \\      
1 & 1 & -1 & 1      
\end{array}      
\right).    
$$     
Similarly,      
$$      
\cN_c(g_2)=\left(\begin{array}{rrrr}       
1 & 1 & -1 & 1 \\      
1 & -2 & -1 & -2 \\      
-1 & -1 & 1 & -1 \\      
1 & -2 & -1 & -2      
\end{array}      
\right), \       
\cN_c(b)=\left(\begin{array}{rrrr}       
1 & -2 & -1 & -2 \\      
-2 & 1 & 2 & 1 \\      
-1 & 2 & 1 & 2 \\      
-2 & 1 & 2 & 1      
\end{array}      
\right),      
$$      
$$      
\cN_c(ab)=\left(\begin{array}{rrrr}       
-1 & 2 & 1 & 2 \\      
-1 & -1 & 1 & -1 \\      
1 & -2 & -1 & -2 \\      
-1 & -1 & 1 & -1      
\end{array}      
\right), \      
\cN_c(ba)=\left(\begin{array}{rrrr}       
-1 & -1 & 1 & -1 \\      
2 & -1 & -2 & -1 \\      
1 & 1 & -1 & 1 \\      
2 & -1 & -2 & -1      
\end{array}      
\right).      
$$      
Note that $\cN_c(g_1)=\cN_c(g_3)=\cN_c(a)$     
and $\cN_c(g_4)=\cN_c(g_2)$.      
The matrix $\cN_c(g_i)$ can also be obtained from $\cN_c(1)$     
via the gluing formula      
$$      
\cN_c(g_i)=\cN_c(1) - \cN_c(1) E_i \cN_c(1).      
$$    
The action of the fundamental group on the homology of the       
fiber $M=\{1,2,3\}$ obviously factors through the permutation       
group $\mathfrak{S}_3$ and hence the function      
$\cN_c:\Gamma\to \Z^{4\times 4}$ is uniquely determined by its values on      
$1$, $g_2$, $a$, $b$, $ab$ and $ba$. One can verify directly that $\cN_c$       
satisfies~(\ref{eq:N1}-\ref{eq:N3}), for example $\cN_c(ba)=\cN_c(ab)^T$    
and       
$$      
\cN_c(b)=\cN_c(ag_1b)=\cN_c(ab)-\cN_c(a)E_1\cN_c(b).      
$$     
   
   
\appendix   
   
\section{Lefschetz fibrations}\label{app:L} 
 
In this section we review some basic facts about 
Lefschetz fibrations (see~\cite[Chapter~I]{AGV}). 
Let $X$ be a compact K\"ahler manifold of complex dimension $n+1$ 
and ${f:X\to\Sigma}$ be a Lefschetz fibration with a finite set 
$Z\subset\Sigma$ of critical values. Let $m:=\#Z$ be the number 
of critical values. We assume that each critical fiber $f^{-1}(z)$ 
contains precisely one (Morse) critical point $x_z$. 
Choose a regular value $z_0\in\Sigma$ and denote $M := f^{-1}(z_0)$. 
Then the fundamental group 
$$ 
\Gamma := \pi_1(\Sigma\setminus Z,z_0) 
$$ 
acts on the middle dimensional homology 
$$ 
H_n(M):=H_n(M;\Z)/\mathrm{torsion} 
$$ 
via parallel transport. 
We denote the action by
$$
\rho:\Gamma\to\Aut(H_n(M)).
$$
Each path $c:[0,1]\to\Sigma$ connecting $z_0=c(0)$ to a critical 
value $z=c(1)$, avoiding critical values for $t<1$, 
and satisfying $\dot c(1)\ne 0$, 
determines a vanishing cycle $L\subset M$ and an element 
$g\in\Gamma$.  Geometrically, $L$ is the set of all points 
in $M$ that converge to the critical point $x_z\in f^{-1}(z)$
under parallel transport along~$c$. The orientation of $L$ is 
not determined  by the path and can be chosen independently. 
The element $g$ is the homotopy class of the path obtained 
by following $c$, encircling $z$ once counterclockwise, 
and then following $c^{-1}$. Thus $\rho(g)$ acts on $H_n(M)$ 
by the Dehn--Arnold--Seidel twist $\psi_L$ about $L$ with its framing 
as a vanishing cycle. (The framing is a choice of a 
diffeomorphism from $L$ to the $n$-sphere, see Seidel~\cite{S1}.) 
Now choose $m$ such paths $c_i$, $i=1,\dots,m$, 
connecting $z_0$ to the critical values, 
denote by $L_i\in H_n(M)$ the homology classes of the 
resulting vanishing cycles with fixed choices of orientations, 
and by $g_i\in\Gamma$ the resulting elements of the 
fundamental group. 
    
\begin{remark}\label{rmk:vc}\rm   
For $i=1,\dots,m$ the action of $g_i$ on $H_n(M)$    
is given by     
\begin{equation}\label{eq:dehn}    
\rho(g_i)\alpha = \alpha - (-1)^{n(n+1)/2}\inner{L_i}{\alpha}L_i,
\end{equation}
where $\inner{\cdot}{\cdot}$ denotes the intersection form.
The intersection numbers
$$
n_{ij}(g):=\inner{L_i}{\rho(g)L_j}, \qquad i,j=1,\dots,m,
$$
satisfy~(\ref{eq:N1}-\ref{eq:N3}).        
Equation~\eqref{eq:N1} follows from the (skew-)symmetry    
of the intersection form, equation~\eqref{eq:N2} is a general fact 
about the self-intersection numbers of Lagrangian spheres, 
and~\eqref{eq:N3} follows from~\eqref{eq:dehn}. 
Equation~\eqref{eq:dehn} follows from Example~\ref{ex:dehn} below. 
\end{remark} 
   
\begin{remark}\label{rmk:horizontal}\rm    
The intersection number    
$    
n_{ij}(g) = \inner{L_i}{\rho(g)L_j}    
$    
can be interpreted as the algebraic number of horizontal lifts of the      
path $c_jg^{-1}c_i^{-1}$ (first $c_i^{-1}$ then $g^{-1}$ then $c_j$)    
connecting the critical points $x_i$ and~$x_j$.  That these numbers     
continue to be meaningful (at least conjecturally)     
in certain infinite dimensional    
settings (where the vanishing cycles only exist in some heuristic sense)    
is one of the key ideas in Donaldson--Thomas theory~\cite{DT}.    
\end{remark} 
    
\begin{example}\label{ex:dehn}\rm
The archetypal example of a Lefschetz fibration is the    
function $f:\C^{n+1}\to\C$ given by    
$$    
f(z_0,\dots,z_n) := z_0^2 + z_1^2 + \cdots + z_n^2.      
$$    
Denote the fiber over $1$ by       
$$      
M:=f^{-1}(1) = \left\{x+iy\in\C^{n+1}\,|\,      
\Abs{x}^2-\Abs{y}^2=1,\,\inner{x}{y}=0\right\}.      
$$      
The vanishing cycle obtained by parallel transport along the real axis    
from $0$ to $1$ is the unit sphere $L:=M\cap\R^{n+1}$.     
The manifold $M$ is symplectomorphic to the cotangent      
bundle of the $n$-sphere      
$$      
T^*S^n = \left\{(\xi,\eta)\in\R^{n+1}\times\R^{n+1}\,|\,      
\Abs{\xi}=1,\,\inner{\xi}{\eta}=0\right\}      
$$      
via      
$      
M\to T^*S^n:(x,y)\mapsto(x/\Abs{x},\Abs{x}y).      
$      
The monodromy around the loop      
$t\mapsto e^{2\pi it}$ is the symplectomorphism      
$\psi:T^*S^n\to T^*S^n$ given by      
\begin{equation}\label{eq:psi}      
\psi(\xi,\eta)      
=\left(-\cos(\theta)\xi-\sin(\theta)\frac{\eta}{\Abs{\eta}},      
\sin(\theta)\Abs{\eta}\xi-\cos(\theta)\eta\right)      
\end{equation}      
where $\theta := 2\pi|\eta|/\sqrt{1+4|\eta|^2}$.       
Let $\bar M$ be the sphere bundle over $S^n$ obtained       
by the one point compactification of each fiber of the       
cotangent bundle.  Its middle dimensional homology       
$H_n(\bar M;\Z)$ has two generators, namely the       
homology class of the zero section and the homology      
class of the fiber, and $\psi$ acts on $H_n(\bar M;\Z)$      
by      
$$      
\psi_*\alpha     
= \alpha - (-1)^{n(n+1)/2}\inner{L}{\alpha}L.      
$$      
To see this fix $\xi$ in~(\ref{eq:psi}) and project to $S^n$       
to get a map of degree $(-1)^{n+1}$. An additional factor       
$(-1)^{n(n-1)/2}$ arises by comparing the orientation of the       
cotangent bundle $T^*L$ with the complex orientation of $M$.      
In the even case the value of the factor follows from the       
observation that the self-intersection number of the       
vanishing cycle is  $\inner{L}{L}=2(-1)^{n/2}$.    
\end{example}    
       

\bigskip
    
\noindent Authors' addresses:
    
\medskip 
    
\noindent G.~Massuyeau, Institut de Recherche Math\'ematique Avanc\'ee (IRMA), 
UMR 7501, Universit\'e de Strasbourg \& CNRS, 7 rue Ren\'e Descartes, 
67084 Strasbourg, France. {\tt massuyea@math.unistra.fr}    
    
\medskip
    
\noindent A.~Oancea, Institut de Math\'ematiques de Jussieu (IMJ-PRG), 
UMR 7586, Universit\'e Pierre et Marie Curie \& CNRS, 4 place Jussieu, 
75005 Paris, France. {\tt oancea@math.jussieu.fr}  
    
\medskip 

\noindent D.~Salamon, Department of Mathematics, ETH, 
R\"amistrasse 101, 8092 Zurich, Switzerland. {\tt salamon@math.ethz.ch}
    
\end{document}